# Exponential Growth Series and Benford's Law

## ABSTRACT

Exponential growth occurs when the growth rate of a given quantity is proportional to the quantity's current value. Surprisingly, when exponential growth data is plotted as a simple histogram disregarding the time dimension, a remarkable fit to the positively skewed k/x distribution is found, where the small is numerous and the big is rare, and with a long tail falling to the right of the histogram. Such quantitative preference for the small has a corresponding digital preference known as Benford's Law which predicts that the first significant digit on the left-most side of numbers in typical real-life data is proportioned between all possible 1 to 9 digits approximately as in LOG(1 + 1/digit), so that low digits occur much more frequently than high digits in the first place. Exponential growth series are nearly perfectly Benford given that plenty of elements are considered and that (for low growth) order of magnitude is approximately an integral value. An additional constraint is that the logarithm of the growth factor must be an irrational number. Since the irrationals vastly outnumber the rationals, on the face of it, this constraint seems to constitute the explanation of why almost all growth series are Benford, yet, in reality this is all too simplistic, and the real and more complex explanation is provided in this article. Empirical examinations of close to a half a million growth series via computerized programs almost perfectly match the prediction of the theoretical study on rational versus irrational occurrences, thus in a sense confirming both, the empirical work as well as the theoretical study. In addition, a rigorous mathematical proof is provided in the continuous growth case showing that it exactly obeys Benford's Law. A non-rigorous proof is given in the discrete case via uniformity of mantissa argument. Finally cases of discrete series embedded within continuous series are studied, detailing the degree of deviation from the ideal Benford configuration.

Alex Ely Kossovsky
akossovsky@gmail.com

# PART 1:
# BENFORD'S LAW

# PART 2:
# EXPONENTIAL GROWTH SERIES

# PART 1:
# BENFORD'S LAW

# [I]   The First Digit on the Left Side of Numbers

It has been discovered that the first digit on the left-most side of numbers in real-life data sets is most commonly of low value such as {1, 2, 3} and rarely of high value such as {7, 8, 9}.
As an example serving as a brief and informal empirical test, a small sample of 40 values relating to geological data on time between earthquakes is randomly chosen from the data set on all global earthquake occurrences in 2012 – in units of seconds. Figure A depicts this small sample of 40 numbers. Figure B emphasizes in bold and black color the 1st digits of these 40 numbers.

| | | | |
|---|---|---|---|
| 285.29 | 185.35 | 2579.80 | 27.11 |
| 5330.22 | 1504.49 | 1764.41 | 574.46 |
| 1722.16 | 815.06 | 3686.84 | 1501.61 |
| 494.17 | 362.48 | 1388.13 | 1817.27 |
| 3516.80 | 5049.66 | 2414.06 | 387.78 |
| 4385.23 | 2443.98 | 2204.12 | 1224.42 |
| 1965.46 | 3.61 | 1347.30 | 271.23 |
| 3247.99 | 753.80 | 1781.45 | 593.59 |
| 1482.64 | 1165.04 | 4647.39 | 1219.19 |
| 251.12 | 7345.52 | 1368.79 | 4112.13 |

**Figure A**: Sample of 40 Time Intervals between Earthquakes

| | | | |
|---|---|---|---|
| **2**85.29 | **1**85.35 | **2**579.80 | **2**7.11 |
| **5**330.22 | **1**504.49 | **1**764.41 | **5**74.46 |
| **1**722.16 | **8**15.06 | **3**686.84 | **1**501.61 |
| **4**94.17 | **3**62.48 | **1**388.13 | **1**817.27 |
| **3**516.80 | **5**049.66 | **2**414.06 | **3**87.78 |
| **4**385.23 | **2**443.98 | **2**204.12 | **1**224.42 |
| **1**965.46 | **3**.61 | **1**347.30 | **2**71.23 |
| **3**247.99 | **7**53.80 | **1**781.45 | **5**93.59 |
| **1**482.64 | **1**165.04 | **4**647.39 | **1**219.19 |
| **2**51.12 | **7**345.52 | **1**368.79 | **4**112.13 |

**Figure B**: The First Digits of the Earthquake Sample

Clearly, for this very small sample, low digits occur by far more frequently on the first position than do high digits. A summary of the digital configuration of the sample is given as follows:

Digit Index:                                          {  1,    2,    3,    4,    5,  6,  7,  8,  9 }
Digits Count totaling 40 values:          { 15,    8,    6,    4,    4,  0,  2,  1,  0 }
Proportions of Digits with '%' sign omitted:   {38,   20,  15,  10,  10,  0,  5,  3,  0 }

Assuming (correctly) that these 40 values were collected in a truly random fashion from the large data set of all 19,452 earthquakes occurrences in 2012; without any bias or attempt to influence first digits occurrences; and that this pattern is generally found in many other data sets, one then may conclude with the phrase "not all digits are created equal", or rather "not all <u>first</u> digits are created equal", even though this seems to be contrary to intuition and against all common sense.

The focus here is actually on the <u>first meaningful digit</u> – counting from the left side of numbers, excluding any possible encounters of zero digits which only signify ignored exponents in the relevant set of powers of ten of our number system. Therefore, the complete definition of the **First Leading Digit** is the <u>first non-zero digit</u> of any given number on its left-most side. This digit is the first significant one in the number as focus moves from the left-most position towards the right, encountering the first non-zero digit signifying some quantity; hence it is also called the **First Significant Digit**. For 2365 the first leading digit is 2. For 0.00913 the first leading digit is 9 and the zeros are discarded; hence even though strictly-speaking the first digit on the left-most side of 0.00913 is 0, yet, the first significant digit is 9. For the lone integer 8 the leading digit is simply 8. For negative numbers the negative sign is discarded, hence for -715.9 the leading digit is 7. Here are some more illustrative examples:

**6**,719,525    →   digit 6
0.0000**7**61    →   digit 7
-0.**2**81264    →   digit 2
**8**75         →   digit 8
**3**           →   digit 3
**-5**          →   digit 5

For a data set where all the values are greater than or equal to 1, such as in the sample of the earthquaqe data, the first digit on the left-most side of numbers is also the First Leading Digit and the First Significant Digit, and necesarily one of the nine digits {1, 2, 3, 4, 5, 6, 7, 8, 9}; while digit 0 never occurs first on the left-most side.

## [II]   Benford's Law and the Predominance of Low Digits

Benford's Law states that:

$$\text{Probability[First Leading Digit is d]} = \text{LOG}_{10}(1 + 1/d)$$

$\text{LOG}_{10}(1 + 1/1) = \text{LOG}(2.00) = 0.301$
$\text{LOG}_{10}(1 + 1/2) = \text{LOG}(1.50) = 0.176$
$\text{LOG}_{10}(1 + 1/3) = \text{LOG}(1.33) = 0.125$
$\text{LOG}_{10}(1 + 1/4) = \text{LOG}(1.25) = 0.097$
$\text{LOG}_{10}(1 + 1/5) = \text{LOG}(1.20) = 0.079$
$\text{LOG}_{10}(1 + 1/6) = \text{LOG}(1.17) = 0.067$
$\text{LOG}_{10}(1 + 1/7) = \text{LOG}(1.14) = 0.058$
$\text{LOG}_{10}(1 + 1/8) = \text{LOG}(1.13) = 0.051$
$\text{LOG}_{10}(1 + 1/9) = \text{LOG}(1.11) = 0.046$
---------
1.000

Figure C depicts the distribution. Figure D visually depicts Benford's Law as a bar chart. This set of nine proportions of Benford's Law is sometimes referred to in the literature as **'The Logarithmic Distribution'**. Remarkably, Benford's Law is confirmed in almost all real-life data sets with high order of magnitude, such as in data relating to physics, chemistry, astronomy, economics, finance, accounting, geology, biology, engineering, governmental census data, and many others.

| Digit | Probability |
|-------|-------------|
| 1 | 30.1% |
| 2 | 17.6% |
| 3 | 12.5% |
| 4 | 9.7% |
| 5 | 7.9% |
| 6 | 6.7% |
| 7 | 5.8% |
| 8 | 5.1% |
| 9 | 4.6% |

**Figure C**:  Benford's Law for First Digits

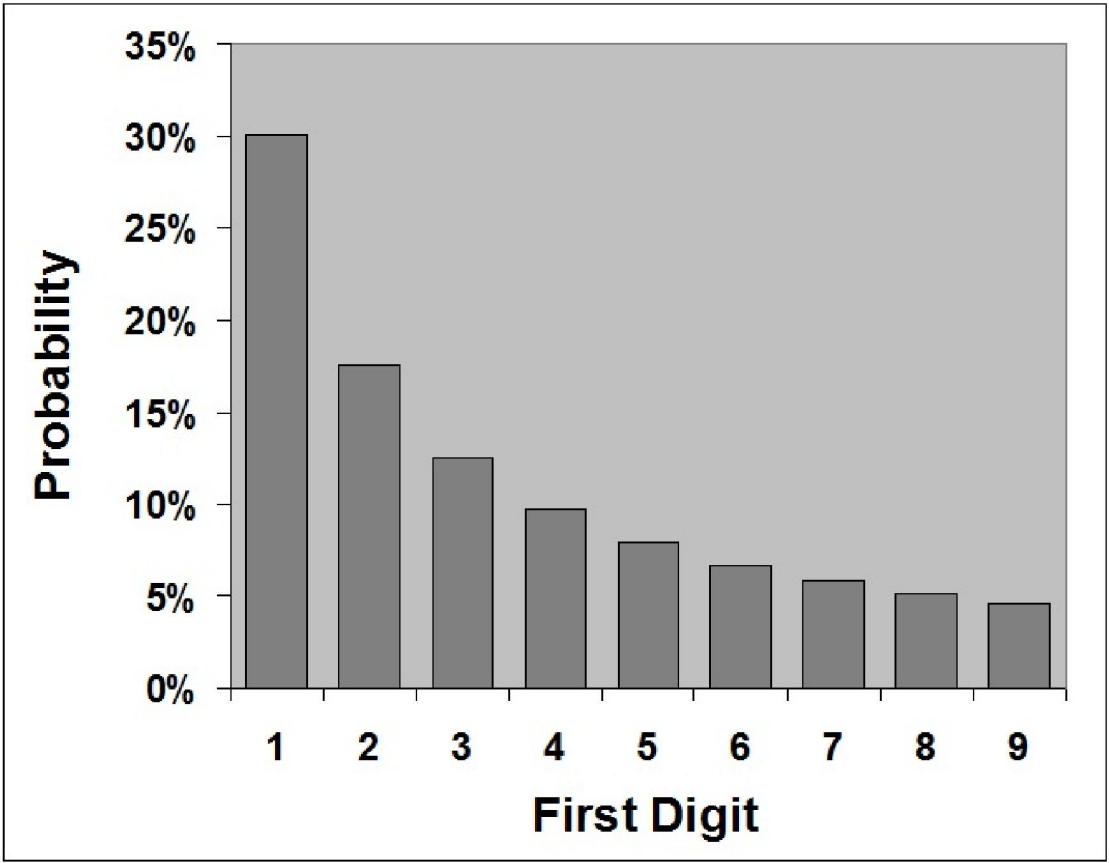

**Figure D**:  Benford's Law – Probability of First Leading Digit Occurrences as a Bar Chart

## [III]  Sum of Squares Deviation Measure (SSD)

It is necessary to establish a standard measure of 'distance' from the Benford digital configuration for any given data set. Such a numerical measure could perhaps tell us about the conformance or divergence from the Benford digital configuration of the data set under consideration. This is accomplished with what is called **Sum Squares Deviations (SSD)** defined as the sum of the squares of the 'errors' between the Benford expectations and the actual/observed values (in percent format – as opposed to fractional/proportional format):

$$ SSD = \sum_{1}^{9} \left( \text{Observed \% of d} - 100 * \text{LOG}(1 + \frac{1}{d}) \right)^2 $$

For example, for observed 1st digits proportions of {31.1, 18.2, 13.3, 9.4, 7.2, 6.3, 5.9, 4.5, 4.1} with '%' sign omitted, SSD measure of distance from the logarithmic is calculated as:

**SSD** = $(31.1 - \mathbf{30.1})^2 + (18.2 - \mathbf{17.6})^2 + (13.3 - \mathbf{12.5})^2 + (9.4 - \mathbf{9.7})^2 +$
$+ (7.2 - \mathbf{7.9})^2 + (6.3 - \mathbf{6.7})^2 + (5.9 - \mathbf{5.8})^2 + (4.5 - \mathbf{5.1})^2 + (4.1 - \mathbf{4.6})^2 = \mathbf{3.4}$

SSD generally should be below 25; a data set with SSD over 100 is considered to deviate too much from Benford; and a reading below 2 is considered to be ideally Benford.

## [IV]  Integral Powers of Ten (IPOT)

Integral Powers of Ten (IPOT) play a crucial role in the understanding of Benford's Law.

An integral power of ten is simply $10^{\text{INTEGER}}$ with either negative or positive integer, as well as zero power. For example, the IPOT numbers 0.001, 0.01, 0.1, 1, 10, 100, are directly derived from $10^{-3}, 10^{-2}, 10^{-1}, 10^{0}, 10^{1}, 10^{2}$.

Adjacent integral powers of ten are a pair of two neighboring and consecutive IPOT numbers $10^{\text{INTEGER}}$ and $10^{\text{INTEGER} + 1}$ such as 1 & 10, or 100 & 1000, and so forth.

# [V]  Benford's Law as Uniformity of Mantissa

In the most simplistic way, mantissa could be described as 'the fractional part of the log', although this definition is not true for numbers less than 1 having negative log values.

The formal definition of the mantissa of any positive number X is that unique solution to $X = 10^C * 10^{MANTISSA}$. Alternatively, mantissa of X is the solution to $X = 10^{C + MANTISSA}$. Here C which is called the 'characteristic' is obtained by rounding down LOG(X) to the nearest integer, namely the largest integer less than or equal to LOG(X). Equivalently, C is the first integer to the left on the log-axis, and regardless whether LOG(X) is negative or positive.

Taking log to both sides of the above mantissa-defining equation we obtain the simple result:
LOG(X) = LOG($10^{C + MANTISSA}$)
LOG(X) = C + MANTISSA

For X = 870, log is 2.93952, the characteristic is 2, and mantissa is 0. 93952.
LOG(X) = 2.93952 = 2 + 0. 93952

For X = 0.063, log value is −1.20066, the characteristic is -2, and mantissa is 0.79934.
LOG(X) = −1.20066 = −2 + 0.79934

For X ≥ 1, mantissa is the fractional part of LOG(X), and the characteristic C is the integral part of LOG(X). For 0 < X < 1, mantissa is 1 minus the fractional part of the absolute value of LOG(X). More generality, mantissa of X can be viewed as the distance on the log-axis between LOG (X) and the integer immediately to the left of it.

Let us consider the following four numbers and their corresponding log values:

8.5274      log = 0.9308

85.274      log = 1.9308

852.74      log = 2.9308

8527.4      log = 3.9308

The fractional part of the logarithm [i.e. mantissa] is the same for all four numbers, namely 0.9308. Only the integral part of the logarithm changes by one integer at a time. This is so since each consecutive number is obtained by multiplying the previous number by a factor of 10 (i.e. moving the decimal point once to the right), and which does not change digital configuration in any way. The first digit is always 8 for all of these four numbers. The same is true for the second digit being 5, third digit being 2, fourth digit being 7, and fifth digit being 4.

To summarize, mantissa cycles from 0 to 1 on each adjacent integral powers of ten sub-interval, such as [1, 10), [10, 100), [100, 1000), and so forth. Mantissa repeats itself over and over again between those crucial points on the x-axis. But so do all the digits! They repeat themselves there!

In other words, digital configuration (all orders considered) and mantissa are basically two distinct ways to indicate the same concept. Clearly, mantissa has a one-to-one correspondence with digital configuration (all orders considered). Therefore, instead of stating in great details how all the digital orders are distributed, one might as well consider [more concisely] how mantissa is distributed, and as a consequence all digital orders are determined in one fell swoop!

As discussed in Kossovsky (2014) chapters 61, 62, and 63, the main result here is:

**Benford's Law implies uniformity of mantissa.**
**Uniformity of mantissa implies Benford's Law.**

Hence, converting any large Benford data set from normal numbers and into a set of pure mantissa values, would yield an approximate uniform (flat and horizontal) distribution on the mantissa space of (0, 1), where probability is roughly equal for all mantissa values.

The entire area under the density curve of the mantissa must be set to 1, as in all statistical distributions, and here for mantissa the range is $(1 - 0) = 1$, thus if the data set is Benford and mantissa is uniformly distributed, then density height is a constant at 1, namely PDF(m) = 1, and areas [representing digital probabilities for the whatever digits] are calculated directly and simply via the reading of the widths on the M-axis. It follows that:

If a given data set or distribution obeys Benford's Law then:
**Probability($M_1$ < Mantissa < $M_2$) = $M_2 - M_1$**

Uniformity of mantissa is often referred to as the "General Form of Benford's Law". This is so since uniformity of mantissa directly implies 1st, 2nd, 3rd, and all higher orders digit destitutions, while LOG(1 + 1/d) only refers to the 1st order digit distribution. In addition, the general law also implies all the probability dependencies between the digital orders.

The expression for the first digits in Benford's Law could be re-written as:

LOG(1 + 1/d)
LOG(d/d + 1/d)
LOG((d + 1)/d)
LOG(d + 1) – LOG(d)

Probability[First Digit is 1] = LOG(2) – LOG(1)
Probability[First Digit is 2] = LOG(3) – LOG(2)
Probability[First Digit is 3] = LOG(4) – LOG(3)

Hence, the probability set {30.1%, 17.6%, 12.5%, 9.7%, 7.9%, 6.7%, 5.8%, 5.1%, 4.6%} of the first order in Benford's Law can be thought of as differences in the values of the logarithms of the natural numbers 1, 2, 3, 4, 5, 6, 7, 8, 9, 10, namely differences in the set {log(1), log(2), log(3), log(4), log(5), log(6), log(7), log(8), log(9), log(10)}, which are the differences in the set {0.000, 0.301, 0.477, 0.602, 0.699, 0.778, 0.845, 0.903, 0.954, 1.000}. The sequence of these fractions above is also the set of the cumulative values of LOG(1 + 1/d).

There are nine compartments within the [0, 1) mantissa space, where each compartment points to a unique first digit, and whose probability is directly proportion to its width on the M-axis. These nine compartments are: [0, 0.301), [0.301, 0.477), [0.477, 0.602), [0.602, 0.699), [0.699, 0.778), [0.778, 0.845), [0.845, 0.903), [0.903, 0.954), and [0.954, 1.000).

Figure E depicts the nine distinct compartments of mantissa corresponding to the nine possible first digits. The mix of the two notations 'LOG' as well as "Mantissa" for the horizontal axis is correct in the case where data happens to fall within [1, 10), so that related log of data falls within [0, 1) and in which case the two terms are interchangeable (since here the fractional part of the log is also the log itself).

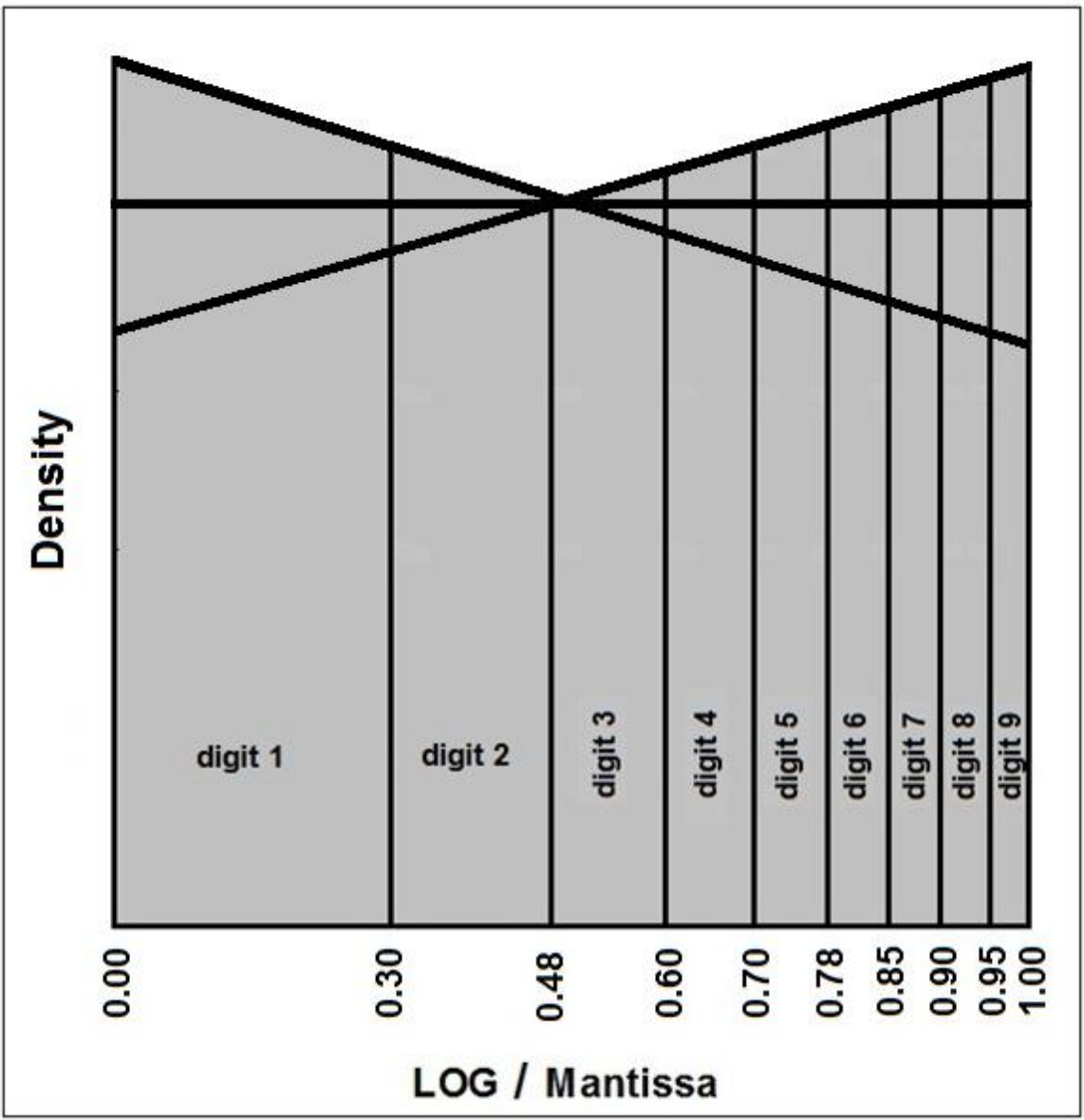

**Figure E**:  Three Hypothetical Densities for Mantissa – Rising, Uniform, and Falling

As an example, digit 2 leads whenever $2 \leq X < 3$, or $20 \leq X < 30$, or $200 \leq X < 300$, and so forth. Taking log of the three parts of the inequality implies that $LOG(2) \leq LOG(X) < LOG(3)$, or $LOG(20) \leq LOG(X) < LOG(30)$, or $LOG(200) \leq LOG(X) < LOG(300)$, and so forth.

These infinite series of inequalities could be written as $LOG(2) \leq LOG(X) < LOG(3)$, or $LOG(2*10) \leq LOG(X) < LOG(3*10)$, or $LOG(2*100) \leq LOG(X) < LOG(3*100)$, and so forth, and therefore they could all be condensed into a singular mantissa inequality of the form $LOG(2) \leq M < LOG(3)$ via a subtraction by the respective characteristic in each series. It follows that digits 2 occurs whenever $0.301 \leq M < 0.477$, so that mantissa related to digit 2 occurrences as the first digit resides within the second compartment of $[0.301, 0.477)$.

As it happened (i.e. Benford's Law), the probability of digit 2 occurring as the first digit is easily calculated as $(0.477 – 0.301) = 0.176$, namely the width of the compartment, it's that simple!

Figure E also depicts three possible densities: rising density, uniform density, and falling density. The case of a rising density corresponds to a more balanced and even distribution of all the nine digits. A much sharper rise in mantissa density than the one seen in Figure E could correspond approximately to digital equality. The flat and horizontal density corresponds to the Benford digital configuration. The case of a falling density corresponds to an extreme digital inequality in favor of low digits, where digit 1 typically leads by 40% or even 60% - much more so than its 30% allocation according Benford's Law.

## [VI]  Uniqueness of K/X Distribution and Connection to Exponential Growth

Of all the many known distributions in mathematical statistics, there is but only one unique distribution whose logarithmically-transformed values are uniformly distributed - and that distribution is k/x. Indeed, one can define k/x as such. The implication is that Monte Carlo computer simulations of say 100,000 values from the k/x distribution, all transformed into their 100,000 log equivalences, would yield nearly uniform and even distribution on the log-axis. Since the Benford configuration seeks uniformity of mantissa, and since mantissa is just the fractional part of the logarithm when $X \geq 1$, then perhaps the uniformity of the logarithm of the k/x distribution might constitute a reasonable explanation for its Benford configuration - assuming some particular ranges on the x-axis in its definition.

As discussed in detail in Kossovsky (2014) chapters 60, 61, 62, 72, and 80, the k/x distribution is the only density that perfectly obeys Benford's Law for a range standing between two adjacent IPOT points, such as (1, 10), (10, 100), or (100, 1000), and so forth. For such adjacent IPOT ranges, there exists no other distribution that perfectly obeys Benford's Law (with all higher orders considered) except k/x distribution. On such particular intervals k/x is unique.

The feature that makes the k/x distribution so unique in this context is that the density of the logarithms values of k/x distribution is uniformly and evenly distributed; and this is so regardless of the range k/x is defined over, hence mantissa could also be uniformly distributed and Benford's Law perfectly obeyed - assuming properly defined ranges. If the range is restricted to two adjacent IPOT such as (1, 10) or (100, 1000) say, then the prerequisite for Benford behavior in having uniformity of mantissa is uniquely achieved via the k/x distribution, since its log distribution is always uniform, and the fractional part of its uniform log (i.e. mantissa) starts at 0.0000 and ends at 0.9999 on these adjacent IPOT intervals.

More generally, the k/x distribution is the only density that perfectly obeys Benford's Law for a range (P, Q) such that log difference LOG(Q) – LOG(P) is unity. Such a range can be written as (P, 10*P), and LOG(10*P) – LOG(P) = LOG(10) + LOG(P) – LOG(P) = LOG(10) = 1.

For example, over the intervals (4.379, 43.79) or (851.23, 8512.3), k/x is that the unique distribution that yields perfect Benford behavior; no other distribution on this range is Benford.

For any generic adjacent IPOT interval, namely $(10^{\text{INTEGER}}, 10^{\text{INTEGER}+1})$, log difference is $\text{LOG}(10^{\text{INTEGER}+1}) - \text{LOG}(10^{\text{INTEGER}}) = (\text{INTEGER} + 1)*\text{LOG}(10) - (\text{INTEGER})*\text{LOG}(10) = (\text{INTEGER} + 1)*1 - (\text{INTEGER})*1 = 1$.

It should be noted that k/x is also perfectly Benford whenever it is defined between any two points A and B such that log difference LOG(B) – LOG(A) is an integer greater than 1, such as say the interval (1.22835, 12283.5) where log difference is the integral value of 4, but k/x is not unique on such wider interval, and there are in principle infinitely many other distributions that are perfectly Benford as well.

As it happened, log values of an exponential growth series are also uniformly spread along the log-axis, albeit in a discrete fashion, and where distances between log values are constant. Let us prove the above assertion using the following notations for exponential progression: B for the base value/quantity at time 0, and F as the factor of growth, such as F = 1.05 for 5% yearly growth. The logarithm of the growing quantity is taken at the beginning of each year or period of growth as follows:

At the 1st year:  LOG(B)
At the 2nd year: LOG(B*F)
At the 3rd year:  LOG(B*F*F)
At the 4th year:  LOG(B*F*F*F)
At the 5th year:  LOG(B*F*F*F*F)

At the 1st year:  LOG(B)
At the 2nd year: LOG(B) + LOG(F)
At the 3rd year: LOG(B) + LOG(F) + LOG(F)
At the 4th year:  LOG(B) + LOG(F) + LOG(F) + LOG(F)
At the 5th year:  LOG(B) + LOG(F) + LOG(F) + LOG(F) + LOG(F)

Clearly distances between log values here are separated by a fixed value, namely by LOG(F), and this implies that the discrete distribution or concentration of log of exponential growth is uniformly spread throughout the log-axis. Hence the (surprising) intimate connection between k/x distribution and exponential growth series! The continuous distribution of the logarithm of k/x, as well as the discrete distribution of the logarithm of exponential growth, are both uniformly and evenly spread throughout the log-axis.

# PART 2:
# EXPONENTIAL GROWTH SERIES

# [1] Normal Versus Anomalous Exponential Growth Series

Exponential growth occurs when the growth of a given quantity [the added value due to growth] is proportional to the quantity's current value, so that 5% exponential growth for example yields only 5 when the quantity is 100, but it yields 50 when the quantity is 1000.

Discrete Exponential growth series are expressed as:

$$\{B, BF, BF^2, BF^3, \ldots, BF^N\}$$

where B is the base (initial) value, P is the constant percent growth rate, and F is the constant multiplicative factor relating to growth rate as in F = (1 + P/100). Throughout this article we shall operate under the assumption of our decimal base 10 number system, although of course everything can be generalized to other bases.

For example, for initial base value B = 10, and 7% growth rate per period (or year), the implied F factor is F = (1 + 7/100) = 1.07. We let quantity 10 grow for 34 years, hence the series is:

$$\{10, \; 10*1.07, \; 10*1.07^2, \; 10*1.07^3, \; 10*1.07^4, \; \ldots, \; 10*1.07^{34}\} =$$

| | | | | | | | | | | | |
|---|---|---|---|---|---|---|---|---|---|---|---|
| {10.0 | 10.7 | 11.4 | 12.3 | 13.1 | 14.0 | 15.0 | 16.1 | 17.2 | 18.4 | 19.7 | 21.0 |
| 22.5 | 24.1 | 25.8 | 27.6 | 29.5 | 31.6 | 33.8 | 36.2 | 38.7 | 41.4 | 44.3 | 47.4 |
| 50.7 | 54.3 | 58.1 | 62.1 | 66.5 | 71.1 | 76.1 | 81.5 | 87.2 | 93.3 | 99.8} | |

First digits are emphasized in bold font:

| | | | | | | | | | | | |
|---|---|---|---|---|---|---|---|---|---|---|---|
| {**1**0.0 | **1**0.7 | **1**1.4 | **1**2.3 | **1**3.1 | **1**4.0 | **1**5.0 | **1**6.1 | **1**7.2 | **1**8.4 | **1**9.7 | **2**1.0 |
| **2**2.5 | **2**4.1 | **2**5.8 | **2**7.6 | **2**9.5 | **3**1.6 | **3**3.8 | **3**6.2 | **3**8.7 | **4**1.4 | **4**4.3 | **4**7.4 |
| **5**0.7 | **5**4.3 | **5**8.1 | **6**2.1 | **6**6.5 | **7**1.1 | **7**6.1 | **8**1.5 | **8**7.2 | **9**3.3 | **9**9.8} | |

The choice of 34 years was deliberate, so as to get as close to Benford as possible by ensuring that the last term is approximately 10-fold the initial base value of 10. In other words, that the exponent difference between maximum 99.8 and minimum 10.0 is LOG(99.8) – LOG(10.0) = 1.9991 – 1.0000 = 0.9991, and which is very close to an integral value (namely integer 1). The phrase 'exponent difference' refers to 'log difference', namely the difference between the log values of the minimum (first element) and the maximum (last element). For example, the exponent of 100 is 2, because LOG(100) is 2, and because $10^2$ = 100, meaning that the exponent necessary to raise 10 to resulting in 100 is exactly 2. For example, for the series {8, 16, 32, 64, 128, 256, 512, 1024}, exponent difference is LOG(1024) – LOG(8) = 3.0103 – 0.9031 = 2.1072. The term Exponent Difference is better known as Order of Magnitude in the literature.

Digit distribution regarding the 7% exponential growth series above for the vector of nine digits {1, 2, 3, 4, 5, 6, 7, 8, 9} with '%' sign omitted is:

7% Growth from 10 for 34 Periods: {31.4, 17.1, 11.4, 8.6, 8.6, 5.7, 5.7, 5.7, 5.7}
Benford's Law 1st Digits Order:    {30.1, 17.6, 12.5, 9.7, 7.9, 6.7, 5.8, 5.1, 4.6}

Quoting from Kossovsky (2014) chapter 98 titled "Singularities in Exponential Growth Series": "Due to their discrete and finite nature, real-life exponential growth series can never be formally considered as exactly logarithmic by pure mathematicians, only in the limit possibly". Discreteness in general causes deviations from the logarithmic, and could possibly be remedied by considering infinitely many elements and the limiting case. For mere mortals such as statisticians and scientists, the pressing issue is not some highly theoretical and mathematically exact digital behavior, but rather the practical consideration of an approximate Benford behavior. Does it matter much that digit 1 has earned 29.8% or 30.5% while digit 9 has earned only about say 4.2%? No! Because such is the typical behavior of all real-life data sets obeying Benord's Law, and no one should be disturbed at all if the exact theoretical 30.1% proportion of Benford is not found. Hence, in extreme generality, what is [practically] necessary for exponential growth series to behave approximately or nearly logarithmically is having plenty of elements, although for low growth rates there is an additional requirement, namely that exponent difference between minimum and maximum should be as close as possible to an integer. One straightforward way of achieving this last requirement for low growth rates is to make sure that the series starts and ends approximately at integral powers of ten, such as 10 & 100, 1 & 1000, or 100 & 1000000, and so forth. Additional discussions about this issue are given in later chapters of this article.

Let us further explore exponential growth series in general by examining how logarithm values of the series progress forward at each new period. As a consequence, this exploration would indirectly lead to better understanding of what happens to mantissa and its resultant distribution after numerous periods. The decimal log series of any exponential growth series is simply:

{ LOG(B),
   LOG(B) + 1*LOG(F),
   LOG(B) + 2*LOG(F),
   LOG(B) + 3*LOG(F),

        , … ,

   LOG(B) + N*LOG(F) }

It is certainly proper to think of LOG(F) namely LOG(1 + P/100) as being a fraction, as it is so in any case up to 900% growth, else the integral part can be ignored as far as digital configuration is concerned. For example, for 2% growth, LOG(1 + 2/100) = 0.009. For 50% growth, LOG(1 + 50/100) = 0.176. Even for 180% growth, LOG(1 + 180/100) = 0.447 which is a fraction. Therefore, related log series of the exponential growth series is simply a series with constant additions (accumulation) of LOG(F) from an initial base of LOG(B). While this related log series grows ever larger, mantissa on the other hand is constantly oscillating between 0 and 1, as it typically (for low growth) takes many small steps forward, then suddenly one large leap backwards whenever log overflows an integer, and so forth. This is so since mantissa is obtained by constantly removing the whole part of the log *(whenever growth series ≥ 1, that is whenever LOG(series) is positive or zero, an assumption which can be taken for granted here)*. Therefore, mantissa is clearly seen as being uniform on (0, 1) as more and more points of newly minted mantissa keep falling there on a variety of points, covering an ever increasing 'portion' of the entire (0, 1) range. Even though the process of exponential growth series is truly deterministic, if one were to visually follow these 'rapid' mantissa additions onto (0, 1) space without getting severely dizzy, it would all seem quite random, disorganized, and highly chaotic, and it is precisely this nature of mantissa creation that yields uniformity to its final overall distribution! As an example, we examine an exponential 30% growth series, starting from an initial base 3, having factor F = 1.30, and with the implied LOG(F) = 0.1139433. Figure 1 depicts part of that series in details. What should be carefully noted here is that mantissa always re-enters into (0, 1) interval at different (newer) locations, namely at 0.047, 0.072, 0.098, and so on. This is a necessary condition for logarithmic behavior, as it guarantees that mantissa is well-spread over the entire (0, 1) range in an even and uniform manner, covering all corners and segments. In a sense it guarantees that mantissa creation is random and unorganized with respect to (0, 1) space, and thus that (almost) always new mantissa values are being created.

| Series | Log | Mantissa |
|---|---|---|
| 3.0 | 0.477 | 0.477 |
| 3.9 | 0.591 | 0.591 |
| 5.1 | 0.705 | 0.705 |
| 6.6 | 0.819 | 0.819 |
| 8.6 | 0.933 | 0.933 |
| 11.1 | 1.047 | 0.047 |
| 14.5 | 1.161 | 0.161 |
| 18.8 | 1.275 | 0.275 |
| 24.5 | 1.389 | 0.389 |
| 31.8 | 1.503 | 0.503 |
| 41.4 | 1.617 | 0.617 |
| 53.8 | 1.730 | 0.730 |
| 69.9 | 1.844 | 0.844 |
| 90.9 | 1.958 | 0.958 |
| 118.1 | 2.072 | 0.072 |
| 153.6 | 2.186 | 0.186 |
| 199.6 | 2.300 | 0.300 |
| 259.5 | 2.414 | 0.414 |
| 337.4 | 2.528 | 0.528 |
| 438.6 | 2.642 | 0.642 |
| 570.1 | 2.756 | 0.756 |
| 741.2 | 2.870 | 0.870 |
| 963.6 | 2.984 | 0.984 |
| 1252.6 | 3.098 | 0.098 |
| 1628.4 | 3.212 | 0.212 |
| 2116.9 | 3.326 | 0.326 |

**Figure 1:** Normal Logarithmic Exponential Series with 30% Growth

Figure 2 depicts the way the 30% growth series almost always re-enter differently after each integer mark on the log-axis.

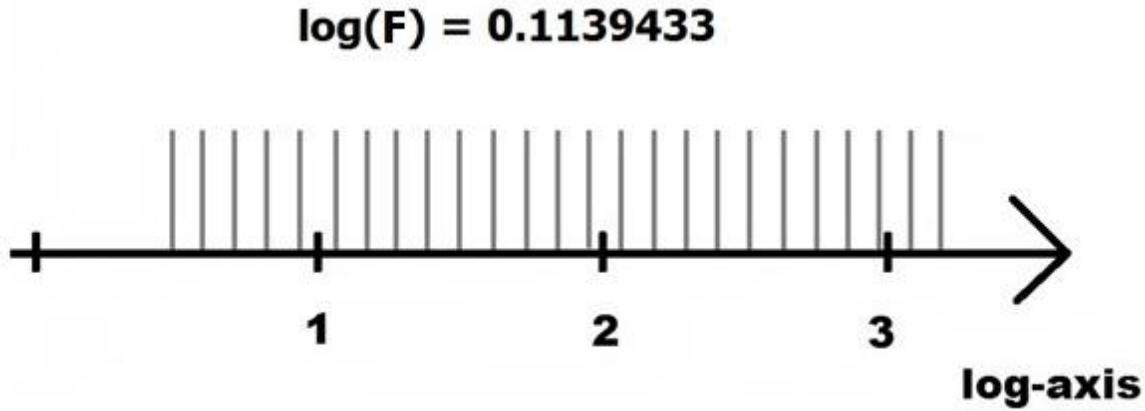

**Figure 2:** The March on the Log-Axis of Normal Logarithmic Series with 30% Growth

Such logarithm and mantissa vista leads to the detection of some peculiar digital singularities. The argument above falls apart whenever the fraction LOG(F) happens to be such that exactly T whole multiples of it add up to unity, as in the fractions 0.50, 0.25, 0.10, 0.125, 0.05, and so forth, in which case constant LOG (F) additions lead to re-entering (0, 1) always at the same point, and taking the same type of steps over and over again focusing only on a few selected fortunate points having some very strong concentration (density), all the while ignoring all the other points or sections on the interval (0, 1). Such growth rate results in some quite uneven mantissa distribution on (0, 1) and yields decisively non-logarithmic digital distribution for the exponential series itself. In other words, beyond the creation of a few initial mantissa points, the series does not create any new mantissa, but just repeats those few old ones over and over again. Algebraically, the series is non-logarithmic and **rebellious** (**anomalous**) whenever there exists an integer T such that LOG(F)*T = 1. One should always be reminded that only uniformity of mantissa yields logarithmic behavior.

As an example, we examine one such anomalous exponential growth series starting from the initial base value of 10, growing at 77.8279% per period, and thus having factor F = 1.778279 and the implied problematic LOG(F) = LOG(1.778279) = 0.25 = 1/4, where exactly 4 (called T) multiples of it add up to 1. Figure 3 depicts part of that series in details. Mantissa always re-enters into (0, 1) interval at the same location, namely at 0.000, and then always takes the same subsequent 'long' steps landing at 0.250, 0.500, and 0.750, intentionally skipping 'numerous' points in between them. Obviously such state of affairs cannot result in any logarithmic behavior because mantissa is not uniform, and the series is anomalous. This conclusion is so regardless of the initial value (base) of the series, and the fact that here the series starts at 10 is irrelevant. For example, for the 77.8279% anomalous growth series above, had the base value been 8 instead of 10, then mantissa would always re-enter into (0, 1) interval at 0.153, and then land exclusively on the points 0.403, 0.653, 0.903, over and over again. These points are all 0.25 units apart.

| Series | Log | Mantissa |
|---|---|---|
| 10.0 | 1.000 | 0.000 |
| 17.8 | 1.250 | 0.250 |
| 31.6 | 1.500 | 0.500 |
| 56.2 | 1.750 | 0.750 |
| 100.0 | 2.000 | 0.000 |
| 177.8 | 2.250 | 0.250 |
| 316.2 | 2.500 | 0.500 |
| 562.3 | 2.750 | 0.750 |
| 1000.0 | 3.000 | 0.000 |
| 1778.3 | 3.250 | 0.250 |
| 3162.3 | 3.500 | 0.500 |
| 5623.4 | 3.750 | 0.750 |
| 10000.0 | 4.000 | 0.000 |
| 17782.8 | 4.250 | 0.250 |
| 31622.8 | 4.500 | 0.500 |
| 56234.1 | 4.750 | 0.750 |
| 100000.0 | 5.000 | 0.000 |
| 177827.9 | 5.250 | 0.250 |
| 316227.8 | 5.500 | 0.500 |
| 562341.3 | 5.750 | 0.750 |

**Figure 3:** Anomalous Non-Logarithmic Series with 77.8279% Growth

Figure 4 depicts the 77.8279% anomalous growth series with base value of 10 as it re-enters the log-axis always at the same location after each integer mark.

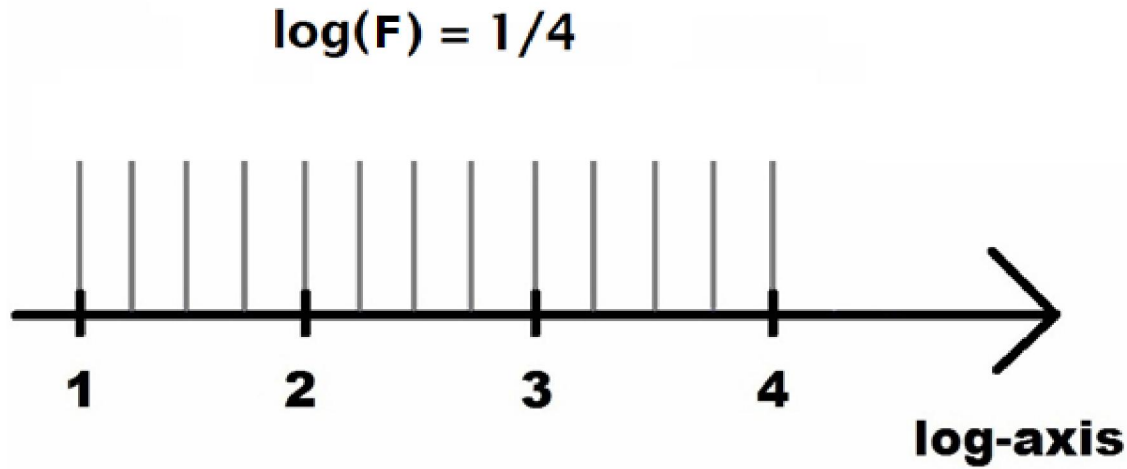

**Figure 4:** The March on the Log-Axis of Anomalous 77.8279% Growth Series - Basic Type

Yet, even for those rebellious non-logarithmic rates, some comfort and relief can be found whenever LOG(F) [the width of the steps by which log of the series advances along the log-axis] is sufficiently small as compared with unit length, because then no matter how repetitive, peculiar, and picky mantissa chooses the points upon which to stamp on (0, 1), each step is still so tiny that it has no choice but to walk almost all over the interval covering most corners and locations of (0, 1). Simply put: the creature is such that its legs are too short to be truly picky, so it cannot jump and skip much, and it is reduced to walking almost all over the interval, willingly or unwillingly. Only a walking creature with long legs can be effectively picky and successfully avoid certain segments lying on the ground. A giraffe can successfully avoid a 50-centimeter hole or gap on the ground, but a tiny ladybug cannot, no matter how carefully it walks.

Therefore, as the value of LOG(F) decreases, becoming very small, say less than 0.01 (a rational 1/100, designated as 1/T), its spread over (0, 1) is fairly good and quite even, so that its digit distribution is very nearly logarithmic, even though upon closer examination it still stamps cautiously and discretely on the log/mantissa-axis in even but tiny steps of 0.01-width each.

Figure 5 depicts such a scenario with LOG(F) = 1/100 where the growth rate is 2.329%. The percent growth rate is calculated via $1/100 = \text{LOG}(1 + \text{Percent}/100)$, then taking 10 to the power of both sides of the equation we get: $10^{1/100} = 10^{\text{LOG}(1 + \text{Percent}/100)}$, therefore $10^{1/100} = (1 + \text{Percent}/100)$, and $10^{1/100} - 1 = \text{Percent}/100$, so that $100[10^{1/100} - 1] = 2.329\%$. Digit distribution is {30.0, 17.9, 12.1, 10.0, 8.0, 6.0, 6.0, 5.0, 5.0} for such a series starting at an initial base 3 and growing for 1000 periods; where SSD comes out as 1.1, and such very low SSD value indicates nearly perfect Benford behavior.

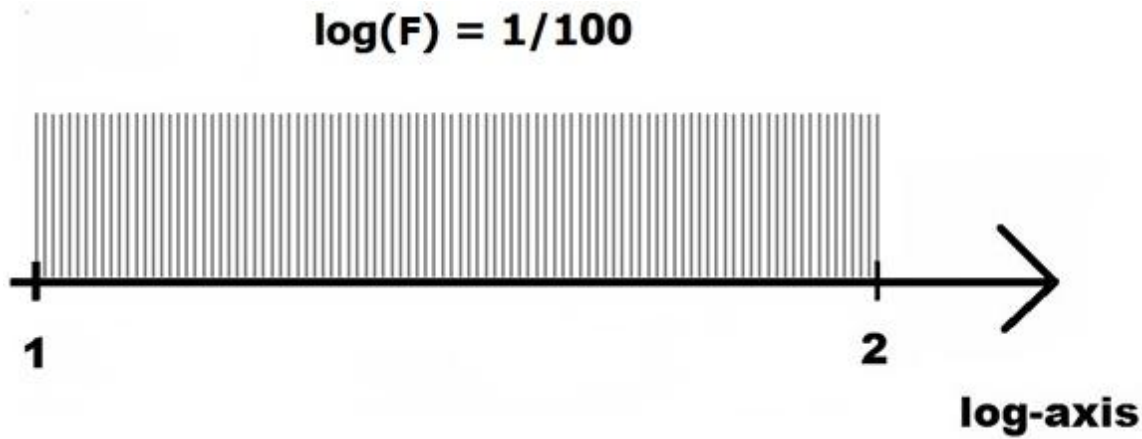

**Figure 5:** Uniformity of Mantissa is Achieved via Small Rational Steps

Clearly, the series with LOG(F) = 1/N in the limit as N goes to infinity is perfectly logarithmic **by definition**, since this is indeed the very meaning of uniformity of mantissa!

For example, let us consider the rational series with LOG(F) = 1/1,000,000. Here for each tiny sub-interval on the log-axis with width of 0.0005 say, there are always (that is: uniformly) 500 points. Hence the interval on (0.0005, 0.0010) has 500 points; the interval on (0.0010, 0.0015) has 500 points; the interval on (0.8005, 0.8010) has 500 points as well, and so forth for any interval of 0.0005 width, and which is the very meaning of the density being uniform! No excuses, explanations, or reference to exceptions, are needed in stating this, except that formally the definition is stated in a continuous – infinitely refined – sense. Hence surely, for the above series with LOG(F) = 1/1,000,000, if one insists on scrutinizing density fanatically on even tinier sub-intervals such as (0.0000000000001, 0.0000000000002) then mantissa is not uniform, although it is uniform indeed for all practical purposes.

Figure 6 depicts various anomalous growth rates by varying integer T from 1 to 20, as well as evaluating T at 25, 35, 40, 50, 100, and 200. For each anomalous growth rate in this table, an actual computer simulation (calculation rather) of the relevant exponential series is run, using the first 1000 elements, all starting from the initial base value of 3. Digital results from such computer runs are displayed in the table, along with their associated SSD values. Figure 6 demonstrates that eventually when T gets to be over 100 or so, deviations from the logarithmic become very small and insignificant.

| T | LOG(F) | % Growth | 1 | 2 | 3 | 4 | 5 | 6 | 7 | 8 | 9 | SSD |
|---|---|---|---|---|---|---|---|---|---|---|---|---|
| 1 | 1 | 900.000% | 0.0 | 0.0 | 100.0 | 0.0 | 0.0 | 0.0 | 0.0 | 0.0 | 0.0 | 9155.8 |
| 2 | 0.5 | 216.228% | 0.0 | 0.0 | 50.0 | 0.0 | 0.0 | 0.0 | 0.0 | 0.0 | 50.0 | 4947.6 |
| 3 | 0.333 | 115.443% | 33.3 | 0.0 | 33.3 | 0.0 | 0.0 | 33.3 | 0.0 | 0.0 | 0.0 | 1701.8 |
| 4 | 0.25 | 77.828% | 25.0 | 22.7 | 2.3 | 0.0 | 25.0 | 0.0 | 0.0 | 0.0 | 25.0 | 1063.3 |
| 5 | 0.2 | 58.489% | 40.0 | 0.0 | 20.0 | 20.0 | 0.0 | 0.0 | 20.0 | 0.0 | 0.0 | 926.9 |
| 6 | 0.167 | 46.780% | 16.6 | 16.6 | 16.7 | 16.7 | 0.0 | 16.7 | 0.0 | 0.0 | 16.7 | 619.8 |
| 7 | 0.143 | 38.950% | 28.6 | 14.2 | 14.3 | 14.3 | 14.3 | 0.0 | 0.0 | 14.3 | 0.0 | 262.9 |
| 8 | 0.125 | 33.352% | 25.0 | 24.5 | 0.5 | 12.5 | 12.5 | 0.0 | 12.5 | 0.0 | 12.5 | 424.9 |
| 9 | 0.111 | 29.155% | 33.3 | 11.1 | 22.3 | 0.0 | 11.1 | 11.1 | 0.0 | 11.1 | 0.0 | 362.6 |
| 10 | 0.1 | 25.893% | 30.0 | 10.0 | 20.0 | 10.0 | 10.0 | 0.0 | 10.0 | 0.0 | 10.0 | 236.7 |
| 11 | 0.091 | 23.285% | 36.4 | 14.1 | 13.1 | 9.1 | 9.1 | 9.1 | 0.0 | 9.1 | 0.0 | 130.3 |
| 12 | 0.083 | 21.153% | 24.9 | 16.6 | 16.8 | 8.4 | 8.4 | 8.3 | 8.3 | 0.0 | 8.3 | 97.4 |
| 13 | 0.077 | 19.378% | 30.8 | 22.8 | 7.9 | 7.7 | 7.7 | 7.7 | 7.7 | 7.7 | 0.0 | 84.8 |
| 14 | 0.071 | 17.877% | 28.4 | 20.9 | 7.7 | 14.4 | 7.2 | 7.2 | 0.0 | 7.1 | 7.1 | 103.6 |
| 15 | 0.067 | 16.591% | 33.2 | 19.4 | 7.2 | 13.4 | 6.7 | 6.7 | 6.7 | 6.7 | 0.0 | 80.3 |
| 16 | 0.063 | 15.478% | 31.0 | 12.4 | 12.6 | 12.6 | 6.3 | 6.3 | 6.3 | 6.3 | 6.2 | 43.5 |
| 17 | 0.059 | 14.505% | 35.3 | 11.6 | 17.7 | 5.9 | 11.8 | 5.9 | 5.9 | 5.9 | 0.0 | 141.9 |
| 18 | 0.056 | 13.646% | 27.5 | 16.5 | 16.8 | 5.6 | 11.2 | 5.6 | 5.6 | 5.6 | 5.6 | 56.6 |
| 19 | 0.053 | 12.884% | 31.4 | 15.6 | 15.9 | 10.6 | 5.3 | 5.3 | 10.6 | 5.3 | 0.0 | 71.0 |
| 20 | 0.05 | 12.202% | 30.0 | 15.0 | 15.0 | 10.0 | 10.0 | 5.0 | 5.0 | 5.0 | 5.0 | 21.2 |
| 25 | 0.04 | 9.648% | 28.0 | 16.0 | 16.0 | 8.0 | 8.0 | 8.0 | 4.0 | 4.0 | 8.0 | 40.1 |
| 35 | 0.029 | 6.800% | 28.1 | 16.8 | 14.5 | 8.7 | 8.7 | 5.8 | 5.8 | 5.8 | 5.8 | 13.1 |
| 40 | 0.025 | 5.925% | 30.0 | 17.5 | 12.5 | 10.0 | 10.0 | 5.0 | 7.5 | 5.0 | 2.5 | 14.5 |
| 50 | 0.02 | 4.713% | 30.0 | 16.0 | 14.0 | 10.0 | 8.0 | 6.0 | 6.0 | 4.0 | 6.0 | 8.8 |
| 100 | 0.01 | 2.329% | 30.0 | 17.9 | 12.1 | 10.0 | 8.0 | 6.0 | 6.0 | 5.0 | 5.0 | 1.1 |
| 200 | 0.005 | 1.158% | 30.0 | 17.5 | 12.5 | 10.0 | 8.0 | 6.5 | 6.0 | 5.0 | 4.5 | 0.2 |
| Ben | ===== | ====== | 30.1 | 17.6 | 12.5 | 9.7 | 7.9 | 6.7 | 5.8 | 5.1 | 4.6 | 0.0 |

**Figure 6:** Benford is Eventually Found  in Small Steps of Anomalous Rational Series

Anomalous series of the form described above, where the fraction LOG(F) happens to be such that exactly whole multiples of it add up to <u>unity</u>, are actually just one particular type. More generally, whenever whole multiples of LOG(F) add up <u>any integer</u>, be it 2, 3, or any other integral number, we encounter the same dilemma of having uneven mantissa distribution on (0, 1), resulting in non-logarithmic digital distribution for the exponential series itself. To recap, **General Types** of anomalous exponential growth rates are found when exactly **T** whole multiples of LOG(F) add up exactly to any integer **L**, and not just to unity.

For example, 5 whole multiples of 0.4 yield exactly the value of 2 units of distance spanning log/mantissa interval, hence its related 151.1886% growth series has a non-logarithmic digit distribution. The percent growth rate is calculated via $0.4 = 2/5 = LOG(1 + Percent/100)$, which leads to $2/5 = LOG(1 + (151.1886)/100) = LOG(2.511886)$. Figure 7 depicts in details part of that series starting from the initial base value of 10. Figure 8 depicts the selective and repetitive way the series marches along the log-axis, in coordination and synchronization with integer marks. Clearly this rebellious growth series is not Benford at all.

| Series | Log | Mantissa |
|---|---|---|
| 10.0 | 1.000 | 0.000 |
| 25.1 | 1.400 | 0.400 |
| 63.1 | 1.800 | 0.800 |
| 158.5 | 2.200 | 0.200 |
| 398.1 | 2.600 | 0.600 |
| 1000.0 | 3.000 | 0.000 |
| 2511.9 | 3.400 | 0.400 |
| 6309.6 | 3.800 | 0.800 |
| 15848.9 | 4.200 | 0.200 |
| 39810.7 | 4.600 | 0.600 |
| 100000.0 | 5.000 | 0.000 |
| 251188.6 | 5.400 | 0.400 |
| 630957.3 | 5.800 | 0.800 |
| 1584893.2 | 6.200 | 0.200 |
| 3981071.7 | 6.600 | 0.600 |
| 10000000.0 | 7.000 | 0.000 |
| 25118864.3 | 7.400 | 0.400 |
| 63095734.4 | 7.800 | 0.800 |
| 158489319.2 | 8.200 | 0.200 |
| 398107170.6 | 8.600 | 0.600 |

**Figure 7:** Anomalous 151.1886% Growth Series

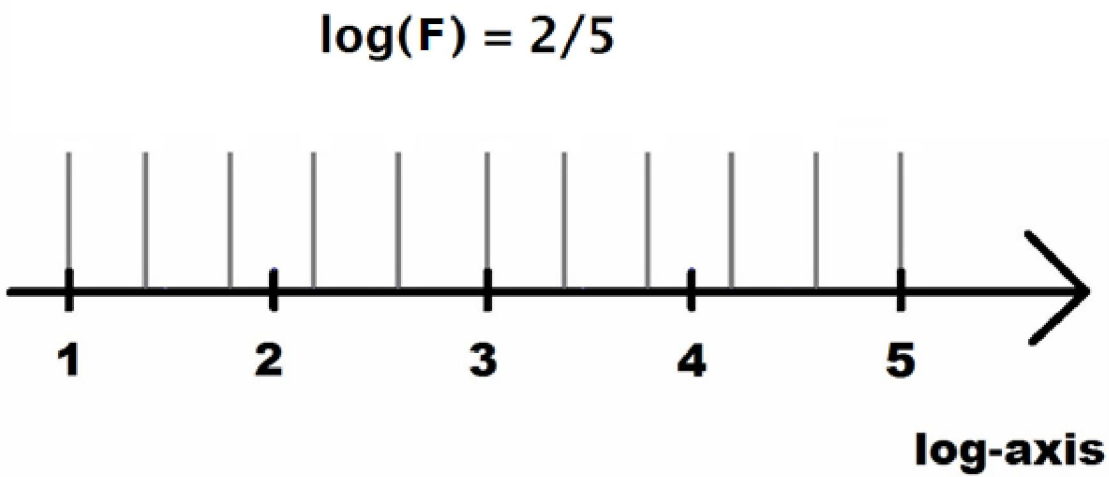

**Figure 8:** The March on the Log-Axis of Anomalous 151.1886% Series - General Type

The general rule regarding Benford behavior for any exponential growth series is as follows: whenever the fraction LOG(1 + P/100) equals the rational number L/T, non-logarithmic behavior for the series with P% growth rate is found. In other words, for both L & T being positive integers, the exponential series is rebellious whenever:

# LOG(F) = L/T

The **Basic Types** anomalous exponential series discussed earlier, such as the one of Figures 3 and 4, are those where L = 1 and where the rationality rule there is that the series is not Benford whenever LOG(F) = 1/T with T being a positive integer.

Conceptually it helps keeping in mind that for rebellious rates, either of the general type or of the basic type, the interpretation of the rationality condition can be stated as follows:

# LOG(F) = [in L whole log units] / [we cycle exactly T periods]

The set of first digits obtained in one cycle of T periods, will repeat itself in the next cycle of T periods, and also in the next one, and so forth. This is so since the same series values are repeated in each cycle of T periods, except for the decimal point which moves once to the right (namely that mantissa values are perfectly repeated in each cycle of T periods). This is why T plays such a decisive role in determining the magnitude of the deviation from Benford, while L is almost totally irrelevant to the magnitude of deviation.

For large value of T, say over 100, the series has already a fairly large number (namely T) of distinct mantissa values and a nearly logarithmically decent set of first digits, and repeating it on each new cycle of T periods does not diminish its approximate Benfordness.

For a very small value of T, say 5, the series repeats these 5 first digits over and over again on each cycle of T periods, and there is no hope of ever breaking out of this vicious digital cycle.

It should be noted that L > T cases are distinct real possibilities, although they all correspond to extremely high growth rates over 900%. Cases where L = T correspond to LOG(F) = 1, which implies that F = 10, so that (1 + P/100) = 10, and therefore P = 900%.

In addition, care should be exercised always to avoid non-reduced forms such as L/T = 2/214 and L/T = 500/11200, which are really of 1/107 and 5/112 reduced forms respectively.

The non-reduced form 4/10 for example, with L = 4 and T = 10, is associated with the interpretation that L/T = (in 4 whole log units)/(we cycle 10 periods), but upon careful examination of how such a series marches along the log-axis one could note that the more basic interpretation is L/T = (in 2 whole log units)/(we cycle 5 periods), and that the 4/10 interpretation is simply a pair of two distinct 2/5 cycles mistakenly considered as one. Non-reduced forms should be avoided as they do not represent any new anomalous possibilities distinct from the sufficient and complete set of all reduced formed possibilities.

## [2]  The Simplistic Explanation of the Rarity of Anomalous Series

On the face of it, we can expect not merely most, but rather practically all exponential growth series to behave logarithmically. This appears to be so simply because that for growth rate to be rebellious the fraction $LOG(1 + P/100)$ must equal exactly some rational number $L/T$, and this is quite rare. In fact, Mathematics teaches us that when a number is picked at random from any continuous set of real numbers, there is zero probability of obtaining a rational number and 100% probability of obtaining an irrational number. Figure 9 humorously dramatizes the argument.

Admittedly, quite often, financial rates, economics-related rates, and others, are quoted rationally as fractions, such as say 5½ % yearly interest rate on a 30-year bond investment, but one should not confuse the set of growth rates with the set of LOG operations performed on these rates. Such rational quotes in finance do not usually lead to anomalous series, because it's the equality $LOG(1 + \text{Rational-Rate}/100) = (\text{Rational } L/T)$ which should hold. A possible set of rates on say the real continuous interval (0%, 100%) is mapped into a different set of $LOG(1 + \text{Rate}/100)$ on the real continuous interval $(0, LOG(2))$, and with the same resultant rarity of the rationals.

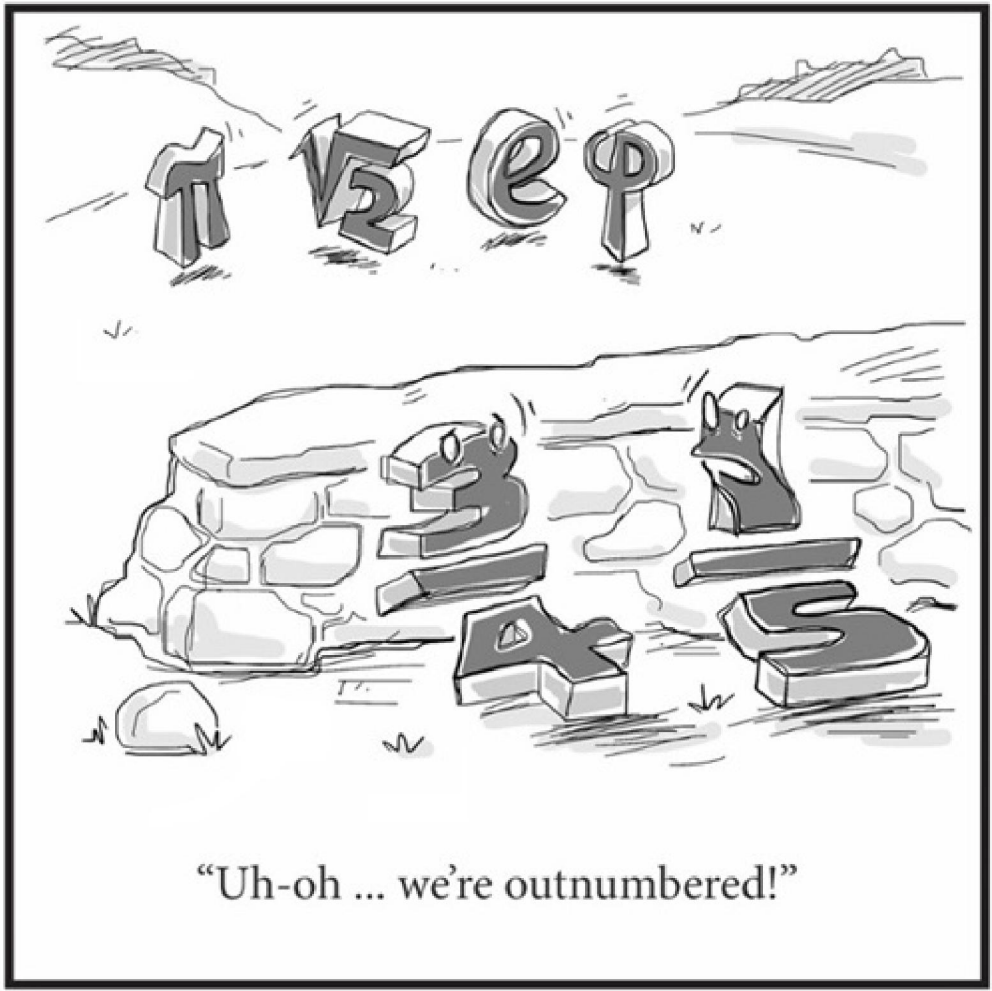

**Figure 9:** The Simplistic Explanation of the Rarity of Anomalous Series

The above argument attempts to utilize the fact that the irrationals vastly outnumber the rationals. Back in the time of Pythagoras around 500 BC, mathematicians refused to believe that irrational numbers even existed. A major principle or dogma of the Pythagorean Sect was the claim that whole numbers ruled the world, and that everything could be explained in terms of integers or fractions of integers. When the Pythagoreans were able to prove that the diagonal of the square is irrational it caused a major intellectual crisis. Many centuries later, at the end of the 1800's, the mathematician Georg Cantor made the startling discovery that the irrational numbers are actually more numerous than the rational numbers! The rational fractions are countably infinite whereas the irrationals are uncountably infinite. The Pythagoreans would have been quite shocked and further dismayed had they learnt of Georg Cantor's discoveries!

## [3] The Indirect Damage Done to Normal Series by Anomalous Series

Yet, the above argument of why exponential growth series are rarely anomalous is too simplistic, as it indirectly relies on the assumption of highly abstract and idealized infinitely-long exponential growth series having infinite number of growth periods, and which is irrelevant to finite real-life exponential growth series. This is so because for an irrational normal growth series that is finite, merely being in the neighborhood of a rational rebellious series [namely having nearly the same growth rate as shown in Figures 10 and 11] is problematic and logarithmic behavior is severely disrupted, unless a large number of growth periods is considered to overcome the closeness to the series (and that number depends on the intensity of the rebelliousness of the anomalous series and on the closeness to it). The only remedy for an irrational normal series of finite length residing near a rational rebellious series is increasing its length, i.e. considering many more growth periods. Such a cure may be too demanding though in terms of the number of growth periods necessary for good behavior, and sometimes even standard personal computers are not sufficiently powerful to perform this task. Hence, for any irrational normal growth series of finite number of growth periods, one can theoretically find a nearby rational rebellious series in the neighborhood - extremely close to it! Would this supposedly render many or most finite exponential growth series rebellious!?

For example, for the rational rebellious series with a growth rate of **93.069773%** and where $LOG(1 + 93.069773/100) = 0.285714286 = 2/7$, starting at the initial base value of 3, and growing for 300 periods, first digits are {28.6, 14.3, 14.3, 14.3, 14.3, 0.0, 0.0, 14.3, 0.0} resulting in high SSD value of 261.7.

Just being in the neighborhood of this 93.069773% rebellious series causes problems. For the irrational normal series **93.15%**, starting at the initial base 3, and growing for 300 periods, first digits are with SSD = 166.5, hence this short and finite series is not Benford. But let us increase the number of periods in order to 'remedy' this deviation and manifest its true Benfordian nature:

300 periods: {28.3, 14.3, 14.3, 14.3, 4.0, 10.3, 0.0, 12.7, 1.7}   SSD = 166.5
500 periods: {28.6, 14.2, 14.4, 12.4, 4.2, 10.6, 1.4, 7.6, 6.6}   SSD =  83.6
600 periods: {30.3, 14.3, 14.3, 10.3, 6.0,  8.8, 3.5, 6.3, 6.0}   SSD =  31.6
700 periods: {30.3, 16.3, 14.1,  9.0, 7.1,  7.6, 5.0, 5.4, 5.1}   SSD =   7.4
900 periods: {29.9 , 17.2, 12.8, 10.2, 8.3,  6.2, 5.1, 6.2, 4.0}   SSD =   3.0

Hence with 900 periods we were able to cure and bring the 93.15% series back to Benford. Yet, getting even closer to the 93.069773% rebellious series, say for the series **93.09%**, starting at base 3, would spell serious trouble even with 900 periods. This is so since this series 'mimics' even better the rebellious series. Here digits are {28.6, 14.3, 14.3, 14.2, 5.3, 9.0, 0.0, 14.2, 0.0} with SSD = 186.5, even after 900 periods. Since the 93.09% series is even closer to the 93.069773% rebellious series as compared with the 93.15% series above, many more periods are needed to remedy it, and merely 900 periods are simply not sufficient for it.

Let us explain why an irrational series in close proximity to a rational one should suffer a similar non-Benfordian fate, at least in the beginning, for the first tens of periods, or perhaps for the first hundreds or so periods.

In order to give a general explanation, let us consider the log values of the above rebellious series 93.069773% having the rationality condition of LOG(F) = 2/7, as well as log values of the above irrational normal series 93.15%, both starting at an initial base 3. The goal is to follow how these two series develop their log values, and by extension to examine how they develop their mantissa values, in the beginning, for the initial 22 periods. The following two sets of data refer to their mantissa values as well as to their first digits (the 93.069773% series is shown on top and the 93.15% series is shown at the bottom below it):

| 1 | 2 | 3 | 4 | 5 | 6 | 7 | 8 | 9 | 10 | 11 | 12 | 13 |
|---|---|---|---|---|---|---|---|---|----|----|----|----|
| 0.477 | 0.763 | 0.049 | 0.334 | 0.620 | 0.906 | 0.191 | 0.477 | 0.763 | 0.049 | 0.334 | 0.620 | 0.906 |
| 0.477 | 0.763 | 0.049 | 0.335 | 0.621 | 0.907 | 0.192 | 0.478 | 0.764 | 0.050 | 0.336 | 0.622 | 0.908 |

| 14 | 15 | 16 | 17 | 18 | 19 | 20 | 21 | 22 |
|----|----|----|----|----|----|----|----|----|
| 0.191 | 0.477 | 0.763 | 0.049 | 0.334 | 0.620 | 0.906 | 0.191 | 0.477 |
| 0.194 | 0.480 | 0.766 | 0.051 | 0.337 | 0.623 | 0.909 | 0.195 | 0.481 |

| 1 | 2 | 3 | 4 | 5 | 6 | 7 | 8 | 9 | 10 | 11 | 12 | 13 |
|---|---|---|---|---|---|---|---|---|----|----|----|----|
| 3 | 5 | 1 | 2 | 4 | 8 | 1 | 3 | 5 | 1 | 2 | 4 | 8 |
| 3 | 5 | 1 | 2 | 4 | 8 | 1 | 3 | 5 | 1 | 2 | 4 | 8 |

| 14 | 15 | 16 | 17 | 18 | 19 | 20 | 21 | 22 |
|----|----|----|----|----|----|----|----|----|
| 1 | 3 | 5 | 1 | 2 | 4 | 8 | 1 | 3 |
| 1 | 3 | 5 | 1 | 2 | 4 | 8 | 1 | 3 |

The 93.15% series grows a tiny bit faster than the 93.069773% series, having an extra ≈0.08% growth rate, hence after 22 periods it reaches the mantissa value of 0.481 which is just slightly higher than the mantissa value of 0.477 of the slightly slower 93.069773% series. Clearly, these two series differ only by very little in how they generate their mantissa values early on.

And what is the difference here on leading digits for these initial 22 periods? None whatsoever! They both generate the same (non-Benfordian) first digits! This is so because the mapping from mantissa to first digits is quite flexible and slack, so that tiny differences in mantissa values usually do not affect first digits. Chapter V of Part I titled 'Benford's Law as Uniformity of Mantissa' discusses the nine compartments within the [0, 1) mantissa space, where each compartment points to a unique first digit. These compartments are: [0, 0.301), [0.301, 0.477), [0.477, 0.602), [0.602, 0.699), [0.699, 0.778), [0.778, 0.845), [0.845, 0.903), [0.903, 0.954), and [0.954, 1.000).

For example, for 14th period, the two series generate the distinct mantissa values of 0.191 and 0.194, yet both values are firmly inside the [0, 0.301) compartment belonging to digit 1. Figure E of Part 1 illustrates clearly the slack nature of the mapping of mantissa values into first digits. Yet, since these tiny differences in mantissa values are cumulative, eventually, when many more periods are considered, the 93.15% series begins to generate quite different mantissa values as compared with the 93.069773% series, and consequently its overall digit distribution becomes closer and closer to the Benford configuration after many more additional periods.

For any irrational normal series such as say 40.0% growth rate where LOG(1 + 40/100) = 0.14612803567823238…etc., supposedly one could find a pair of L and T integers such that approximately but very nearly 0.146128035678238 ≈ L/T. A finite irrational growth series in the neighborhood of another rational series having strong deviation from Benford experiences a great deal of deviation itself, and it requires a large number of periods in order to manifest its true and innate Benfordian character. In practical terms, such a cure is not available or not realistic when real-life growth series are concerned.

## [4] Empirical Tests of Digital Configurations of Exponential Growth Series

Do rational series in close proximity to irrational series really adversely affect logarithmic behavior for many or even most finite series? Fortunately they do not! They hardly ever disrupt digital configuration! And this fact can be verified empirically with the aid of the computer. The line of reasoning above, worrying about neighboring rational series adversely affecting many or most of the irrationals series turned out to be mere 'theoretical paranoia'. Let us empirically test a large variety of finite exponential growth series [with length from 822 growth periods to 3000 growth periods] for their compliance or non-compliance with Benford's Law.

Three distinct random checks shall be performed: (1) for rates randomly chosen from 1% to 5%, (2) for rates randomly chosen from 1% to 50%, (3) for rates randomly chosen from 1% to 890%. The rates are chosen with equal probability as in the continuous Uniform(1, 5), Uniform(1, 50), and Uniform(1, 890), respectively. Each randomly chosen growth rate is then utilized to create several exponential growth series with that rate applying a variety of growth periods, from the shortest one of only 822 periods, to the longest one of 3000 periods; all beginning at the same initial base value of 3. The need to create for each random rate several series with distinct growth periods emanates from the general requirement for an approximate integral exponent difference between the first and the last elements in obtaining logarithmic behavior, namely that LOG(maximum) – LOG(minimum) ≈ an integer, as in the k/x distribution.

The computer program then chooses that series with the lowest SSD, and which is often equivalent to choosing the series having exponent difference closest to an integer.
The variety of growth periods that are being employed for each randomly chosen growth rate are: {3000, 2897, 2800, 2697, 2597, 2297, 2284, 2262, 1930, 1759, 1433, 1268, 822}, and the program then registers only the minimum SSD of all these thirteen series, so as to arrive at the most Benford-complying series, and by implication often at that series where LOG(maximum/minimum) is closest to an integer. Such a 'crude' method of avoiding [in the approximate] non-integral exponent difference is not effective for extremely low growth rates between 0% and 1%, and this fact provided the motivation to avoid this challenging range on the percent-axis, and to consider only rates above 1%. Extremely low rates below 1% grow very slowly, and therefore such series often take very many periods (>>3000) to obtain integral exponent difference. Any computer program attempting to determine whether such low rate series are Benford would need to be much more sophisticated and to allow numerous periods.

Let us first outline that part of the empirical results coming out of the computer calculations regarding the supposed damage done to irrational series due to existing in the neighborhood of rational series.

As a dramatic example from the empirical computerized results thus obtained, we focus on the rebellious percent range around the theoretical rationality LOG(F) = 2/3 with its associated growth rate of 364.1588%. This is a highly rebellious growth rate with T = 3, having the exceedingly high SSD value of about 1500, indicating severe deviation from Benford. As can be seen in Figure 10, such a highly disruptive rate affects wide range all around it. Surely, affected range is such only due to the fact that the irrational series around it have no more than 3000 periods as programmed in the computer scheme, and surely substantially longer lengths of, say, one million periods would remedy almost all of them except those residing extremely near 364.1588% [and even those could also be remedied via billions or trillions growth periods]. The closer is the irrational to the disruptive rational, the more periods are needed in order to overcome the bad influence of the rebellious neighboring series, and this is the reason for the approximate bell-shape round curve seen in Figure 10.

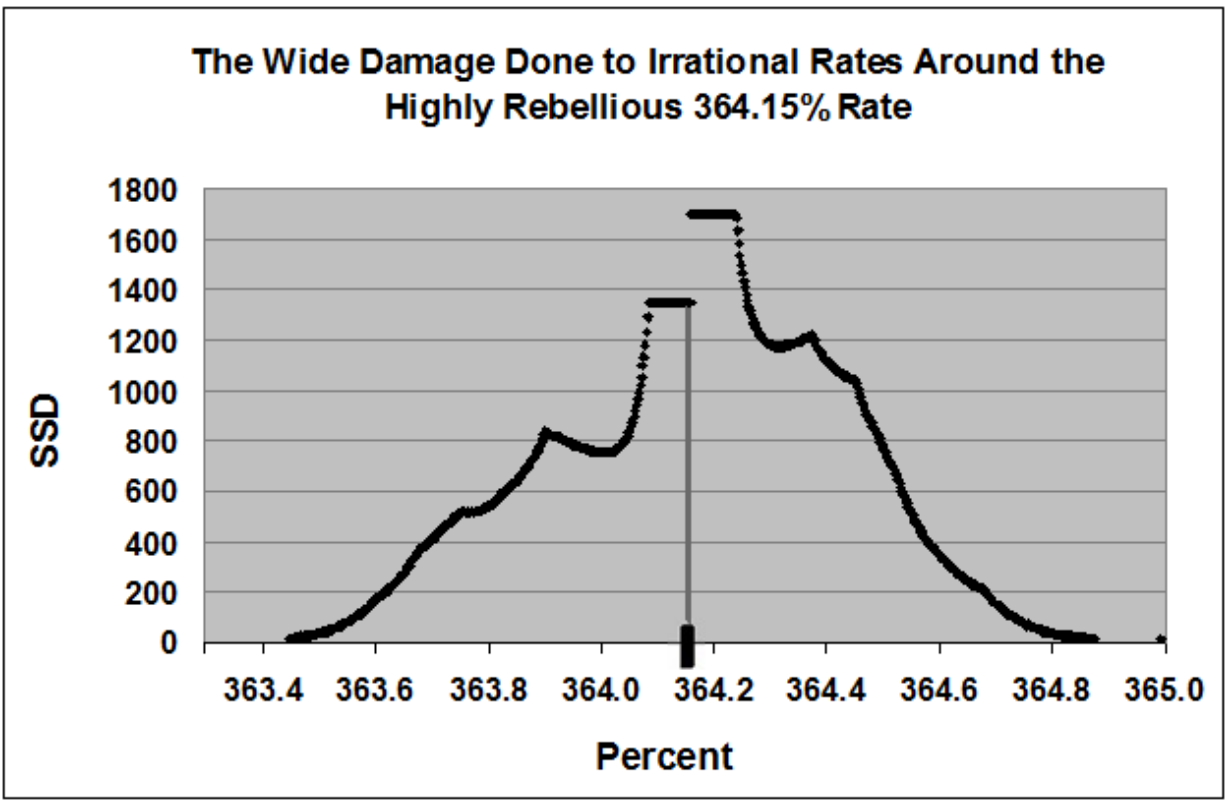

**Figure 10:** Wide Range is Affected Due to this Rational Series, Approximately 1.4% Units

Note: the percent growth rate is calculated via $2/3 = LOG(1 + Percent/100)$, then taking 10 to the power of both sides of the equation we get: $10^{2/3} = 10^{LOG(1 + Percent/100)}$ therefore $10^{2/3} = (1 + Percent/100)$, and $10^{2/3} - 1 = Percent/100$, so that $100[10^{2/3} - 1] = 364.1588\%$.

The lower the T value of the rebellious rate, the more significant is the deviation from Benford, and the wider is affected range. The higher the T value, the milder is the deviation from Benford, and the narrower is affected range. Hence, in sharp contrast to the example of Figure 10 with its affected wide range of about 1.4% units, another example shall be explored focusing on the rebellious range round the theoretical rationality $LOG(F) = 13/19$ and its associated growth rate of 383.2930%. This series has the relatively large T value of 19, thus its deviation from Benford is mild. The empirical result obtained here shows a much narrower affected range for this rate, as well as much reduced intensity of deviation, with SSD consistently below the value of 76. As can be seen in Figure 11 which depicts the empirical results around 383.2930%, such mildly disruptive rate affects only a very narrow range around it of about 0.15% units.

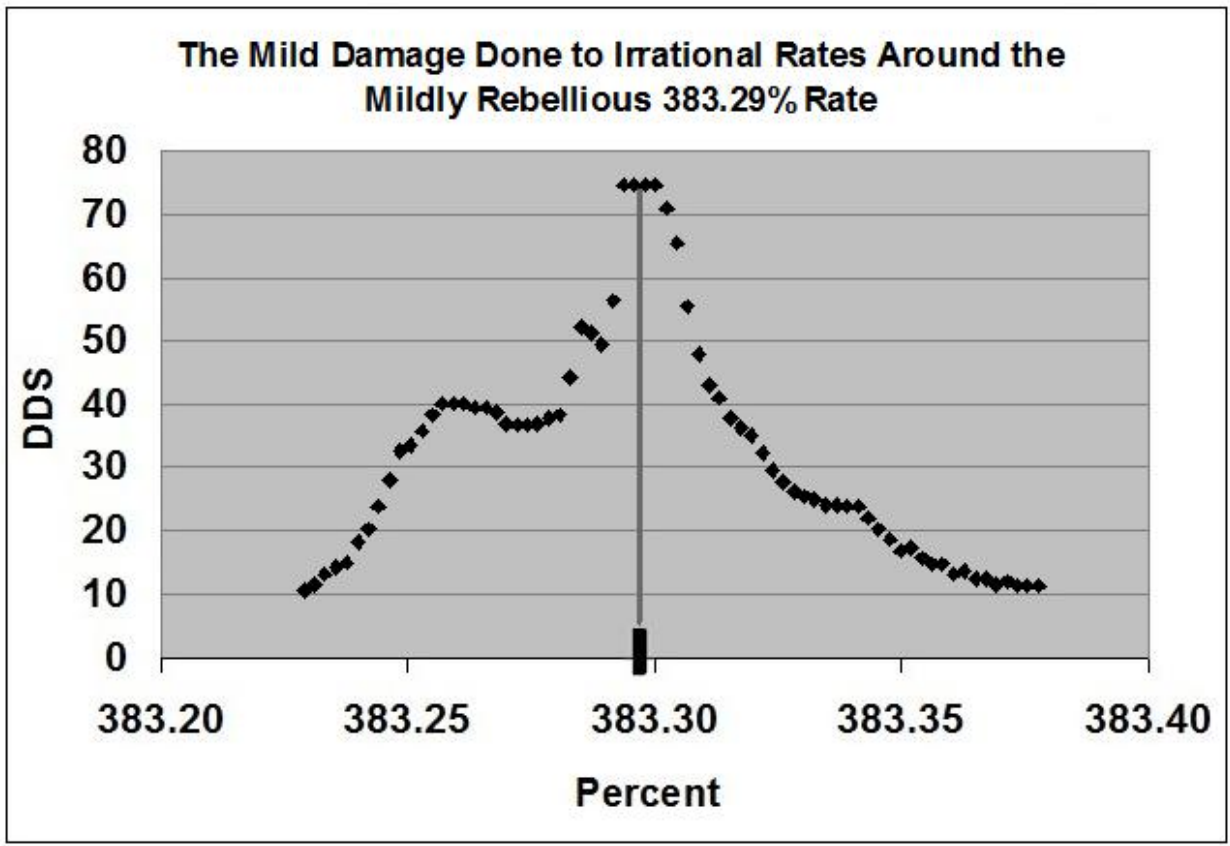

**Figure 11:** Narrow Range is Affected Due to this Rational Series, Only About 0.15% Units

Note: the percent growth rate is calculated via $13/19 = LOG(1 + Percent/100)$, then taking 10 to the power of both sides of the equation we get: $10^{13/19} = 10^{LOG(1 + Percent/100)}$ therefore $10^{13/19} = (1 + Percent/100)$, and $10^{13/19} - 1 = Percent/100$, so that $100[10^{13/19} - 1] = 383.2930\%$.

The [false] alarm is now raised, and the worry that many or most [finite] exponential growth series are disrupted by neighboring rational series now seems all too convincing. The milder case of Figure 11 gives some hope for better logarithmic results, while the severe case of Figure 10 is a bad omen for logarithmic behavior.

Yet in reality, the following thorough empirical tests encompassing growth rates from 1% to 5%, from 1% to 50%, and from 1% to 890%, decisively show that there is not much to worry about, and that almost all exponential growth series are Benford or very nearly so – assuming that these series are of considerable length, say over 1000 periods.

For 30,000 simulations of such random growth rates between 1% and 5%, only two came out with SSD over the mild value of 10! These two mildly rebellious series were 4.81134164367% with 13.9 SSD, and 4.81139819992% with 10.4 SSD; both located at the upper range near the 5% boundary, and both belong to the L/T = 1/49 theoretical rationality. Other rationality-related growth rates found there with even lower deviations were: 4.71269871352% with only 4.11 SSD; also located at the upper range near the 5% boundary, and belonging to the L/T = 1/50 theoretical rationality, as well as 4.91365456231% with only 4.75 SSD; also located at the upper range near the 5% boundary, and belonging to the L/T = 1/48 theoretical rationality. Average SSD for all of these 30,000 low growth rates is the highly satisfactory and extremely low value of 0.12.  The median gives an even more satisfactory and incredibly low value of 0.03.

For 30,000 simulations of random growth rates between 1% and 50%, we get the following results: Out of 30,000 rates, only 278 (namely 0.9%) came out with SSD over the mild value of 10. Out of 30,000 rates, only 101 (namely 0.3%) came out with SSD over the relatively high value of 50. Out of 30,000 rates, only 55 (namely 0.2%) came out with SSD over the high value of 100. Average SSD for all of these 30,000 growth rates is the satisfactory very low value of 0.9.  The median (being a much more robust measure of centrality) gives the extremely low value of 0.02. As known from Descriptive Statistics, the average usually exaggerates and overestimates centrality, while the median could at times underestimates centrality.

For 30,000 simulations of random growth rates between 1% and 890%, we get the following results: Out of 30,000 rates, 2,515 (namely 8.4%) came out with SSD over the mild value of 10. Out of 30,000 rates, 1070 (namely 2.5%) came out with SSD over the relatively high value of 50. Out of 30,000 rates, 455 (namely 1.5%) came out with SSD over the high value of 100. Average SSD for all of these 30,000 growth rates is the moderate value of 10.0. The median gives the very low value of 0.6. The median in this case might be underestimating centrality.

**Conclusion:** For finite series of considerable length with at least 822 periods up to 3000 periods: not only most, but the vast majority of exponential growth series are nearly perfectly Benford! Moreover, LOW rates with growth below approximately 5% are practically all nearly perfectly Benford with the exceptions of only very few and extremely rare mildly rebellious cases. For HIGH rates with growth approximately over 5% there exists a small non-Benford minority, although the vast majority of them do satisfactorily obey the law of Benford. [Note: the choice of 3 as the base value for all of these empirical tests is practically irrelevant to this conclusion.]

## [5]  The Proper Explanation of the Rarity of Anomalous Series

Let us now give the proper and more profound explanation of why the vast majority of exponential growth series are rarely rebellious and do obey the law of Benford nearly perfectly:

The reason why almost all finite exponential growth series with considerable length [of over 1000 periods say] are Benford is because of a pair of limitations placed on the values of T and L.

Limitation on the value of T:

Deviation due to coordinated and rational walk along the log-axis is limited to series with relatively low T values, while series with high T values are very close to Benford approximately in spite of their rationality. As seen in Figure 6, whenever T is roughly over 50, deviation from the logarithmic is fairly small. If one is willing to tolerate SSD values up to say approximately 10 or 15 - considering all series with SSD value lower than say 10 or 15 as sufficiently near Benford - then any rational series with T over 50 should not be considered as truly rebellious, but rather as approximately Benford.

Hence T is limited as in: $\mathbf{T \leq 50}$.

Limitation on the value of L:

For any fixed interval on the percent-axis, as in the exclusive consideration of the range $(0\%, P_{MAX}\%)$, and the exploration of all the rebellious rates that might be laying there, T and L values are restricted as in $L/T \leq LOG(1 + P_{MAX}/100)$.

Hence L is limited as in: $\mathbf{L \leq T*LOG(1 + P_{MAX}/100)}$.

For any possible T value, the value of L is restricted by this ceiling on the upper side.

This pair of severe limitations on T and L values places a firm cap on how many rational L/T series with significant deviation from Benford exist on any given percent-axis interval. And it renders the set of all rational series with significant deviation from Benford not merely finite, but also very small in size! Moreover, as shall be demonstrated in the next two numerical examples, this pair of limitations on T and L values implies that the vast majority of T and L rational combinations are of relatively high T values around 30 to 50, relating to mild deviations from Benford, and that only a minority of such rational combinations are of relatively low T values, relating to significant deviations from Benford.

In some limited sense, high T values correspond with low growth, and low T values correspond with high growth – but this is so only for L being fixed at a constant value. Yet, the true nature of T value is not about low growth rate versus high growth rate, but rather about the number of periods the rebellious series go through in order to return to the same mantissa value. In the same vein, and in some limited sense, low L values correspond with low growth, and high L values correspond with high growth – but this is so only for T being fixed at a constant value.

As a consequence of the above discussion, it is noted that the T = 50 and L = 1 case yields the rebellious series with the lowest possible growth rate [under this pair of limitations on T and L values], and that lowest possible rebellious growth rate is 4.71285480509%. For this series, LOG(1 + 4.71285480509/100) = 0.02 = 1/50.

For all low growth rates below 4.71285480509%, deviation is very small, and thus they are considered to be nearly Benford. This fact alone saves all growth rates within the (0%, 4.7128%) percent-axis interval from the fate of being anomalous [assuming that SSD values below 10 or 15 are considered mild and that therefore such series are considered to be practically Benford].

Under this pair of limitations on T and L values, all combinations of T and L values other than the T = 50 and L = 1 case, yield growth rates higher than 4.71285480509%. To prove this assertion one needs only to note that all other cases of T = 50 where L > 1 such as 2/50, 3/50, 4/50, and so forth, represent higher growth rates than that of the 1/50 case. In addition, all other cases of L = 1 where T < 50, growth rates are higher as well. For example, 1/49, 1/48, 1/47 are all of higher growth rate than that of the 1/50 case. The same argument applies for example to 2/49, 3/48, 4/47, 5/46 and so forth, which are all of higher growth rate than that of the 1/50 case.

In a nutshell:

**(1)/(50) ≤ (any L ≥ 1)/(any T ≤ 50)**

and this is so for any L and T combination.

## [6] A Brief Example of the Consequences of L and T Limitations

The first brief example regarding the very limited number of pairs of T and L possibilities which are effectively disruptive to Benord behavior is given for the exclusive consideration of the percent-axis interval on **(0%, 100%)**. Here T and L values are restricted as follows:

$T \le 50$ and $L/T \le LOG(1 + 100/100)$
$T \le 50$ and $L \le T*LOG(2)$
$T \le 50$ and $L \le T*0.30103$

The set of all possible T and L disruptive pairs is then:

$T = 50$ and $L \in \{1, 2, 3, \ldots 12, 13, 14, 15\}$
because $L \le 50*0.30103$, namely $L \le 15.05$

$T = 49$ and $L \in \{1, 2, 3, \ldots 12, 13, 14\}$
because $L \le 49*0.30103$, namely $L \le 14.75$

$T = 48$ and $L \in \{1, 2, 3, \ldots 12, 13, 14\}$
because $L \le 48*0.30103$, namely $L \le 14.45$

$T = 47$ and $L \in \{1, 2, 3, \ldots 12, 13, 14\}$
because $L \le 47*0.30103$, namely $L \le 14.15$

$T = 46$ and $L \in \{1, 2, 3, \ldots 12, 13\}$
because $L \le 46*0.30103$, namely $L \le 13.85$

*and so forth*

$T = 14$ and $L \in \{1, 2, 3, 4\}$
because $L \le 14*0.30103$, namely $L \le 4.21$

$T = 13$ and $L \in \{1, 2, 3\}$
because $L \le 13*0.30103$, namely $L \le 3.91$

*and so forth*

$T = 7$ and $L \in \{1, 2\}$
because $L \le 7*0.30103$, namely $L \le 2.11$

*and so forth*

$T = 4$ and $L \in \{1\}$
because $L \le 4*0.30103$, namely $L \le 1.20$

T cannot fall below 4. For example, T cannot be 3, because then $L \leq 3*0.30103$, and this implies that $L \leq 0.90$, but this is a contradiction since L must be a positive integer, namely that at a minimum L must be integer 1. Hence the possibilities of $T = 3$, $T = 2$, $T = 1$ are eliminated here.

Therefore, the complete set of all possible {T, L} disruptive pairs is:

{50, 1}, {50, 2}, {50, 3}, … {50, 13},{50, 14},{50, 15}
{49, 1}, {49, 2}, {49, 3}, … {49, 13},{49, 14}
{48, 1}, {48, 2}, {48, 3}, … {48, 13},{48, 14}
{47, 1}, {47, 2}, {47, 3}, … {47, 13},{47, 14}
{46, 1}, {46, 2}, {46, 3}, … {46, 13}

  … etc. …

{14, 1}, {14, 2}, {14, 3}, {14, 4}
{13, 1}, {13, 2}, {13, 3}
{12, 1}, {12, 2}, {12, 3}
{11, 1}, {11, 2}, {11, 3}
{10, 1}, {10, 2}, {10, 3}
{9, 1}, {9, 2}
{8, 1}, {8, 2}
{7, 1}, {7, 2}
{6, 1}
{5, 1}
{4, 1}

This yields the very small set of only 360 disruptive pairs. Yet, there are many pairs that are not in reduced form, and they should be eliminated.

For example:

{50, 2}   is really {25, 1}       2/50 = 1/25
{50, 15} is really {10, 3}       15/50 = 3/10
{12, 2}   is really   {6, 1}       2/12 = 1/6

All of which yields an even smaller set of only 232 {T, L} disruptive pairs in purely reduced form on the percent-axis interval (0%, 100%).

# [7]  A Detailed Example of the Consequences of L and T Limitations

A second example fully analyzed and detailed - regarding the very limited number of T and L pair possibilities which are effectively disruptive to Benford behavior - is given for the exclusive consideration of the percent-axis interval on **(0%, 900%)**. Here T and L values are restricted as follows:

T ≤ 50   and   L/T ≤ LOG(1 + 900/100)
T ≤ 50   and   L ≤ T*LOG(10)
T ≤ 50   and   L ≤ T*1
T ≤ 50   and   L ≤ T

Therefore, the complete set of all possible {T, L} disruptive pairs is:

{50, 50}, {50, 49}, {50, 48},  … {50, 3}, {50, 2}, {50, 1}
{49, 49}, {49, 48}, {49, 47},  … {49, 2}, {49, 1}
{48, 48}, {48, 47}, {48, 46},  … {48, 1}

  … etc.  …

{3, 3}, {3, 2}, {3, 1}
{2, 2}, {2, 1}
{1, 1}

This yields the small set of only $(N^2 + N)/2 = (50^2 + 50)/2 = 1275$ pairs.

It should be noted that pairs with low T values [high intensity deviation] are by far fewer than pairs with high T values [low intensity deviation], therefore most of these pairs cause relatively mild deviation from Benford and only on relatively narrow and limited neighboring ranges. For example, there are only $(N^2 + N)/2 = (10^2 + 10)/2 = 55$ pairs with low T values from 1 to 10 causing intense and wide damage, while the vast majority of them, namely $1275 - 55 = 1220$ pairs, are with T values of over 10, causing relatively milder and narrower damage. In other words, the upside down pyramid-like structure of the above set of all possible {T, L} pairs implies that most of the pairs are around the top near T = 50 case; that fewer are around the center near the T = 25 case; and that only a tiny minority are around the bottom near the T = 1 case. Such a structure for the set of {T, L} pairs bodes well for Benford compliance in general, because low T values where most significant deviations occur are much fewer and rarer than high T values where only insignificant and minor deviations occur.

The above set of 1275 pairs contains numerous redundant cases. This is so since many pairs are not in reduced form and should be eliminated. This purge yields an even smaller set of only 774 reduced form disruptive {T, L} pairs. For example, the first pair {50, 50} is really the last pair {1, 1}; the pair {50, 2} is really {25, 1}, and so forth.

Hence only 774 explosive mines are laid out there quite sparsely on the large minefield of all the real numbers from 0% to 900%. These 774 land mines are literally few and far between, and they are rarely encountered. Most of the mines are of distinct sizes and gunpowder concentration (namely distinct T values), hence usually they cause different types of damage when activated, killing or maiming to a large or small degree. Is it really very dangerous for a normal exponential growth series to walk over this relatively sparse minefield? Well, for finite and relatively short series without many growth periods, the answer depends on the width of the damaged range around each land mine. For infinitely long series with infinite number of growth periods, walking over this minefield is quite safe, because affected range around each land mine approaches zero in the limit as the length of the series goes to infinity, rendering these mines as a non-threatening set of 774 real points on the infinitely uncountable real range (0%, 900%).

Let us contemplate the danger of walking over this minefield for finite series of considerable length with 1000 to 3000 periods approximately.

If we assume that all of these 774 mines are of the low intensity type as seen in Figure 11 relating to the 13/19 theoretical rationality where only a narrow range of about 0.15 units was affected, then simple-minded calculations lead to the conclusion that this would supposedly cover only 0.15*774 = 116.1 units on the percent-axis range, out of a total of 900 units, assuming that there is no overlapping of affected ranges. This still represents a significant portion of affected ranges, namely 13% of the total range.

If we assume that all of these 774 mines are of the high intensity type as seen in Figure 10 relating to the 2/3 theoretical rationality where a wider range of about 1.4 units was affected, then simple-minded calculations lead to the conclusion that this would supposedly cover a whapping 1.4*774 = 1083.6 units on the percent-axis range, out of a total of 900 units, assuming that there is no overlapping of affected ranges. This represents more than the entire possible range of 0% to 900%! The implication of this scenario would be that all finite exponential growth series with about 1000 to 3000 periods are rebellious!

Yet, the empirical results decisively showed that nearly all finite series of considerable length with 1000 to 3000 periods approximately are Benford! This implies that the fear of the impact due to rational series is exaggerated, and that the damage done on the percent-axis to irrational series is overestimated.

In a nutshell, the reason for the very mild and very limited damage done on the entire range by these 774 land mines is that the vast majority of T values here are of relatively high value of around 25 to 50. These high-T-value series cause very mild and very narrow damage, and progressively so as one nears the 50 value mark. The scenario of Figure 10 relates to T = 3, and this is rather a very rare case within the set of 774 land mines. In addition, the scenario of Figure 11 relates of T = 19, and even this case is not of the most frequent types within the set of 774 land mines. The shape of the entire set of {T, L} disruptive pairs shown in the previous page is that of an upside down pyramid-like structure, therefore, overall, out of 774 cases, those with relatively low T values causing significant and wide damage constitute only a very small minority, representing truly only very few cases. All this renders the (0%, 900%) range very much Benford-like.

Let us define the theoretical magnitude of deviation from Benford as **500/T**. The motivation behind this definition emanates from the fact that the larger the value of T the less deviation from Benford is found for the rebellious rational series, namely that deviation is inversely proportional to T. Since L practically does not play any role in the magnitude of deviation, it is left out of this definition.

The smallest possible magnitude of theoretical deviation is of the T = 50 case, and it is calculated as 500/50, namely 10. The largest possible magnitude of theoretical deviation is of the T = 1 case, and it is calculated as 500/1, namely 500, yet this case is discarded and excluded, and only deviations from the 2nd largest T = 2 case are to be considered. This top value of theoretical deviation is calculated as 500/2, namely 250.

The reason for this exclusion is that T = 1 case relates to the highest possible growth rate where P = 900%, and perhaps it is wise to avoid the inclusion of the boundary of the percent-axis interval of (0%, 900%). [Since L ≤ T it follows that L ≤ 1, so that L must be 1 as well, therefore L/T = 1/1 = LOG(1 + P/100), so that 10 = (1 + P/100), and finally we obtain P = 900%.] The other reason for this [singular] exclusion is to be able to give a good visual presentation of all the other [773] less dramatic cases where 500/T ≤ 250, as shown clearly in Figure 12. This unique case of T = 1 where 500/T is 500 can be considered as an outlier in some sense, and including it would needlessly distort and expand the scale for all the other cases in Figure 12.

Figure 12 depicts these 773 (namely 774 – 1) theoretical rebellious rates, where the vertical axis expressing 500/T is plotted versus the horizontal axis expressing the percent growth. The count of the points in the entire chart of Figure 12 where 500/T values are 50 or higher shows that there are only 31 such points out of a total of 773 points, namely that only 4% of all the points are with truly strong deviation from Benford. Clearly, the vast majority of the points in Figure 12 are with low 500/T < 50 deviation values.

In addition, it should be noted that the overall frequency of the points (crude 'density' of sorts) falls off as focus shifts from low percent growth to higher percent growth, getting more diluted on the right. For example, from 0% up to 100% there are 232 points, from 100% up to 200% there are 138 points, and from 800% up to 900% there are only 32 points. The vector of frequencies per 100% interval is {232, 138, 94, 77, 62, 52, 45, 41, 32}, falling off monotonically to the right. In contrast, average magnitude of deviation - namely average of 500/T per 100% interval - is nearly steady, and it's approximately of the constant value 19.5.

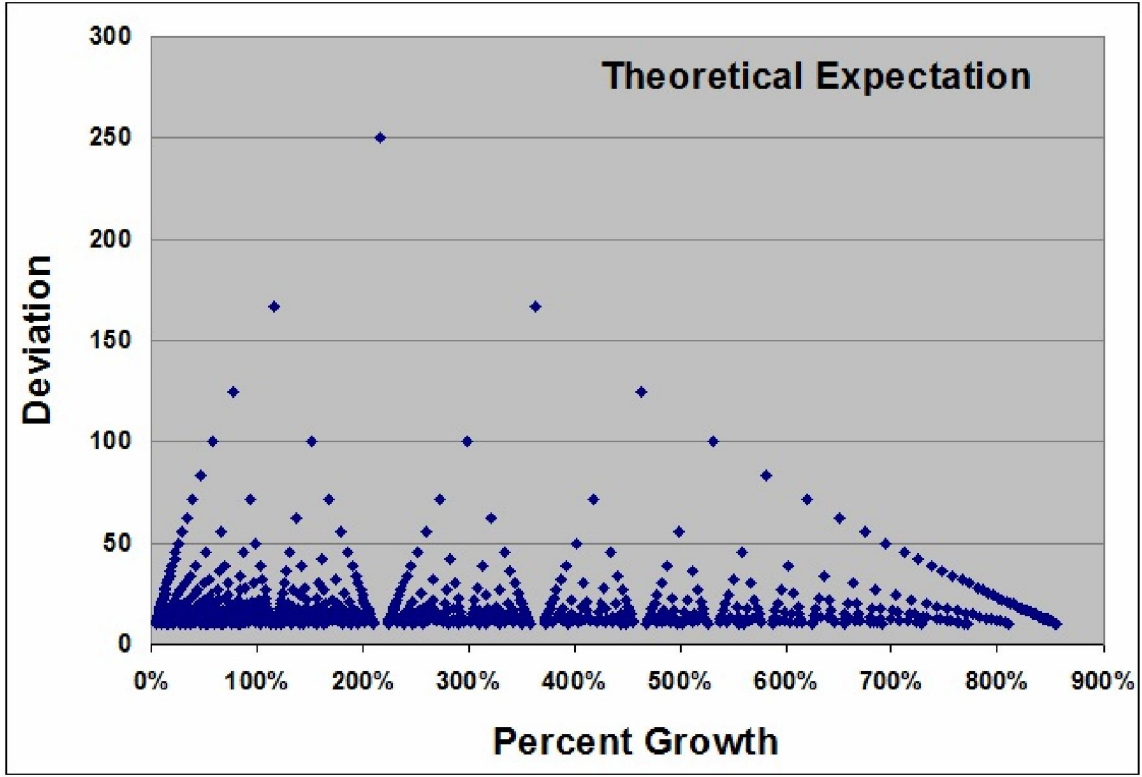

**Figure 12:** Theoretical Expectation of Rebellious Rational Rates from 0% up to 900%

## [8]  Excellent Fit between Empirical Results and Theoretical Expectation

Consequently, another computer program is run, empirically checking exponential growth series for any possible deviation from Benford. Here for this 2nd empirical check, instead of choosing growth rates randomly as was done previously in the 1st empirical check, this program methodically examines growth rates in order, from 1% all the way to 890%, advancing forward in tiny % increments. This is so in order to cover for sure and without any random glitches all corners and all sub-intervals within the entire interval of (0%, 900%) – except for the problematic part at the end of (890%, 900%). The aim is to show excellent correspondence between what is expected theoretically due to any possible rational march along the log-axis and what is actually found empirically. If such good fit is found, then as a consequence, this could indirectly imply that no other causes or reasons for deviation from Benford are observed for exponential growth series in extreme generality except for the argument discussed in this article [assuming that sufficient number of elements of the series are considered for good compliance, as well as adherence to the constraint of an integral exponent difference for low rates].

The program starts at 1% and ends at 890%, in tiny refine increments of 0.00215714598%. It checks 412,118 rates in total, namely nearly half a million exponential growth series. It uses the quantity 3 as the initial base value for all the series. The program only registers those series with deviation larger than the somewhat arbitrary value of 8.88 SSD. Any obedient series with SSD less than 8.88 is deleted and forgotten. The motivation for the choice of the cutoff SSD value of 8.88 is to aim just slightly below the round number of 10, in an attempt to correspond to the theoretical list of rebellious rational rates with a cap of T at 50. That cap at 50 for the value of T restricts SSD values for all the theoretical rational series to over 10 or over 15 approximately.

Just as was done for the random selection of growth rates, here too, in order to allow low rates to comply with the logarithmic requirement that exponent difference between the maximum and the minimum should be approximately close to an integer, a variety of growth periods are used via the set {3000, 2897, 2800, 2697, 2597, 2297, 2284, 2262, 1930, 1759, 1433, 1268, 822}. This is the same set of 13 growth periods that was used in the 1st program. The program then registers only the minimum SSD for all these series, so as to arrive at the most complying series.

The program has registered 37,453 rebellious growth rates with SSD over 8.88 – out of 412,118 in total. In other words, about 9.1% of the series were found to deviate by registering over the 8.88 SSD threshold value. Figure 13 depicts those empirical rebellious rates with SSD over 8.88 on the interval (1%, 890%). It should be emphasized that there are various gaps and empty sub-ranges here for the obedient series, in spite of the fact that the chart appears continuous.

**Note:** The earlier simulations in chapter 4 of random growth rates between 1% and 890% resulted in having 8.4% of all the series with SSD over 10. In comparison, here the SSD threshold is lowered from 10 to 8.88, allowing a few more series between 8.88 and 10 to be considered rebellious, and this is why slightly higher 9.1% ratio is obtained here instead of 8.4%.

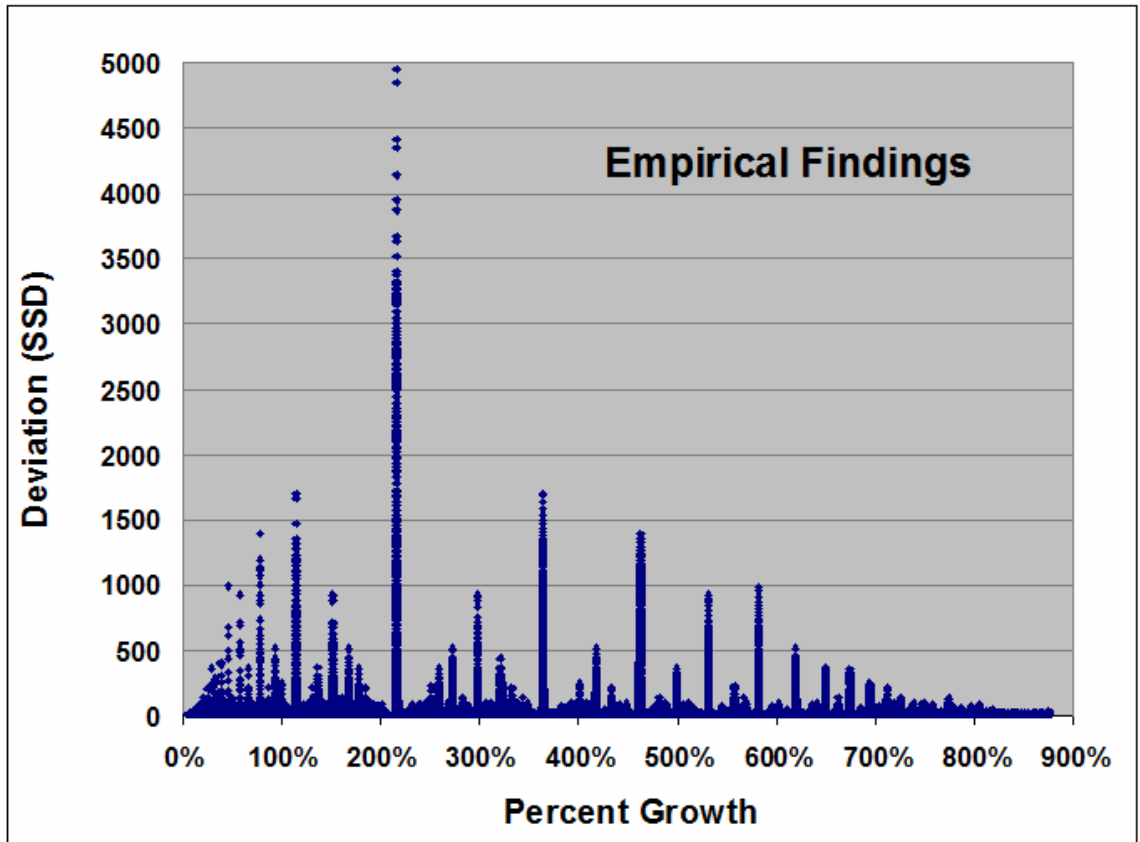

**Figure 13:** Empirical Findings of Rebellious Exponential Growth Series from 1% to 890%

The moment of truth has arrived. Let us now check compatibility between theoretical and empirical results by putting these two charts next to each other as shown in Figure 14.

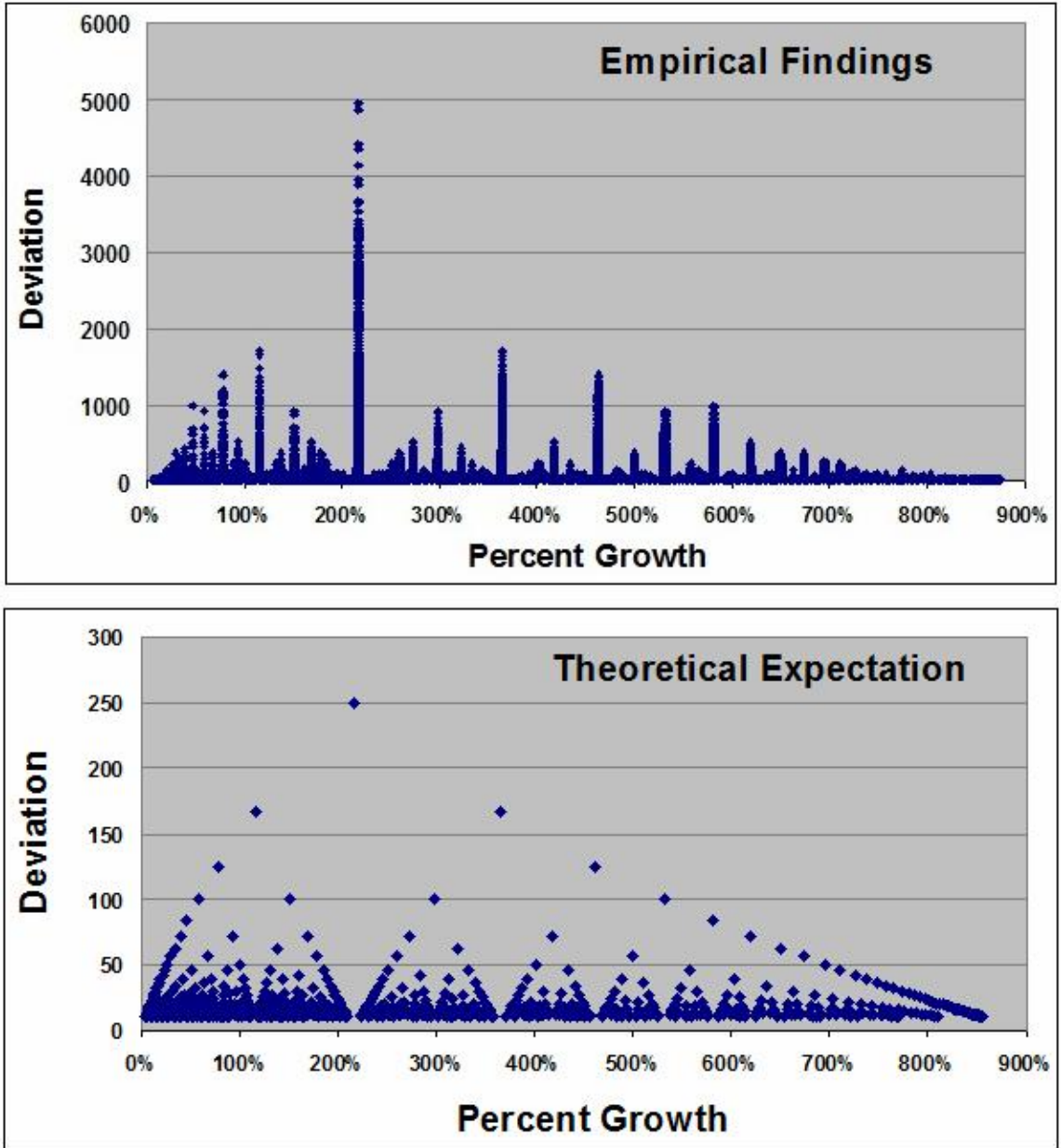

**Figure 14:** Compatibility between Theoretical & Empirical Rebellious Series 1% to 890%

Hooray! The empirical and the theoretical agree!

Furthermore, a closer scrutiny of various smaller sub-intervals within the entire 1% to 890% range reveals a near perfect fit between the theoretical and the empirical. As for but one example of these careful examinations, Figure 15 depicts the excellent fit for rates between 1% and 50%. Comparisons of various other sub-intervals show the same good fit as well.

Figure 15 clearly reveals the various gaps and empty sub-ranges that exist for the obedient series in the empirical computerized study. This helps to repudiate the false appearance of continuity in Figure 13.

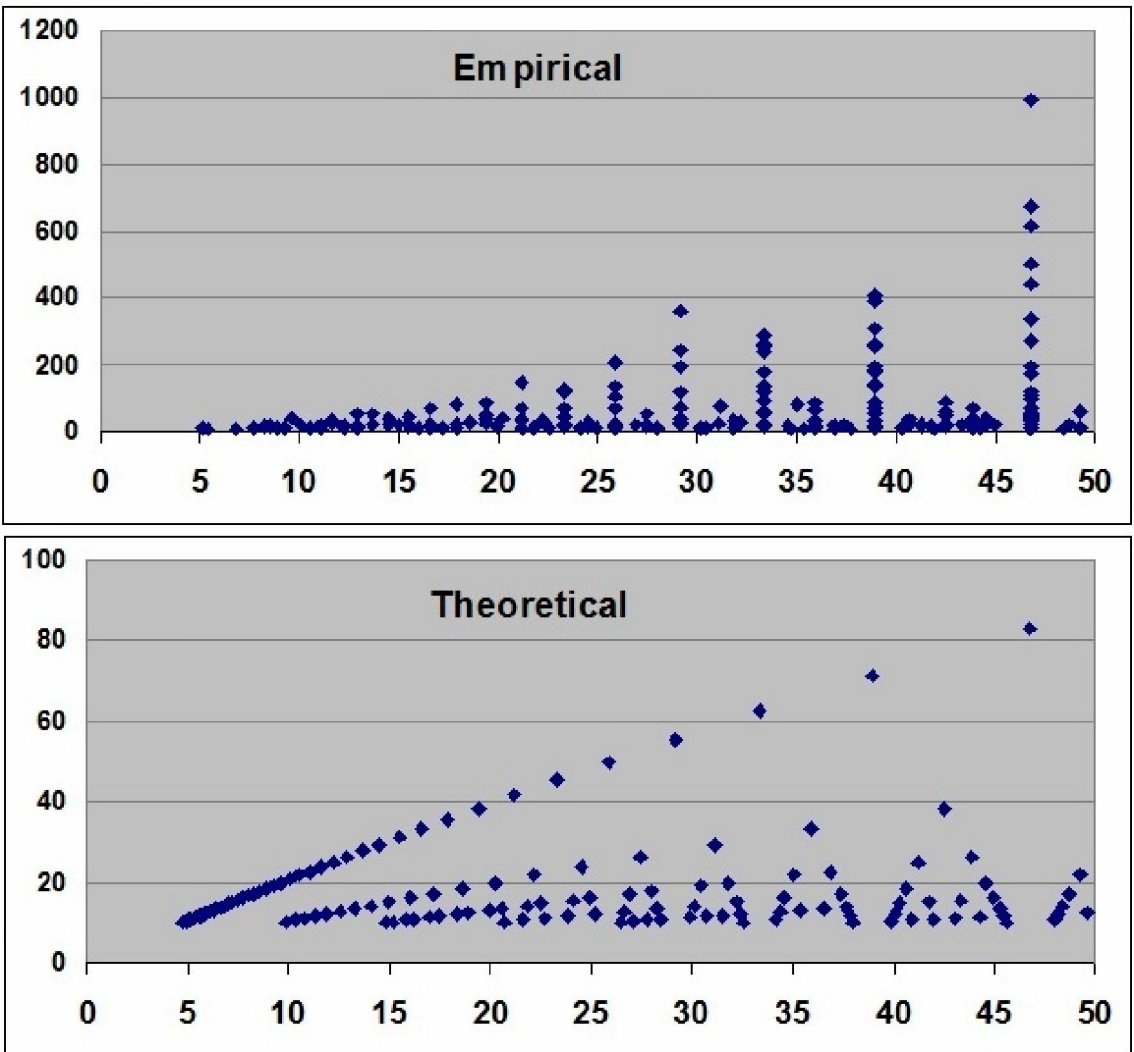

**Figure 15:** Compatibility between Theoretical and Empirical Rebellious Series 1% to 50%

Empirically, for the entire region below 5%, not a single rebellious rate was detected; as all growth rates from 0% to 5% had their SSD values below 8.88. This empirical result is almost in perfect harmony with the theoretical threshold of approximately 5% that was calculated earlier (actually 4.71285480509% to be precise), and which showed that there exist only 3 rational rebellious series just below 5% - under the assumption of the theoretical pair of limitations on T and L values. These 3 theoretical rational rates are: 4.71285480509%, corresponding to the theoretical rationality of 1/50; 4.81131341547%, corresponding to the theoretical rationality of 1/49; and 4.91397291363%, corresponding to the theoretical rationality of 1/48.

These 3 theoretical rational series have been missed out by the empirical computer program due to the relatively low SSD generated exactly at these 3 points, as well as around their neighborhood perhaps. Under the parametrical assumptions of the computer program, namely the base of 3 and with these 822 to 3000 periods, at exactly the theoretical series of 4.71285480509%, SSD is 8.77; at exactly the theoretical series of 4.81131341547%, SSD is 8.47; and at exactly the theoretical series of 4.91397291363%, SSD is 7.85. Since the SSD threshold of the computerized empirical scheme is set at 8.88, these theoretical rationalities must have been missed out by a hair.

Figure 16 depicts discrete-like comparisons between all 'large' empirical deviations with SSD values over 100, referring only to the central point of each such rebellious sub-range, together with all those 'large' theoretical deviations with T < 14, namely all those with theoretical magnitude of deviation 500/T > 38.5. These two cutoff points, namely T < 14 and SSD > 100, correspond to each other in the approximate, since SSD of the rationality 1/14 under the assumptions of the empirical parameters is 101.5, thus SSD value of 100 could serve perhaps as a nice round figure for the cutoff point between 'large deviations' and 'non-large deviations' of the empirical results. As can be seen in Figure 16, there is an excellent fit between the theoretical and the empirical for all these 'large' deviations from Benford.

The excellent fit found here between the empirical results and the theoretical expectation strongly confirms the conclusion that there exist no other causes or reasons for deviations of exponential growth series from the Benford configuration except for the rationality of their LOG(F). And even though the next concluding statement may sound highly puzzling and paradoxical to the psychologist, psychiatric, or social scientist, yet it should be believed, namely that: "Any **deviant** digital behavior of a growth series differing from the Benford configuration is always exclusively associated and correlated with the **rationality** of its related LOG(F) fraction". There seems to be nothing else adversely affecting digital behavior except for the argument given in this article regarding the rational march along the log-axis.

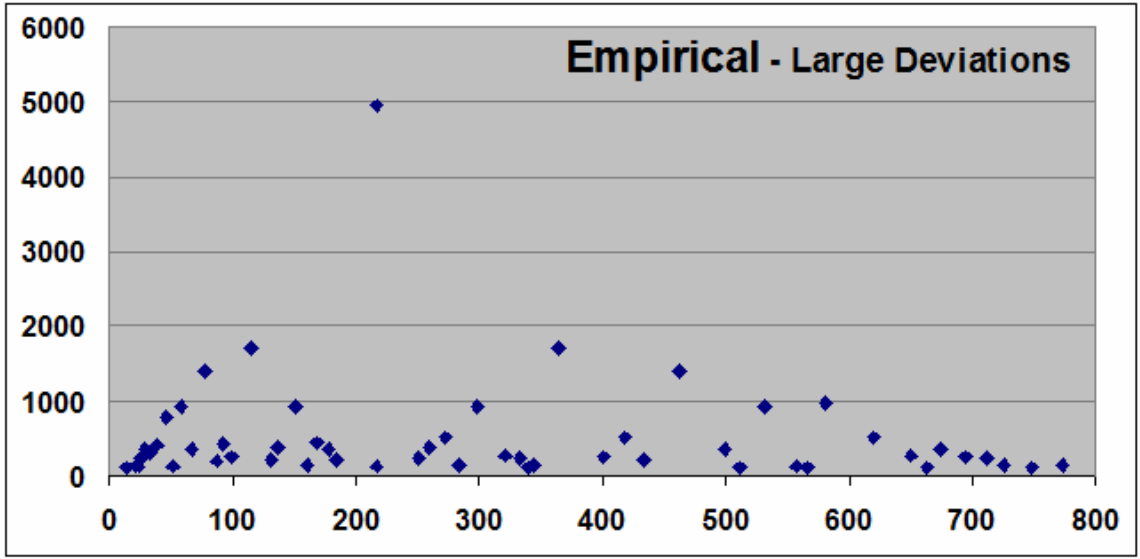

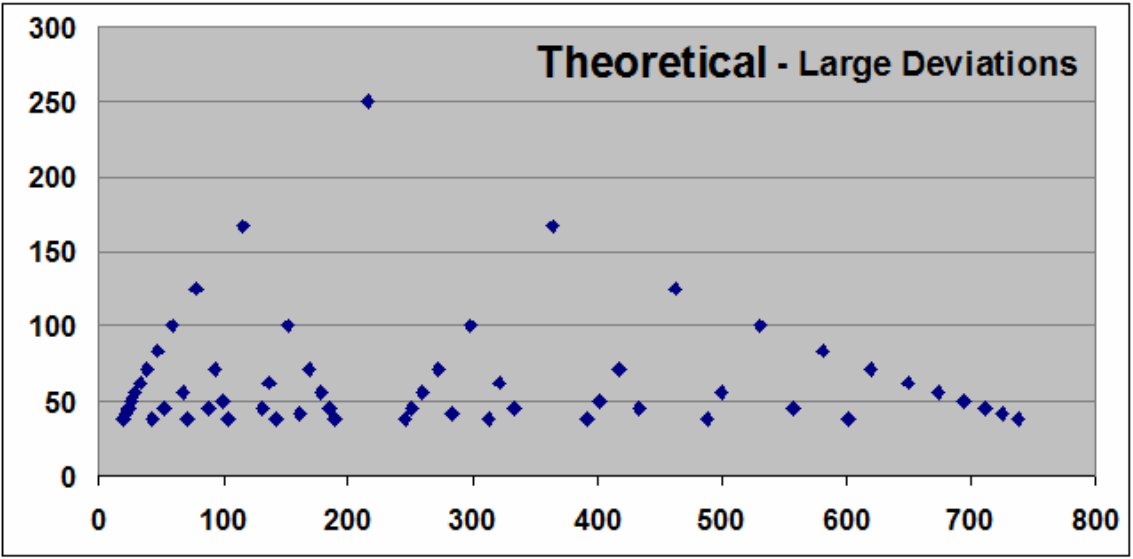

**Figure 16:** Compatibility between Theoretical & Empirical in Large Deviations

## [9]  Wider Damage Done to Short Normal Series by Anomalous Series

The excellent fit between the empirical and the theoretical seen above; and the very positive conclusion of chapters 4 and 8, namely that the vast majority of exponential growth series are nearly perfectly Benford; both assume that the series are of considerable length, with at least 822 periods up to 3000 periods.

For shorter series with length of only100 periods or so for example, the effects of the rational series are felt far and wide, leaving many more series with considerable deviation from Benford.

Consequently, in order to examine how very short series behave digitally, another computer program is run, empirically checking short series for deviation from Benford. The program starts at 2% and ends at 890%, in tiny refine increments of 0.0098375462%. It checks 90,266 exponential growth series in total. It uses the quantity 3 as the initial base value for all the series.

Just as was done in the previous two empirical programs, utilizing a variety of lengths, here too, in order to allow low rates to comply with the logarithmic requirement that exponent difference between the maximum and the minimum should be approximately close to an integer, a variety of growth periods are used via the set {101, 102, 103, 104, 105, 106, 107, 108, 109, 110}. The program then registers only the minimum SSD for all these series, so as to arrive at the most complying series.

The deliberate avoidance of the range (1%, 2%) is due to the fact that the set of the above choices of 10 lengths is still not sufficient for very low growth rates below 1%, and many more choices for the length are needed there in order to comply with the logarithmic requirement that exponent difference should be close to an integer.

The program has registered 27,554 growth rates with SSD over 10, out of the 90,266 total. In other words, about 30.5% of the series were found to deviate by registering over the 10 SSD threshold value. About a third of very short exponential growth series are not really Benford! And this highly pessimistic conclusion is all due to these 773 land mines of the rational series sparsely laying there on the long minefield from 2% to 890%!

Such high **30.5%** portion of affected [short] series with SSD over 10 is compared with only **8.4%** portion of affected [long] series with SSD over 10 as was seen earlier in chapter 4.

Figure 17 depicts the wider damage done to the entire neighborhood of the rational series 364.1588% growth rate - associated with the rationality LOG(F) = 2/3 – affecting many short irrational series all around it due to its presence. This chart is derived from the above computer program utilizing very short series of lengths varying from 100 to 110 periods, and starting from the initial base value of 3. The comparison between Figure 17 and Figure 10 reveals that short series suffer by far greater damage from the 364.1588% rebellious rational series. Here for short series with only 100 to 110 growth periods, much wider range of about 5.5% is adversely affected, as compared with only about 1.4% affected range of the longer series with 822 to 3000 growth periods.

The same conclusions and results are expected to be seen for the neighborhood of the rational series of 383.2930% growth rate - associated with the rationality LOG(F) = 13/19. As seen in Figure 11, only very narrow range of about 0.15% units is affected, but this is so due to the consideration of long finite series with 822 to 3000 growth periods. Much shorter series of about 100 growth periods should show a wider range of affected series than merely 0.15%.

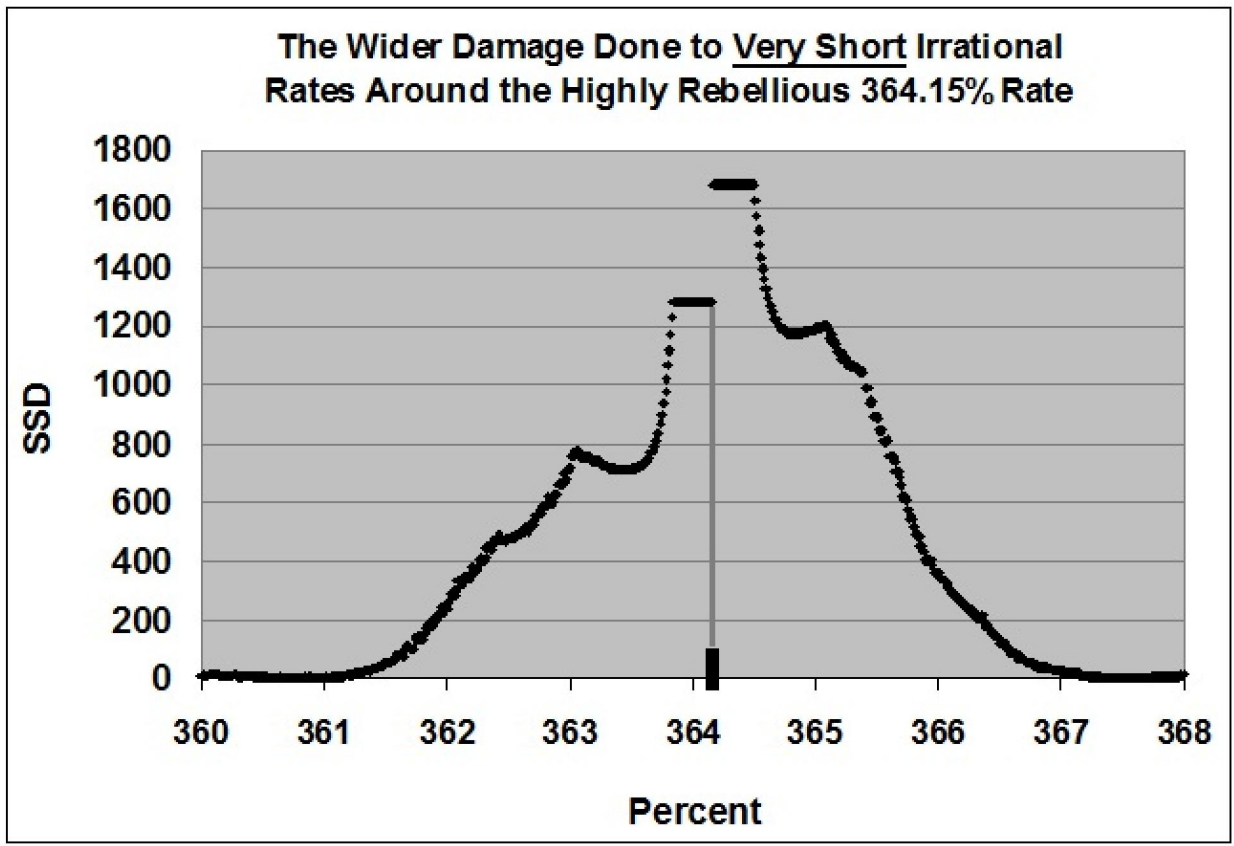

**Figure 17:** Wider Range is Affected for very Short Series, Approximately 5.5% Units

# [10]  Random Exponential Growth Series

Exponential growth series are constructed with a constant (fixed) multiplicative factor $F_{CONSTANT}$. Another possibility to consider is that of the random multiplicative factor $F_{RANDOM}$. For example, random factors constantly chosen from the Uniform on (1.23, 1.67) for each growth period could be utilized, leading to an exponential growth series with haphazard and constantly changing growth rate, randomly fluctuating between 23% and 67%. We shall denote $U_N$ as the Nth realization from the Uniform Distribution, standing for the multiplicative factor $F_{RANDOM}$.

$\{B, BF, BF^2, BF^3, \ldots , BF^N\}$ standard exponential growth series

$\{B, BU_1, BU_1U_2, BU_1U_2U_3, \ldots , BU_1U_2\,U_3{\ldots}U_N \}$ random exponential growth series

Such random selections of growth rates yield an overall uniform log distribution just as the fixed growth rate series usually do, assuming there are plenty of elements in the series. Let us examine how the series advances along the log-axis:

$\{LOG(B),$
  $LOG(B) + LOG(U_1),$
  $LOG(B) + LOG(U_1) + LOG(U_2),$
  $LOG(B) + LOG(U_1) + LOG(U_2) + LOG(U_3),$

  $\ldots$ etc. $\ldots$

  $LOG(B) + LOG(U_1) + LOG(U_2) + LOG(U_3) + \ldots + LOG(U_N) \}$

Clearly log series here is an additive random walk on the log-axis. It can be assumed that each new element gives rise to totally new mantissa value, and this in turns implies an overall uniform and flat 'density' on the mantissa space of (0, 1) albeit in a discrete manner. Hence the random growth series is Benford!

In contrast with the digital pitfalls of deterministic anomalous exponential growth series, random exponential growth series on the other hand are always nicely logarithmic. Here we never stumble upon the perils of LOG(F) whose multiples always add up exactly to an integral value on the log-axis in a consistent manner (i.e. anomalous series). Repeated additions of random $LOG(F_{RANDOM})$ values result in covering the entire (0, 1) mantissa space fairly and evenly whenever there are plenty of such accumulations.

As an example of a specific random model, log values march forward along the log-axis as in $LOG_{N+1} = LOG_N + Uniform(0, 1)$, as depicted in Figure 18. More generally, log values may be modeled as in: $LOG_{N+1} = LOG_N + [Almost Any Positive Random Variable]$.

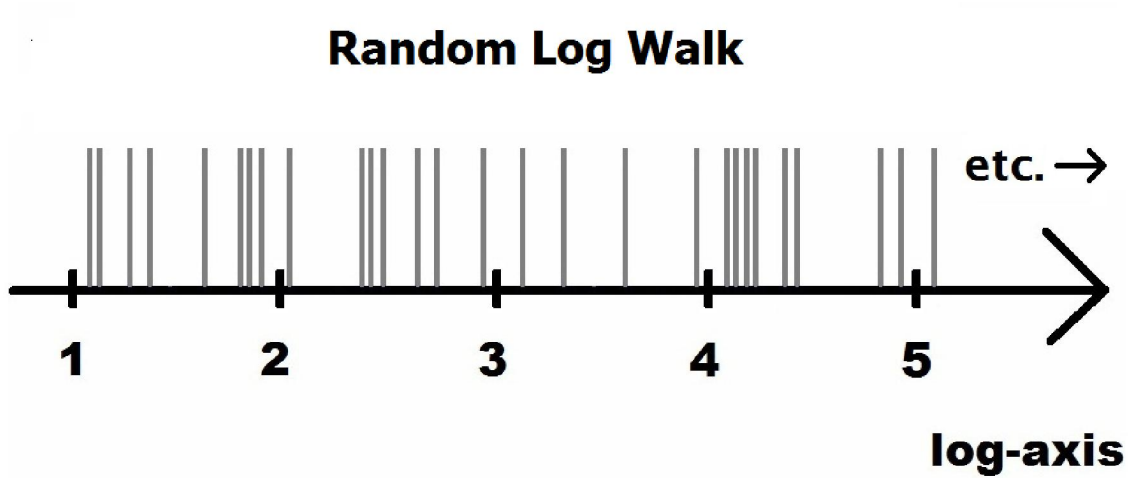

**Figure 18:** Uniformity of Mantissa via Random Additions of Uniform(0, 1)

# [11]  Super Exponential Growth Series

Classic exponential growth series such as $\{B, BF^1, BF^2, BF^3, BF^4, \dots, BF^N\}$ are characterized by having a constant F multiplicative factor, and where the set of factors in going from one element to the next element is the vector $\{F, F, F, F, F, \dots, F\}$ with $F > 1$. Super exponential growth series is a case where the factors themselves are growing exponentially, and where logarithmic behavior is found in spite of the rather odd nature of such series!

The set of the growing factors themselves is defined as $\{F^1, F^2, F^3, F^4, F^5, \dots, F^N\}$ where F stands for the initial factor as well as for the factor by which the factors themselves are growing. The super exponential series itself is then written in terms of its construction, term by term, as:

$\{B, \quad B(F^1), \quad BF^1(F^2), \quad BF^1F^2(F^3), \quad BF^1F^2F^3(F^4), \quad BF^1F^2F^3F^4(F^5), \quad BF^1F^2F^3F^4F^5(F^6), \quad \dots, \\ BF^1F^2F^3F^4F^5F^6\,F^7F^8\dots F^{N-2}\,F^{N-1}\,(F^N)\}$

which is simplified as the super series:

$\{B, BF^1, BF^3, BF^6, BF^{10}, BF^{15}, BF^{21}, \dots, BF^{(N*N + N)/2}\}$

This is essentially of the same format of the classic exponential growth series, but instead of having the exponents of the factors increasing sequentially as in 1, 2, 3, 4, 5, 6, etc., they are expanding more rapidly as in 1, 3, 6, 10, 15, 21, etc., namely as in $(N*N + N)/2$. In one particular computer calculation example, base B is chosen as 100, and the factor of the factors F is chosen as 1.0008. The 1st digits distribution of this super exponential growth series considering the first 1328 elements is {31.4, 16.8, 12.3, 9.9, 7.4, 5.6, 5.9, 5.6, 5.0}, and SSD is a rather low value of 4.3, signifying a great deal of closeness to the logarithmic.

Let us examine how the series advances along the log-axis:

{LOG(B),
 LOG(B) + LOG(F),
 LOG(B) + LOG(F) + 2*LOG(F),
 LOG(B) + LOG(F) + 2*LOG(F) + 3*LOG(F),
 LOG(B) + LOG(F) + 2*LOG(F) + 3*LOG(F) + 4*LOG(F),
 LOG(B) + LOG(F) + 2*LOG(F) + 3*LOG(F) + 4*LOG(F) + 5*LOG(F),
 LOG(B) + LOG(F) + 2*LOG(F) + 3*LOG(F) + 4*LOG(F) + 5*LOG(F) + 6*LOG(F),
… etc. … }

Clearly, **the super growth series speeds up along the log-axis**, namely that distances between consecutive log values of the super series are increasing, assuming the series represents exponential growth so that $F > 1$ and thus LOG(F) is positive, as opposed to exponential decay where $F < 1$ and where LOG(F) is negative.

It is very hard to imagine any scenario here where a particular value of F could lead to a situation where exactly N number of cycles (periods) always yields an integral value thus leading to an anomalous series where only very few selected mantissa values are being generated over and over again by the series. Hence it seems that any possible pitfall of LOG(F) = L/T rational number is totally irrelevant here to logarithmic behavior. Surely, the super exponential growth series does not march in constant steps of LOG(F), but rather in steps of [increasingly] multiple values of LOG(F).

The rationale for the logarithmic behavior of super exponential growth series is that it speeds up along the log-axis in a random and haphazard fashion as far as mantissa is concern. In other words, it is logarithmic because the series tramps upon the log-axis in a disorganized and uneven manner in relation to the log integers, resulting in the accumulation of all sorts of mantissa values. This in turns guarantees that mantissa is approximately uniform since the process gives no preference to any type or any sub-set of mantissa space, but rather keeps picking truly 'random' mantissa values as the series marches along.

Figure 19 depicts one such [imaginary] accelerated march of super exponential series along the log-axis, where log values are clearly speeding up and distances between consecutive log values are constantly increasing. See Kossovsky (2014) chapter 99 for more details and discussion about super exponential growth series.

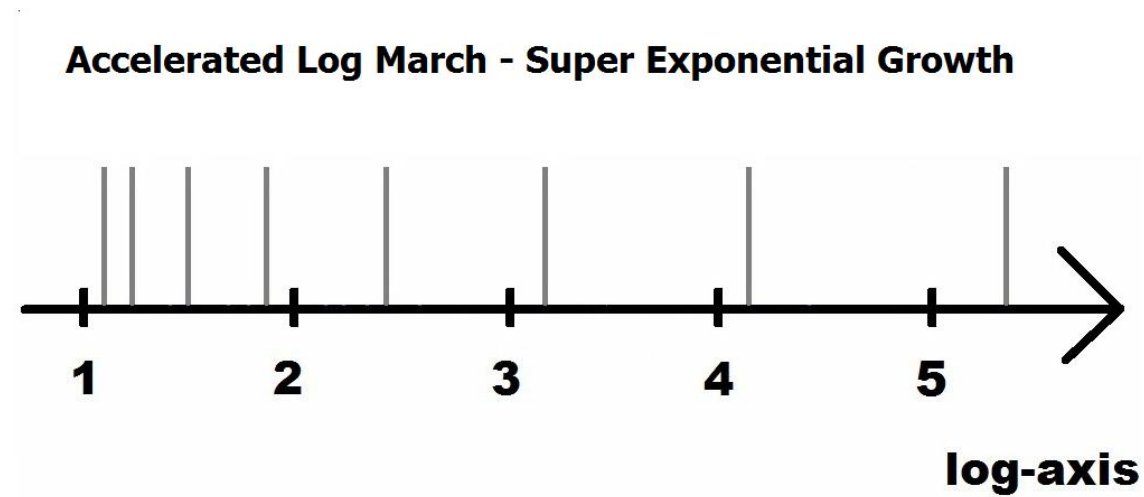

**Figure 19:** Uniformity of Mantissa via Accelerated Log-Axis March

## [12]  The Factorial Sequence as Super Exponential Series

The Factorial Sequence {N!} = {1!, 2!, 3!, 4!, …} is logarithmic in the limit as N goes to infinity, or in a practical sense and approximately so as N becomes quite large. The sequence is considered to be a deterministic multiplicative process; there is nothing random about it.

Examining log distances between consecutive elements of the Factorial Sequence we get:

$LOG(X_{N+1}) - LOG(X_N)$
$LOG((N + 1)!) - LOG(N!)$
$LOG(N!(N + 1)) - LOG(N!)$
$LOG(N!) + LOG(N + 1) - LOG(N!)$
**$LOG(N + 1)$**

Clearly, the Factorial Sequence marches along the log-axis in an accelerated way where distances between consecutive elements are constantly increasing, as in super exponential growth series. This nicely explains the digital logarithmic nature of the sequence! The sequence keeps minting totally new mantissa values in an unorganized way and without any coordination with respect to integral log marks. Figure 20 depicts the accelerated march along the log-axis for the first 19 elements of the Factorial Sequence.

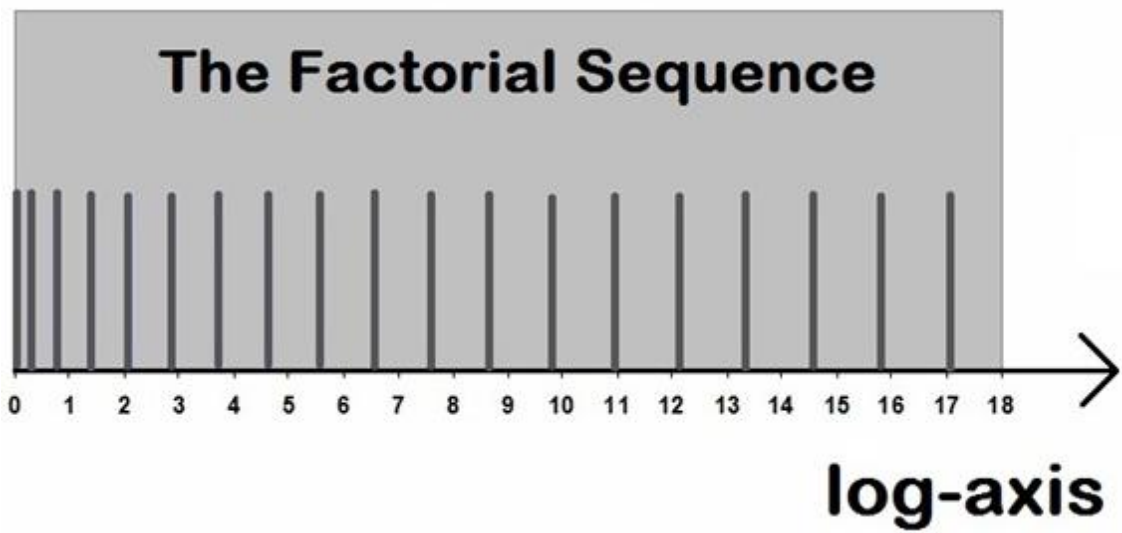

**Figure 20:** Accelerated March Along the Log-Axis – Factorial Sequence

An alternative vista of the Factorial Sequence is to view it simply as an exponential growth with increasing growth factors. This vista is surely justified when written as {N!} = {1!, 2!, 3!, 4!, …} = {1, (1)*2, (1*2)*3, (1*2*3)*4, …}. Hence the sequence can be written as in: $X_{N+1} = N*X_N$. In conclusion: The Factorial Sequence is also a super exponential growth series, and therefore its accelerated march along the log-axis is certainly expected, consistent, and in harmony with all that was discussed in the previous chapter.

Empirically testing the finite and short Factorial Sequence {1!, 2!, 3!,  …  , 168!, 169!, 170!} yields the first digits distribution {31.8, 17.1, 12.9, 7.1, 7.1, 5.9, 3.5, 8.2, 6.5}, with SSD = 30.1.

## [13]  The Self-Powered Sequence as Super Exponential Series

The Self-Powered Sequence $\{N^N\} = \{1^1, 2^2, 3^3, 4^4, \ldots\}$ is logarithmic in the limit as N goes to infinity, or approximately so as N becomes large. This sequence is a deterministic one, yet it differs profoundly from the classic exponential series as there exists no obvious connecting link between consecutive terms by the way of a simple multiplicative factor as in $X_{N+1} = Factor*X_N$.

Examining log distances between consecutive elements of the Self-Powered Sequence we get:

$LOG(X_{N+1}) - LOG(X_N)$
$LOG((N+1)^{(N+1)}) - LOG(N^N)$
$(N+1)*LOG(N+1) - N*LOG(N)$
$N*LOG(N+1) + LOG(N+1) - N*LOG(N)$
$N*LOG(N+1) - N*LOG(N) + LOG(N+1)$
$N*LOG((N+1)/N) + LOG(N+1)$
$N*LOG(1+1/N) + LOG(N+1)$
$LOG((1+1/N)^N) + LOG(N+1)$

$\lim_{N \to \infty} \left(1 + \frac{1}{N}\right)^N = e$, hence in the limit as N goes to infinity, log distance in between is:

### $LOG(e) + LOG(N+1)$

Clearly, the Self-Powered Sequence marches along the log-axis in an accelerated way, and where distances between consecutive elements are constantly increasing, as in super exponential growth series. This nicely explains the digital logarithmic nature of the sequence! The sequence keeps minting totally new mantissa values in an unorganized way and without any coordination with respect to integral log marks. Figure 21 depicts the accelerated march along the log-axis for the first 14 elements of the Self-Powered Sequence.

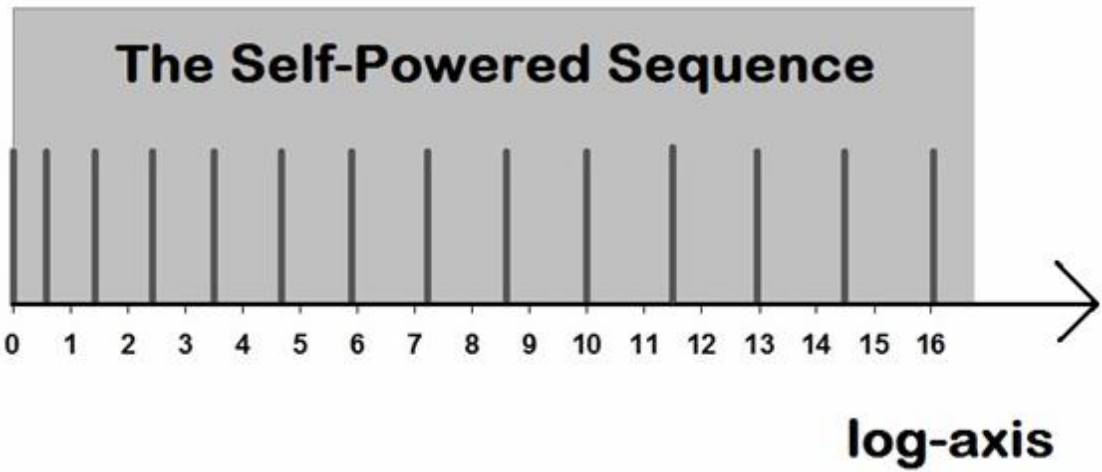

**Figure 21:** Accelerated March Along the Log-Axis – Self-Powered Sequence

Empirically testing the finite and short Self-Powered Seq. $\{1^1, 2^2, 3^3, \ldots, 141^{141}, 142^{142}, 143^{143}\}$ yields the first digits distribution $\{36.4, 14.7, 16.8, 9.1, 6.3, 2.8, 4.2, 7.0, 2.8\}$, with SSD = 93.6.

## [14] Four Distinct Styles of Log-Axis Walk Leading to Benford

Figure 22 summarizes the four distinct ways a sequence could march along the log-axis yielding uniformity of mantissa; (1) Normal exponential growth series with LOG(F) being an irrational number; (2) Anomalous exponential growth series with a rational LOG(F) of a very small value; (3) Super exponential growth series, the Factorial Series, the Self-Powered Series, all characterized by accelerated march along the log-axis where distances between consecutive log values are constantly increasing; (4) Random walk along the log-axis.

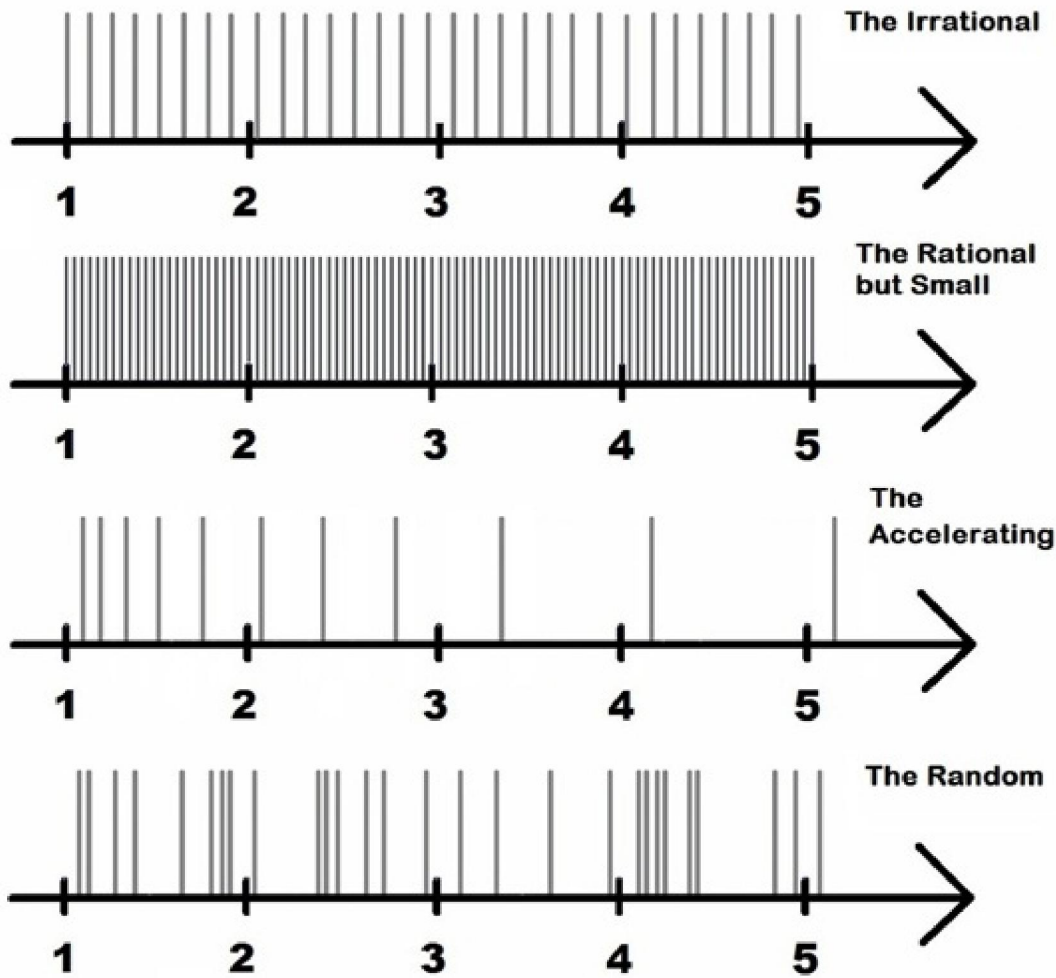

**Figure 22:** Four Distinct Logarithmic Ways for Sequences to March along the Log-Axis

What should we expect a priori - in extreme generality - from the vast majority of mathematical sequences with regards to Benford behavior, without performing the tedious and endless task of examining each and every one of them? The obvious answer is that the vast majority of them are surely Benford, because it's a bit rare for any particular sequence to march along the log-axis in a coordinated and organized way with respect to the integer marks of the log-axis. Such a deliberate march requires a great deal of effort and close attention on the part of the sequence, and very few sequences are willing to work so hard so as to obtain the abstract, intangible, and vague goal of uniformity of mantissa. Most sequences prefer to march merrily, freely, and relaxingly along the log-axis as they see fit, guaranteeing Benford behavior for the corresponding values they generate on the x-axis.

# [15]  Continuous Exponential Growth Series

Thus far we have dealt only with discrete exponential growth series where quantity count takes place at the end of a particular length of time, say at the end of each month, or at the end of each year. Perhaps the quantity actually jumps at an instant over and over again at the end of each period of time, such as a bank account gaining interest at the last minute of the last day of each month, or earning interest at an instant just before midnight on each December 31 day.

For a continuous exponential growth series, the growing variable X is a function of time, and it is expressed as $X(T) = BF^T$ where B is the initial base value at time 0, F is the fixed multiplicative factor of the growth rate per one integral unit time period, and T is the continuous real time variable, comprising not only of integral values such as 1 year, 2 years, 3 years, and so on, but also of the unaccountably infinite values of all the real points in between the integers. Assuming that the time unit in the definition of factor F is measured in years, then quantity X is being continuously recorded say every minute, or even every second if humanly possible, and not only once a year. Hence the data set in focus here is not merely the small collection of say December 31 readings for several years or few decades, but rather the extremely large data set comprising the detailed evolution of quantity X examined each minute, or even each second if possible, and so on, recording X extremely frequently each tiny interval of time. Of interest is the digital configuration of such large data set, and the quest to demonstrate that the data set - having proper span by standing exactly between integral powers of ten or having an integral exponent difference - is perfectly logarithmic in the limit as time measurement shrinks from hours, to minutes, to seconds, and then to infinitesimal small time intervals, approaching zero from above in the limit. In reality for all practical purposes, recording quantity X each minute is considered to be sufficiently small as a time frame for Benford behavior when the growth period in the definition of F is measured in years.

As an example, let us consider 5% exponential growth from X = 1 to X = 10. Here variable X as a function of time is expressed as $X(T) = 1*(1.05^T)$. The motivation for the start at 1 emanates from the fact that at 1 we start anew the natural cycle of the 1st digits. The motivation for the termination at 10 emanates from the fact that at 10 we end the first complete cycle of the 1st digits. Such a deliberate span from 1 to 10 ensures that all nine digits would cycle exactly once and fully so, so that each digit would be allowed to express its full potential; and that all the digits enjoy equal opportunity. Figure 23 depicts this particular continuous exponential growth series from 1 to 10 with F = 1.05, namely 5% growth per period. Clearly, for any X value within the interval [1, 2) the first digit is always 1, such as for X values 1.345, 1.876503, 1.5, 1.0, 1.999, 1.06214, and so on. In general here, all values falling on the sub-interval [d, d + 1) are such that the first leading digit is d (assuming d is an integer from 1 to 9).

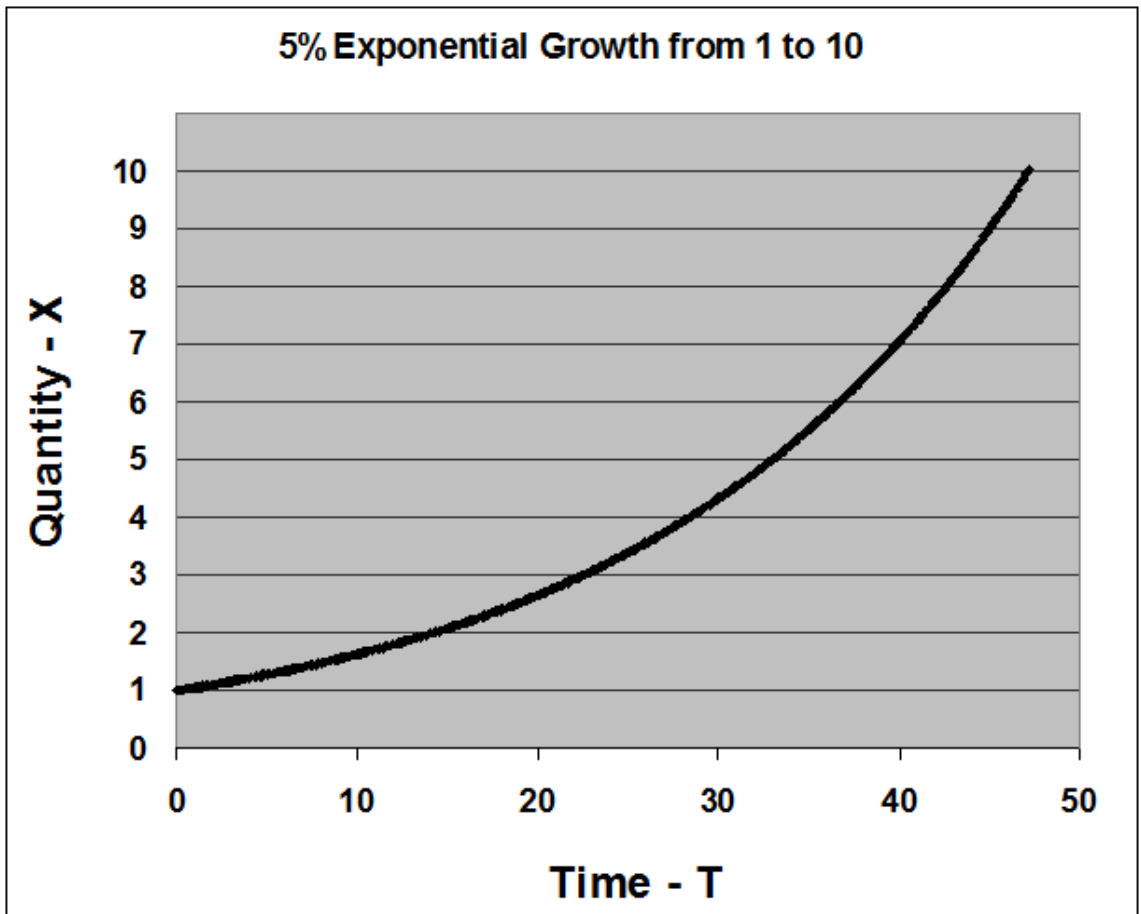

**Figure 23:** Continuous 5% Exponential Growth from 1 to 10

Hence, the amount of time it takes the series to grow from 1 to 2 is equivalent to the amount of time the series spends while digit 1 is leading; and the amount of time it takes the series to grow from 2 to 3 is equivalent to the amount of time the series spends while digit 2 is leading, and so on. What is needed to be calculated here is simply the relative proportions of the time spent on all the nine first digits.

Figure 24 depicts the critical moments when the series encounters those crucial integral values of 1, 2, 3, 4, 5, 6, 7, 8, 9, and 10, when the first digit changes.

Figure 25 depicts the time intervals the series spends on each digit, as well as the proportion of time each digit earns, and which 'happened' to be exactly Benford!

| t - time | X - Quantity |
|----------|--------------|
| 0.0 | 1.0 |
| 14.2 | 2.0 |
| 22.5 | 3.0 |
| 28.4 | 4.0 |
| 33.0 | 5.0 |
| 36.7 | 6.0 |
| 39.9 | 7.0 |
| 42.6 | 8.0 |
| 45.0 | 9.0 |
| 47.2 | 10.0 |

**Figure 24:** Integral Quantity Points & Associated Time

| Change in X | Time Interval | % Time Interval |
|-------------|---------------|-----------------|
| 1 to 2 | 14.2 | 30.1% |
| 2 to 3 | 8.3 | 17.6% |
| 3 to 4 | 5.9 | 12.5% |
| 4 to 5 | 4.6 | 9.7% |
| 5 to 6 | 3.7 | 7.8% |
| 6 to 7 | 3.2 | 6.8% |
| 7 to 8 | 2.7 | 5.7% |
| 8 to 9 | 2.4 | 5.1% |
| 9 to 10 | 2.2 | 4.7% |

**Figure 25:** Relative Time Spent on the Various Nine Digits is Benford

In practical terms and as far as actual data sets are concerned in the context of Benford's Law, what is meant abstractly by the phrase "the relative proportion of the time spent on digit d" is simply equivalent to the proportion of numbers the data analyst has actually written down recording the growth with digit d leading, as compared with all the written numbers recorded while quantity develops from 1 all the way to 10.

For example, let us assume that the period is defined as one year, and that the data analyst measures the growing quantity every minute as it develops from 1 all the way to 10. Assuming a year with 365 days, 24 hours a day, and 60 minutes an hour, then this yields 365*24*60 or 525,600 total minutes per year. The time it takes to reach quantity 10 is 47.2 years, hence the data analyst has recorded (47.2)*(525,600) = 24,808,320 data points (for each minute) in total. The data analyst has observed that from year 33.0 to year 36.7 digit 5 was consistently leading the numbers during these 3.7 years, and which consists of (3.7)*(525,600) = 1,944,720 data points. Hence according to these observations, the proportion of digit 5 leading here is calculated as (1,944,720)/(24,808,320) = 0.0784, and which is in accordance with Benford's Law.

Figures 26 and 27 depict clearly the time intervals on the T-axis corresponding to first digit occurrences on the X-axis.

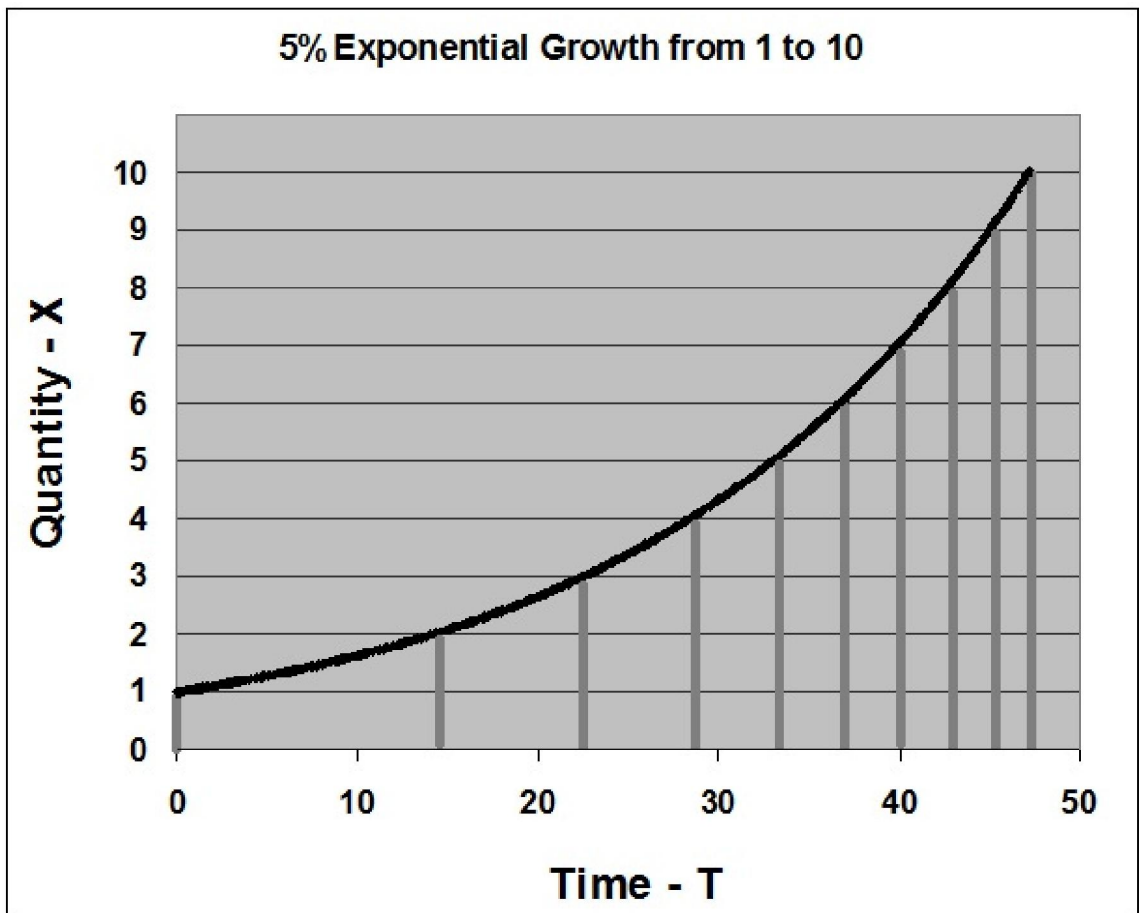

**Figure 26:** Specific Time Intervals Spent on the Nine Digits

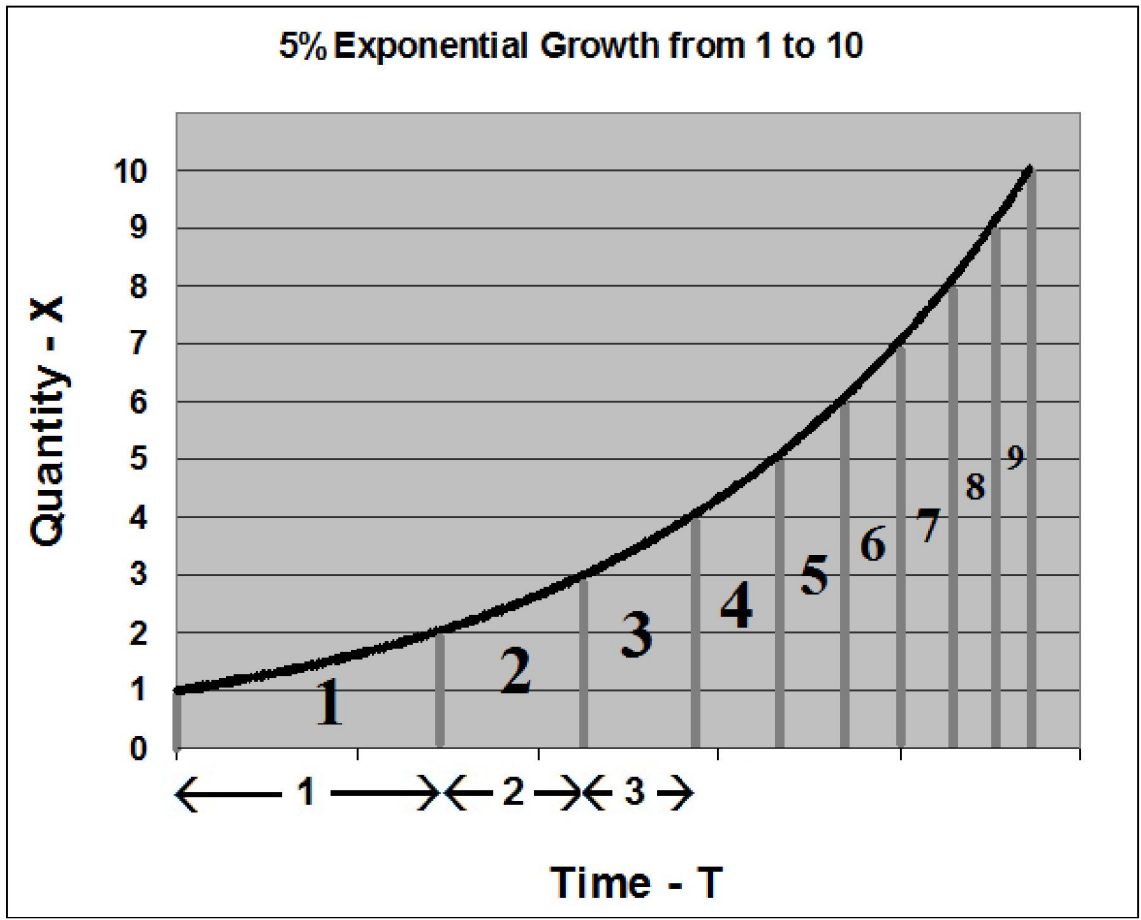

**Figure 27:** Specific Time Intervals Spent on the Nine Digits – Detailed

Let us consider **the general continuous case** of the exponential growth series $X(T) = 1*F^T$ having **any growth rate with factor F**, with an initial quantity 1 at time 0, and limited to the progression from 1 to 10. The time $T_2$ spent to get to quantity 2 is obtained from the expression itself as follows:

$2 = F^{T_2}$

$LOG_F(2) = LOG_F(F^{T_2})$

$LOG_F(2) = T_2 \cdot LOG_F(F)$

$LOG_F(2) = T_2 \cdot 1$

$LOG_F(2) = T_2$

The time $T_3$ spent to get to quantity 3 is obtained from the expression itself as follows:

$3 = F^{T_3}$
$LOG_F(3) = LOG_F(F^{T_3})$
$LOG_F(3) = T_3 \cdot LOG_F(F)$
$LOG_F(3) = T_3 \cdot 1$
$LOG_F(3) = T_3$

The total time $T_{10}$ spent to get to quantity 10 is obtained from the expression itself as follows:

$10 = F^{T_{10}}$
$LOG_F(10) = LOG_F(F^{T_{10}})$
$LOG_F(10) = T_{10} \cdot LOG_F(F)$
$LOG_F(10) = T_{10} \cdot 1$
$LOG_F(10) = T_{10}$

The time $T_{10}$ is the time of one full cycle for all digits 1, 2, 3, 4, 5, 6, 7, 8, and 9.

Hence, time spent between 2 and 3 where 1st digit is always 2 =
Time to get to 3 – Time to get to 2 =
$LOG_F(3) - LOG_F(2) =$
$LOG_F(3/2)$

And as a proportion in comparison to the total time needed for one full cycle for all the digits in the progression from 1 to 10:

Overall time proportion for digit 2 = $LOG_F(3/2) / LOG_F(10)$

Utilizing the logarithmic identity
$LOG_A(X) = LOG_B(X) / LOG_B(A)$
for the numerator and for the denominator we get:

$LOG_F(3/2) / LOG_F(10)$
$[LOG_{10}(3/2) / LOG_{10}(F)] / [LOG_{10}(10) / LOG_{10}(F)]$
$[LOG_{10}(3/2)] / [LOG_{10}(10)]$
$[LOG_{10}(3/2)] / [1]$
$LOG_{10}(3/2)$
$LOG_{10}((2 + 1)/2)$
$LOG_{10}(1 + 1/2)$

Which is the proportion for digit 2 according Benford's Law!

Let us prove this in general for the time spent between (d) and (d+1), where 1st digit is always d:

The time $T_d$ spent to get to quantity d is obtained from the expression itself as follows:

$d = F^{T_d}$

$LOG_F(d) = LOG_F(F^{T_d})$

$LOG_F(d) = T_d \cdot LOG_F(F)$

$LOG_F(d) = T_d \cdot 1$

$LOG_F(d) = T_d$

Similarly the time spent to get to quantity (d + 1) is $LOG_F(d + 1) = T_{d+1}$.

Hence, time spent between (d) and (d + 1) where 1st digit is always d =
Time to get to (d + 1) − Time to get to (d) =
$LOG_F(d + 1) − LOG_F(d) =$
$LOG_F((d + 1)/d)$

And as a proportion in comparison to the total time needed for one full cycle for all the digits in the progression from 1 to 10:

Overall time proportion of digit d = $LOG_F((d + 1)/d) / LOG_F(10)$

Utilizing the logarithmic identity
$LOG_A(X) = LOG_B(X) / LOG_B(A)$
for the numerator and for the denominator we get:

$LOG_F((d + 1)/d) / LOG_F(10)$
$[LOG_{10}((d + 1)/d) / LOG_{10}(F)] / [LOG_{10}(10) / LOG_{10}(F)]$
$[LOG_{10}((d + 1)/d)] / [LOG_{10}(10)]$
$[LOG_{10}((d + 1)/d)] / [1]$
$LOG_{10}((d + 1)/d)$
$LOG_{10}(d/d + 1/d)$
$LOG_{10}(1 + 1/d)$

Which is exactly Benford's Law!

Little reflection is needed to realize that the above proof can be extended to include exponential growth spanning any adjacent integral powers of ten points such as growth from 10 to 100, from 100 to 1000, or from 1000 to 10000, and so on. The transformation of (1, 10) growth into (10, 100) growth is leading-digits-neutral, since digit 3 leads the number 3.987055 say, just as it does lead the number 39.87055. Moreover, even for growth spanning non-adjacent multiple powers of ten points, such as say from 1 to 10000, the same argument as in the above proof holds separately for each sub-interval spanning adjacent integral powers of ten, so that the entire span of (1, 10000) is Benford since it is so on each of its constituent sub-intervals (1, 10), (10, 100), (100, 1000), and (1000, 10000). Furthermore, any continuous exponential growth series spanning integral exponent difference such as those growing say from 3 to 30, from 7.8233 to 78.233, or from 43.129 to 43129, and in general growing from R to R*$10^{\text{INTEGER}}$, are all perfectly Benford as well. In a nutshell, for the example of 3 to 30 growth, viewed in terms of the generic 1 to 10 growth, the 'missing' sub-interval of (1, 3) on the left, is being exactly and perfectly substituted by the 'extra' sub-interval of (10, 30) on the right.

The above proof has been constructed in extreme generality, for any continuous growth rate of factor F, and without requiring LOG(F) to be irrational (except that the series is required to start at 1 and to terminate at 10). Hence, on the face of it, this last result seems to contradict the constraint which excludes anomalous exponential growth series from logarithmic behavior. As an example, for the continuous **5.92537251773%** yearly growth, the yearly steps taken along the log-axis are in units of LOG(1 + 5.92537251773/100) =  LOG(1.0592537251773) =  0.025 = 1/40 = rational number, yet the above proof guarantees logarithmic behavior here!? Nonetheless, harmony and consistency prevail, and our delicate and complex digital edifice does not fall into ruin and contradiction in the least, since what is actually recorded here is not how quantity grows per year, but rather how quantity grows say per second, and that rate is the extremely low **0.0000001825%** growth, calculated as 100*[(60*60*24*365)th root of 1.0592537251773 - 1]. LOG(F) for this very low growth rate per second is calculated as the extremely tiny value of LOG(1 + 0.0000001825/100) = 0.0000000007927447761.

Admittedly, this extremely low rate of 0.0000001825% is also anomalous since its LOG(F) is also rational, yet the exact logarithmic behavior of continuous exponential growth series with rational LOG(F) comes under the protective umbrella of the exception given for anomalous rates whenever LOG(F) is sufficiently small compared with unit length on the log-axis (i.e. growth rate is very low – relating in general to high T values). The fact that continuous exponential growth series is considered here as opposed to the discrete integral case implies that we take the population pulse each very tiny interval of time, say each second, and therefore the growth rate per each such infinitesimal small time interval is extremely small and approaches 0% growth rate from above in the limit. Not only does the creature lack long legs here for skipping undesirable mantissa sections, but it has such extremely short ones that it is literally forced by default to walk all over the (0, 1) mantissa interval, not omitting any spots or corners at all, hence the series is perfectly and exactly logarithmic.

## [16]  The Discrete Embedded within the Continuous

Turning the **continuous** 5% exponential growth from 1 to 10 into a **discrete** series by measuring quantity $X(T) = 1*(1.05^T)$ once each period (say once a year), namely utilizing only the integral values of T = {0, 1, 2, 3, … , 47}, yields 48 elements (including the initial base value considered as the 1st element):

| | | | | | | | | | | | |
|---|---|---|---|---|---|---|---|---|---|---|---|
| 1.00 | 1.05 | 1.10 | 1.16 | 1.22 | 1.28 | 1.34 | 1.41 | 1.48 | 1.55 | 1.63 | 1.71 |
| 1.80 | 1.89 | 1.98 | 2.08 | 2.18 | 2.29 | 2.41 | 2.53 | 2.65 | 2.79 | 2.93 | 3.07 |
| 3.23 | 3.39 | 3.56 | 3.73 | 3.92 | 4.12 | 4.32 | 4.54 | 4.76 | 5.00 | 5.25 | 5.52 |
| 5.79 | 6.08 | 6.39 | 6.70 | 7.04 | 7.39 | 7.76 | 8.15 | 8.56 | 8.99 | 9.43 | 9.91 |

In Figure 28 the black dots represent the quantities for each integral time value. It should be noted that these 48 points are perfectly embedded within the continuous curve of 5% growth of the previous charts as seen in Figures 23 to 27.

Here 1st digits distribution is {31.3, 16.7, 12.5, 8.3, 8.3, 6.3, 6.3, 6.3, 4.2}, and SSD = 6.3, having only slight deviation from the perfect Benford behavior of the continuous growth.

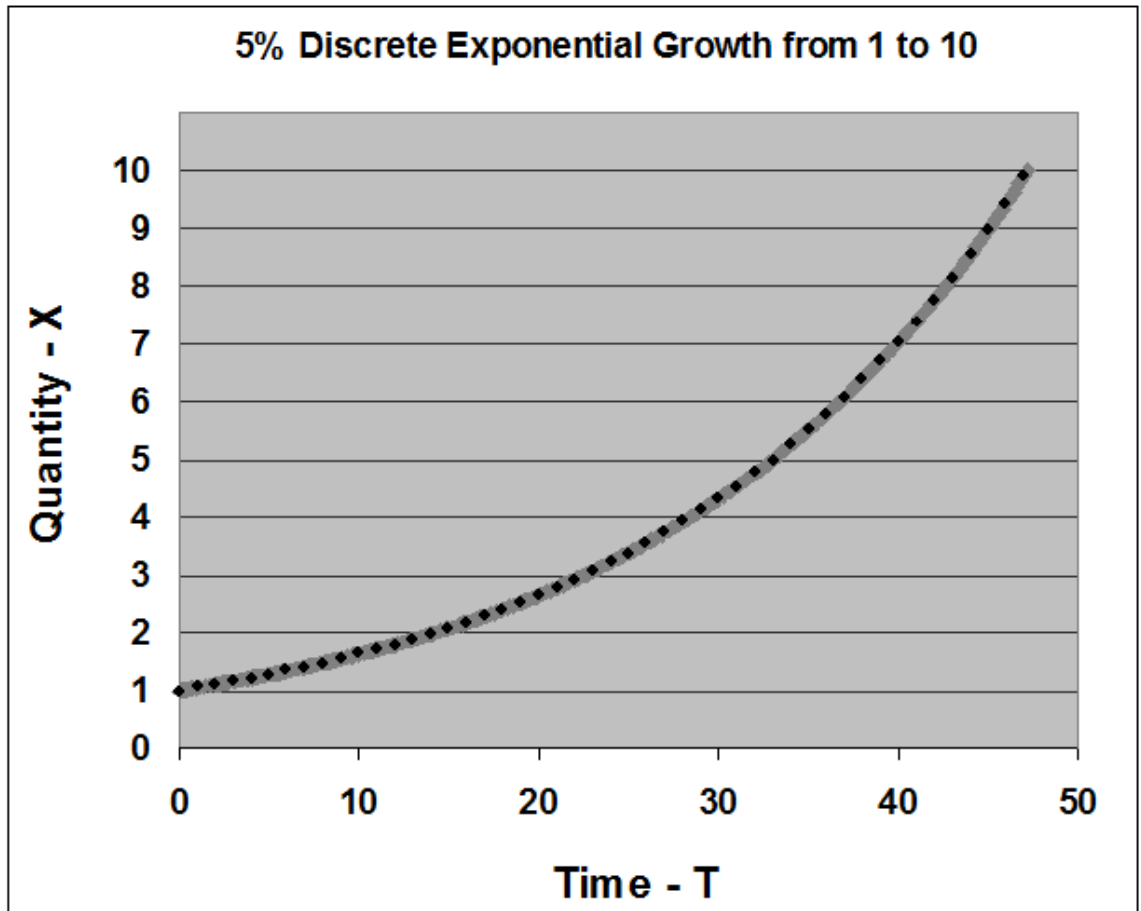

**Figure 28:** The Discrete 5% Growth Embedded within the Continuous 5% Growth

Since this 5% exponential growth series is discrete; and since 5% is not as small a rate as to be considered immune to rebellious behavior, there exists the possibility of it being an anomalous series. As it happened, LOG(1.05) = 0.0211893, and this is not a rational number. Moreover, the nearest rational growth series are the mildly rebellious series of 4.91397291% belonging to the rationality $L/T = 1/48 = 0.0208333$, and 5.02110796% belonging to the rationality $L/T = 1/47 = 0.0212766$. Both of these potentially disruptive series are quiet mild in nature, and they are not as near our 5% normal series, all of which guarantees the approximate logarithmic behavior of the 5% discrete series.

Turning the **continuous** 5% exponential growth from 1 to 10 into a **discrete** series by measuring quantity $X(T) = 1*(1.05^T)$ once every 3 periods (say once each 3-year interval), namely utilizing the particular set of values of T = {0, 3, 6, 9, 12, 15, 18, 21, 24, 27, 30, 33, 36, 39, 42, 45}, yields 16 elements:

| 1.00 | 1.16 | 1.34 | 1.55 | 1.80 | 2.08 | 2.41 | 2.79 | 3.23 | 3.73 | 4.32 | 5.00 |
|------|------|------|------|------|------|------|------|------|------|------|------|
| 5.79 | 6.70 | 7.76 | 8.99 | | | | | | | | |

In Figure 29 the black dots represent the quantities for each 3-year time interval. These 16 points are perfectly embedded within the continuous curve of 5% growth, although the effective (cumulative) growth factor from one point to the next, namely the growth factor for each 3-year period, is $(1.05)^3 = 1.157625$, and thus the associated growth rate here is actually 15.8%.

Here 1st digits distribution is {31.3, 18.8, 12.5, 6.3, 12.5, 6.3, 6.3, 6.3, 0.0}; having some noticeable deviation from the perfect Benford behavior of the continuous; especially having the unusual 0% proportion for digit 9, and with SSD = 58.1. The reason for the deterioration in digital fit to Benford for the higher 15.8% growth of the 3-year measurements program is due to lack of sufficient number of data points, being too-discrete-like. The yearly measurements program of discrete 5% growth enjoys the advantage of having more points, and thus being more continuous-like. Yet it should be emphasized that both of these discrete series represent (approximately) the same overall growth from 1 to about 10 – and both are embedded within the same continuous 5% growth curve. The 3-year measurements program grows from 1 to 8.99, while the yearly measurements program grows from 1 to 9.91, and such similarity in overall growth should be acknowledged and considered.

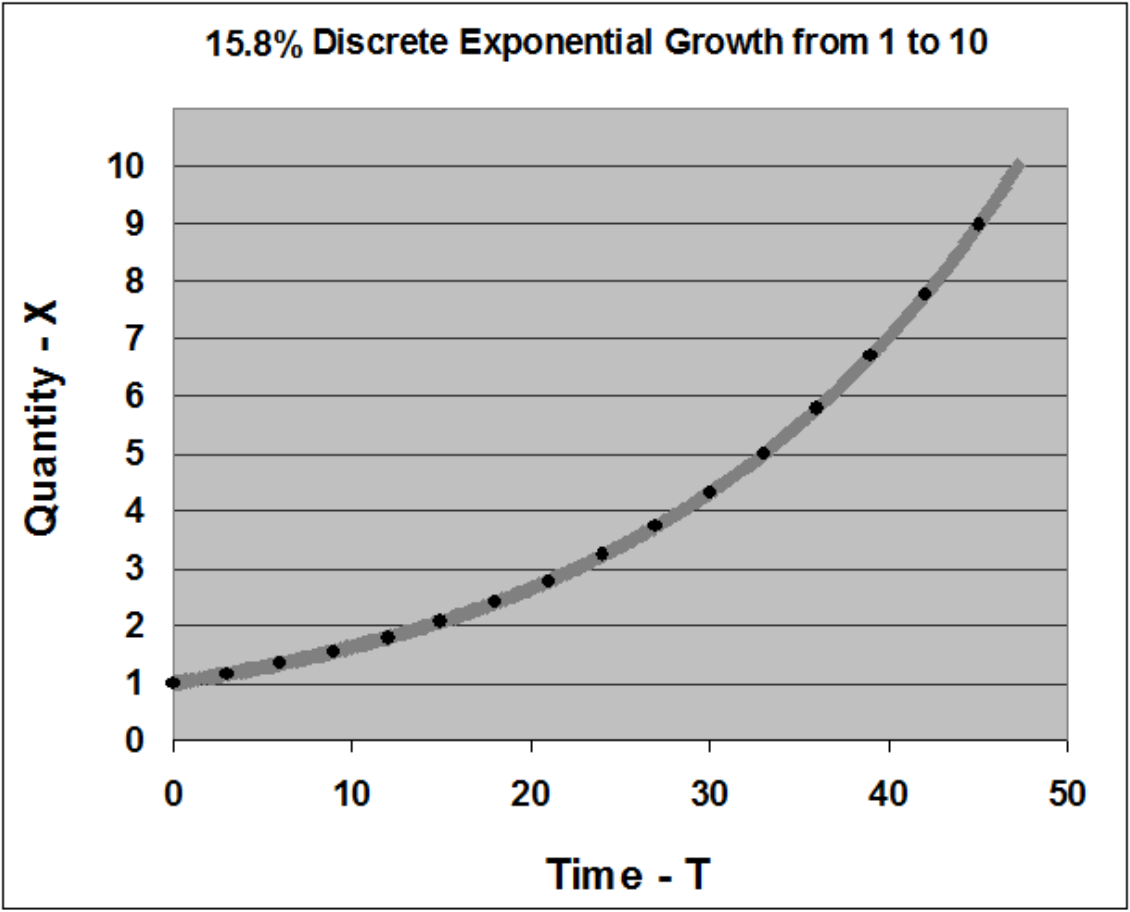

**Figure 29:** The Discrete 15.8% Growth Embedded within the Continuous 5% Growth

Turning the **continuous** 5% exponential growth from 1 to 10 into a **discrete** series by measuring quantity $X(T) = 1*(1.05^T)$ once every 5 periods (once each 5-year interval), namely utilizing the particular set of values of T = {0, 5, 10, 15, 20, 25, 30, 35, 40, 45}, yields the very small set of 10 elements:

1.00   1.28   1.63   2.08   2.65   3.39   4.32   5.52   7.04   8.99

In Figure 30 the black dots represent the quantities for each 5-year time interval. These 10 points are perfectly embedded within the continuous curve of 5% growth, although the effective (cumulative) growth factor from one point to the next, namely the growth factor for each 5-year period, is $(1.05)^5 = 1.276281$, and thus the associated growth rate here is actually 27.6%.

Here 1st digits distribution is {30, 20, 10, 10, 10, 0, 10, 10, 0}, with SSD = 123.6. The markedly larger deviation from the Benford digital configuration of the continuous case can be squarely blamed on the severe lack of sufficient number of data points here. The series is being severely punished for its discreteness. Yet, potentially, this 27.6% series can still remedy itself by simply continuing to march forward along the 5% continuous growth curve, well beyond 8.99 and even well beyond 100, accumulating many more points and passing numerous integral power of ten points along the way, and then successfully becoming nearly logarithmic.

Since $LOG((1.05)^5)$ is not a rational number; the series is then not rebellious, and its re-entrance points into each new integral power of ten interval is random-like and different, generating a variety of mantissa values. Since the 5% territory upon which the 27.6% series stamps upon is itself perfectly Benford, then randomly 'picking' points from it would eventually lead to a set of values which are also Benford.

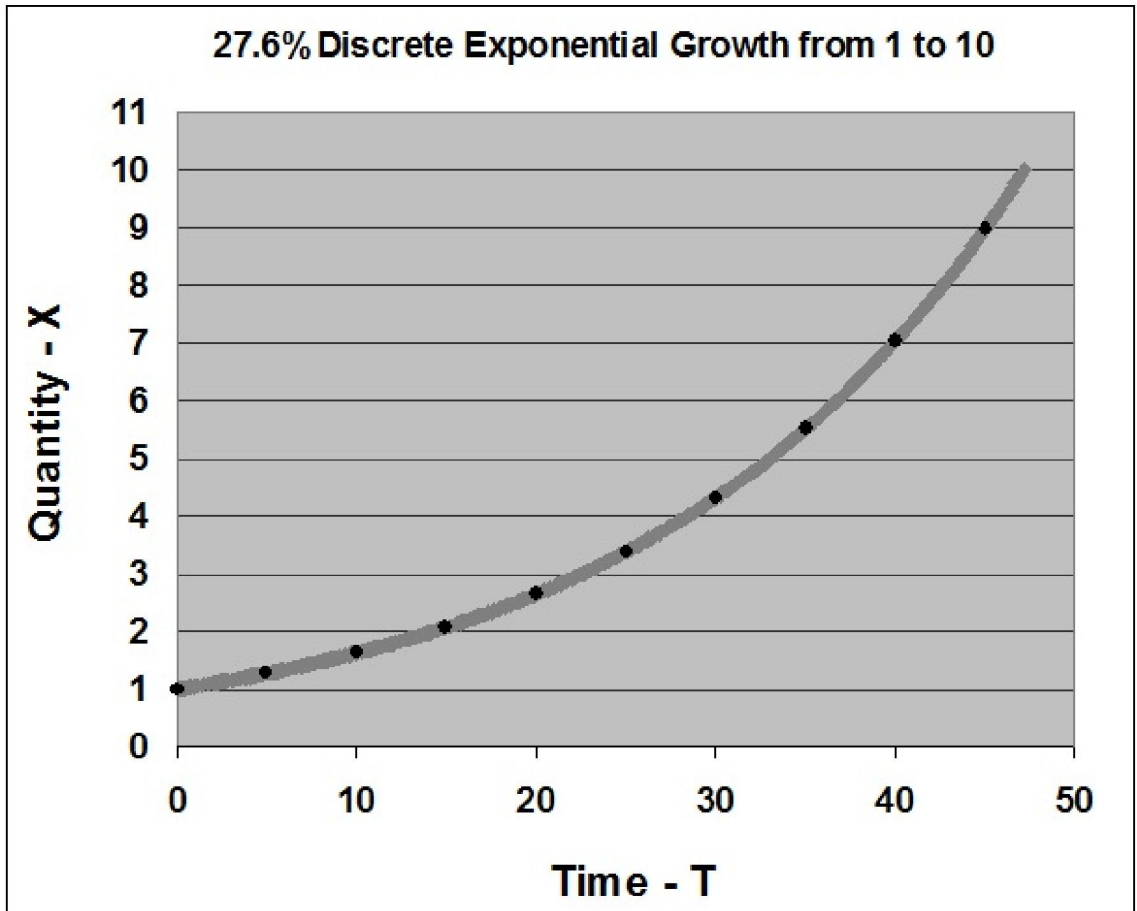

**Figure 30:** The Discrete 27.6% Growth Embedded within the Continuous 5% Growth

It must be emphasized once again that the same physical <u>continuous</u> growth (say a particular bacterial population growing in a lab) yields distinct digital configurations depending on how often <u>discrete</u> readings of population count takes place. As a much more illustrative and dramatic example than the 1-year versus 3-year or 5-year data recording schemes above, let us consider **5% per minute** exponential growth of the laboratory culture bacteria Salmonella Choleraesuis from an approximate initial level of 100,000 as in $\mathbf{X(T) = 100000*(1.05^T)}$, and two data analysts named George and Frank having two distinct measuring styles. George is very fastidious and hard-working, thus he constantly takes measurements of the colony every minute. Frank on the other hand is quite laid-back and relaxed, his mind is focused more on personal matters than work, and he makes the effort to check the colony only once every hour. The general manager at the laboratory supervising George and Frank obliges them to come up with 170 carefully measured data points. This would require hard-working George only little less than 3 hours, while for Frank the task would require about 7 days. For George who stops by the laboratory every <u>minute</u>, the growth rate is 5%, while for Frank who stops by every hour, the growth factor F is $(1.05)^{60}$ = 18.67918, namely a whopping 1767.9% <u>hourly</u> growth rate.

**Frank** fast hourly growth series has the huge 214.86 order of magnitude. The first 5 elements and last 3 elements are shown below:

100000        1867919      34891199     651739184    12173957374 … etc. …
$2.074*10^{217}$    $3.874*10^{218}$    $7.236*10^{219}$

**George** slow minute growth series has the small 3.58 order of magnitude. The first 5 elements and last 3 elements are shown below:

100000        105000        110250        115762       121550 … etc. …
345631204    362912765    381058403

Digital results for George and Frank are as follow:

**George** Series - Low 5% Growth:        {33.5, 19.4, 13.5,  8.2, 6.5, 5.9, 4.7, 4.7, 3.5}
**Frank**  Series  - High 1767.9% Growth:    {30.0, 17.1, 12.9, 10.0, 7.1, 7.1, 6.5, 4.1, 5.3}
**Benford's Law** 1st digits Distribution:    **{30.1, 17.6, 12.5,  9.7, 7.9, 6.7, 5.8, 5.1, 4.6}**

Hard-working George is quite disappointed. In spite of his dedication to the job and his best effort, first digits for him markedly deviate from Benford with SSD = 23.2, while easy-going Frank unfairly earns a much nicer fit to Benford with SSD = 3.4. The reason for the deterioration in digital fit for George's low growth is due to lack of better calibration of the ratio of last to first elements, not being close to an integral power of ten, being forced to abruptly stop immediately after 170 points as ordered by his boss, and such hasty termination of the series determines the value of the last element - disregarding what happens to exponent difference. Hence for George, the first element is 100000 and the last element is 381058403 with ratio of 381058403/100000, or 3810.58, which is unfortunately not close enough to any integral power of ten, such as 1000 from below or 10000 from above. In other words, that the exponent difference of 3.58 between the first and last elements is not close enough to an integer. For Frank, the first element is 100000 and the last element is $7.236*10^{219}$ with ratio of $7.236*10^{219}$/100000, or $7.236*10^{214}$, which admittedly is also not close to any integral power of ten, such as $10^{214}$ from below or $10^{215}$ from above; yet the good fortune of Frank emanates from the extremely high growth rate and the resultant very high order of magnitude of his series which excuses him from any kind of calibration attempt of the first and last elements whatsoever - as in all high growth series.

It should be noted that even though both studies done by George and by Frank trace discrete points on the same continuous 5% curve, and even though both studies are with the same number of data points, namely 170, yet they are indeed of different overall bacterial growth; one spans a shorter time interval and the other spans a much longer time interval; that they are of distinct orders of magnitude from the first to the last elements; and that bacteria grew for Frank much more than it grew for George.

# [17]  The Constraint of an Integral Exponent Difference for Low Growth

Let us vividly demonstrate with actual numerical examples the importance of the constraint regarding an integral exponent difference for exponential growth series; a constraint derived from the connection of exponential growth series to the k/x distribution as shall be discussed in the next chapter.

The phrase 'exponent difference' refers to the difference between the log values of the minimum (first element) and the maximum (last element). For the exponential series {2, 4, 8, 16, 32, 64, 128, 256, 512, 1024}, exponent difference is LOG(1024) – LOG(2) = 3.0103 – 0.3010 = 2.7093.

A typical real-life manifestation of exponential growth series is a single city growing at a constant rate, where the long stream of official population records by its census office measured at the end of each integral period of time (say each year on December 31) converges to Benford in the limit as the number of periods (years) goes to infinity. In practical terms, for those impatient to wait eternally for an infinite number of years until perfect Benford behavior appears, the rule of thumb for an approximate logarithmic behavior contains two distinct requirements: (**A**) Reasonably large (albeit finite) number of periods which could also be just say 100 or perhaps 500, depending on the desired level of accuracy in terms of obtaining a good fit to Benford. (**B**) The second requirement besides length considerations relates to a very particular relationship between the values of the first element (minimum) and the last element (maximum) of the series, namely that exponent difference LOG(Last) – LOG(First) should be as close as possible to an integer, or equivalently stated, that the ratio Last/First should be as close as possible to an integral power of ten such as 10, 100, 1000, and so forth.

In general, low growth series that adheres closely enough to the above two requirements are close enough to the logarithmic. For very high growth series, where the digital cycle is very short (i.e. passing integral powers of ten points quite frequently), having nearly an integral exponent difference is not really necessary; rather the only necessary condition is to have sufficient number of periods, and for such long series the logarithmic is approximately observed.

The combination of discreteness and finiteness of real-life exponential growth series inevitably ruins perfect logarithmic behavior, and this necessitates the consideration of large enough number of periods in order to get a reasonable fit to Benford in the approximate. As mentioned in Kossovsky (2014) in chapter 62 regarding Corollary I, the integral restriction on exponent difference for the k/x distribution becomes somewhat redundant in practical sense and can be ignored if exponent difference is itself quite large. How large? Well, that depends on the desired level of precision. Certainly, fairly large (non-integral) values of exponent difference over 25 or over 35 say, yield digit distribution that is extremely close to the logarithmic for all practical matters. By extension, such a waiver applies to exponential growth series as well, hence growth series with non-integral but fairly large exponent difference are almost Benford.

A brief explanation of the dispensation or exemption to the usual rule requiring an integral exponent difference for a growth series having say the very large 25.7 exponent difference, is that the vast majority of elements in the series are from the minimum (first) element until that element towards the end completing nearly exactly 25 exponent difference, while only a tiny minority beyond that element until the maximum (last) element are with the very odd 0.7 exponent difference. Since the majority of the elements with that nearly 25 exponent difference constitutes already by themselves a nearly logarithmic series, then the inclusion of a tiny extra minority with 0.7 exponent difference does not manage to ruin the logarithmic-ness of the entire series to any significant degree (by virtue of it being a small minority within the entire series).

For **high growth rates** such as say 60%, an integral power of ten is passed by the series relatively quickly, as in 1.00, $1.60^1$, $1.60^2$, $1.60^3$, $1.60^4$, $1.60^5$, namely as in 1.00, 1.60, 2.56, 4.10, 6.55, 10.49, so that in just six periods the series expands by one more order of magnitude, and the first digits manage to complete one full cycle. Therefore if one considers say 100 periods for the exponential 60% growth rate, there is no need really to carefully calibrate the first and the last elements in such a way so as to ensure that exponent difference is nearly an integral value. Here there is almost nothing to gain by such calibration; the series is anyhow very close to Benford after say 100 or so elements; and order of magnitude is already very large. Let us learn a lesson or two from two manifestations of such high 60% growth series:

For an exponential 60% growth series with base B = 5 being the quantity at time 0, and considering the **first 96 elements** of the series, we get:

{5.00, 8.00, 12.80, 20.48, 32.77, … , $1.88*10^{19}$, $3.01*10^{19}$, $4.81*10^{19}$, $7.70*10^{19}$, $1.23*10^{20}$}

and exponent difference is not an integral value; $LOG(1.23*10^{20}) - LOG(5.00) = 19.391$; yet 1st significant digits are nearly Benford, coming at {30.2, 17.7, 12.5, 8.3, 8.3, 8.3, 4.2, 6.3, 4.2}, with SSD = 8.8.

Almost nothing is gained now by adjusting the number of periods so as to obtain a near integral exponent difference. To illustrate this point we shall add 3 more elements, thus obtaining nearly an integral value for the exponent difference, and it shall be observed that first digits have barely nudged from their already near perfect logarithmic configuration.

For an exponential 60% growth series with base B = 5 being the quantity at time 0, and considering the **first 99 elements** of the series, we get:

{5.00, 8.00, 12.80, 20.48, 32.77, … , $7.70*10^{19}$, $1.23*10^{20}$, $1.97*10^{20}$, $3.15*10^{20}$, $5.04*10^{20}$}

and exponent difference now is nearly an integral value; $LOG(5.04*10^{20}) - LOG(5.00) = 20.004$; while 1st digits have barely changed, coming at {30.3, 17.2, 13.1, 8.1, 9.1, 8.1, 4.0, 6.1, 4.0}. In fact, SSD = 10.8 which is actually slightly higher, indicating a slightly less perfect fit to Benford!

For this reason the Fibonacci series is always very nearly perfectly Benford so long as one considers sufficient number of elements. Even though the series is defined in terms of additions, it approaches approximately a repeated multiplication process (of the exponential growth series type) very early on, with the golden ratio 1.61803399 as the F factor. Due to its relatively high growth rate of 61.8%, the Fibonacci series passes through an integral power of ten number each 5 terms approximately, thus the requirement of an integral exponent difference can be easily waived and one need not worry at all where we start and where we end, only that enough elements are considered. As a check: $1.618^5 \approx 11 > 10$, so that in 5 terms the cumulative increase is over ten-fold, and an integral power of ten number is passed.

For **low growth rates** such as say 2%, the consideration of merely 100 or 500 elements for example necessitates a very careful selection and calibration of the end points so that exponent difference is nearly of an integral value.

For an exponential 2% growth series with base B = 7 being the quantity at time 0, and carefully and deliberately considering the **first 118 elements** of the series, we get:

{7.00, 7.14, 7.28, 7.43, 7.58, 7.73, … , 65.60, 66.91, 68.25, 69.62, 71.01}

and exponent difference is nearly an integral value; LOG(71.01) – LOG (7.00) = 1.0062; therefore first digits are nearly Benford coming at {29.7, 16.9, 12.7, 9.3, 7.6, 6.8, 6.8, 5.1, 5.1}. SSD value is low 2.1.

On the other hand, mindlessly considering the **first 174 elements** of 2% exponential growth series with base B = 7 being the quantity at time 0, we get:

{7.00, 7.14, 7.28, 7.43, 7.58, 7.73, … , 198.85, 202.83, 206.88, 211.02, 215.24}

and exponent difference is decisively non-integral, coming at LOG(215.24) – LOG(7.00) = 1.4878; therefore even though there are 56 more elements now for this longer series, yet first digits here strongly deviate from Benford coming at {40.2, 13.8, 8.6, 6.3, 5.2, 4.6, 8.0, 6.3, 6.9}. SSD is much higher now, coming at 167.3.

**[18]  The Connection of Exponential Growth Series to the k/x Distribution**

Kossovsky (2014) discusses two essential results regarding exponential growth series:

<u>Proposition V</u>: Deterministic multiplication process of the exponential growth type is characterized by related log density being uniformly distributed (albeit discretely, not continuously).

<u>Proposition VI</u>: Deterministic multiplication process of the exponential growth type is of the k/x distribution albeit in a discrete fashion, having consistent and steady logarithmic behavior all along its entire range. The logarithmic requirement of an integral span on the log-axis (an integral exponent difference) applies to exponential growth series just as it does for any k/x distribution. The approximate [continuous] 'density' of the [discrete] exponential growth series is seen as falling steadily on the right at the constant logarithmic rate all throughout its entire range mimicking the k/x distribution.

Let us vividly demonstrate with an actual numerical example the intimate relationship between exponential growth series and k/x distribution. We choose monthly discrete readings on population that is growing continuously at 5% per year, so that $(F_{MONTHLY})^{12} = F_{YEARLY} = 1.05$, therefore $F_{MONTHLY}$ = 12th root of 1.05 = 1.004074, namely the very low 0.407% monthly growth. The initial quantity at time 0 is 1, and exactly 567 months (47.25 years) are considered (including the first month before growth begins), after which the quantity has grown to 9.987 (slightly shy of 10). The series is:

| | | | | | | | | | | | |
|---|---|---|---|---|---|---|---|---|---|---|---|
| 1.000 | 1.004 | 1.008 | 1.012 | 1.016 | 1.021 | 1.025 | 1.029 | 1.033 | 1.037 | 1.041 | 1.046 |
| 1.050 | 1.054 | 1.059 | 1.063 | 1.067 | 1.072 | 1.076 | 1.080 | 1.085 | 1.089 | 1.094 | 1.098 |
| 1.103 | 1.107 | 1.112 | 1.116 | 1.121 | 1.125 | 1.130 | 1.134 | 1.139 | 1.144 | 1.148 | 1.153 |
| 1.158 | 1.162 | 1.167 | 1.172 | 1.177 | 1.181 | 1.186 | 1.191 | 1.196 | 1.201 | 1.206 | 1.211 |
| 1.216 | 1.220 | 1.225 | 1.230 | 1.235 | 1.240 | 1.246 | 1.251 | 1.256 | 1.261 | 1.266 | 1.271 |
| 1.276 | 1.281 | 1.287 | 1.292 | 1.297 | 1.302 | 1.308 | 1.313 | 1.318 | 1.324 | 1.329 | 1.335 |
| 1.340 | 1.346 | 1.351 | 1.357 | 1.362 | 1.368 | 1.373 | 1.379 | 1.384 | 1.390 | 1.396 | 1.401 |
| 1.407 | 1.413 | 1.419 | 1.424 | 1.430 | 1.436 | 1.442 | 1.448 | 1.454 | 1.460 | 1.465 | 1.471 |
| 1.477 | 1.483 | 1.490 | 1.496 | 1.502 | 1.508 | 1.514 | 1.520 | 1.526 | 1.533 | 1.539 | 1.545 |
| 1.551 | 1.558 | 1.564 | 1.570 | 1.577 | 1.583 | 1.590 | 1.596 | 1.603 | 1.609 | 1.616 | 1.622 |
| 1.629 | 1.636 | 1.642 | 1.649 | 1.656 | 1.662 | 1.669 | 1.676 | 1.683 | 1.690 | 1.696 | 1.703 |
| 1.710 | 1.717 | 1.724 | 1.731 | 1.738 | 1.745 | 1.753 | 1.760 | 1.767 | 1.774 | 1.781 | 1.789 |
| 1.796 | 1.803 | 1.811 | 1.818 | 1.825 | 1.833 | 1.840 | 1.848 | 1.855 | 1.863 | 1.870 | 1.878 |
| 1.886 | 1.893 | 1.901 | 1.909 | 1.917 | 1.924 | 1.932 | 1.940 | 1.948 | 1.956 | 1.964 | 1.972 |
| 1.980 | 1.988 | 1.996 | 2.004 | 2.012 | 2.021 | 2.029 | 2.037 | 2.045 | 2.054 | 2.062 | 2.070 |
| 2.079 | 2.087 | 2.096 | 2.104 | 2.113 | 2.122 | 2.130 | 2.139 | 2.148 | 2.156 | 2.165 | 2.174 |
| 2.183 | 2.192 | 2.201 | 2.210 | 2.219 | 2.228 | 2.237 | 2.246 | 2.255 | 2.264 | 2.273 | 2.283 |
| 2.292 | 2.301 | 2.311 | 2.320 | 2.330 | 2.339 | 2.349 | 2.358 | 2.368 | 2.377 | 2.387 | 2.397 |
| 2.407 | 2.416 | 2.426 | 2.436 | 2.446 | 2.456 | 2.466 | 2.476 | 2.486 | 2.496 | 2.506 | 2.517 |
| 2.527 | 2.537 | 2.548 | 2.558 | 2.568 | 2.579 | 2.589 | 2.600 | 2.610 | 2.621 | 2.632 | 2.643 |
| 2.653 | 2.664 | 2.675 | 2.686 | 2.697 | 2.708 | 2.719 | 2.730 | 2.741 | 2.752 | 2.763 | 2.775 |
| 2.786 | 2.797 | 2.809 | 2.820 | 2.832 | 2.843 | 2.855 | 2.866 | 2.878 | 2.890 | 2.902 | 2.913 |
| 2.925 | 2.937 | 2.949 | 2.961 | 2.973 | 2.985 | 2.998 | 3.010 | 3.022 | 3.034 | 3.047 | 3.059 |
| 3.072 | 3.084 | 3.097 | 3.109 | 3.122 | 3.135 | 3.147 | 3.160 | 3.173 | 3.186 | 3.199 | 3.212 |
| 3.225 | 3.238 | 3.251 | 3.265 | 3.278 | 3.291 | 3.305 | 3.318 | 3.332 | 3.345 | 3.359 | 3.373 |
| 3.386 | 3.400 | 3.414 | 3.428 | 3.442 | 3.456 | 3.470 | 3.484 | 3.498 | 3.513 | 3.527 | 3.541 |
| 3.556 | 3.570 | 3.585 | 3.599 | 3.614 | 3.629 | 3.643 | 3.658 | 3.673 | 3.688 | 3.703 | 3.718 |
| 3.733 | 3.749 | 3.764 | 3.779 | 3.795 | 3.810 | 3.826 | 3.841 | 3.857 | 3.873 | 3.888 | 3.904 |
| 3.920 | 3.936 | 3.952 | 3.968 | 3.984 | 4.001 | 4.017 | 4.033 | 4.050 | 4.066 | 4.083 | 4.099 |
| 4.116 | 4.133 | 4.150 | 4.167 | 4.184 | 4.201 | 4.218 | 4.235 | 4.252 | 4.270 | 4.287 | 4.304 |
| 4.322 | 4.340 | 4.357 | 4.375 | 4.393 | 4.411 | 4.429 | 4.447 | 4.465 | 4.483 | 4.501 | 4.520 |
| 4.538 | 4.557 | 4.575 | 4.594 | 4.612 | 4.631 | 4.650 | 4.669 | 4.688 | 4.707 | 4.726 | 4.746 |
| 4.765 | 4.784 | 4.804 | 4.823 | 4.843 | 4.863 | 4.883 | 4.903 | 4.922 | 4.943 | 4.963 | 4.983 |
| 5.003 | 5.024 | 5.044 | 5.065 | 5.085 | 5.106 | 5.127 | 5.148 | 5.169 | 5.190 | 5.211 | 5.232 |
| 5.253 | 5.275 | 5.296 | 5.318 | 5.339 | 5.361 | 5.383 | 5.405 | 5.427 | 5.449 | 5.471 | 5.494 |
| 5.516 | 5.538 | 5.561 | 5.584 | 5.606 | 5.629 | 5.652 | 5.675 | 5.698 | 5.722 | 5.745 | 5.768 |
| 5.792 | 5.815 | 5.839 | 5.863 | 5.887 | 5.911 | 5.935 | 5.959 | 5.983 | 6.008 | 6.032 | 6.057 |
| 6.081 | 6.106 | 6.131 | 6.156 | 6.181 | 6.206 | 6.232 | 6.257 | 6.282 | 6.308 | 6.334 | 6.360 |
| 6.385 | 6.411 | 6.438 | 6.464 | 6.490 | 6.517 | 6.543 | 6.570 | 6.597 | 6.623 | 6.650 | 6.678 |
| 6.705 | 6.732 | 6.759 | 6.787 | 6.815 | 6.842 | 6.870 | 6.898 | 6.926 | 6.955 | 6.983 | 7.011 |
| 7.040 | 7.069 | 7.097 | 7.126 | 7.155 | 7.185 | 7.214 | 7.243 | 7.273 | 7.302 | 7.332 | 7.362 |
| 7.392 | 7.422 | 7.452 | 7.483 | 7.513 | 7.544 | 7.575 | 7.605 | 7.636 | 7.667 | 7.699 | 7.730 |
| 7.762 | 7.793 | 7.825 | 7.857 | 7.889 | 7.921 | 7.953 | 7.986 | 8.018 | 8.051 | 8.084 | 8.117 |
| 8.150 | 8.183 | 8.216 | 8.250 | 8.283 | 8.317 | 8.351 | 8.385 | 8.419 | 8.453 | 8.488 | 8.522 |
| 8.557 | 8.592 | 8.627 | 8.662 | 8.697 | 8.733 | 8.768 | 8.804 | 8.840 | 8.876 | 8.912 | 8.949 |
| 8.985 | 9.022 | 9.058 | 9.095 | 9.132 | 9.170 | 9.207 | 9.244 | 9.282 | 9.320 | 9.358 | 9.396 |
| 9.434 | 9.473 | 9.511 | 9.550 | 9.589 | 9.628 | 9.667 | 9.707 | 9.746 | 9.786 | 9.826 | 9.866 |
| 9.906 | 9.946 | 9.987 | | | | | | | | | |

The confluence of such meticulous calibration of the first and the last terms where exponent difference is very nearly an integral value [integer 1], together with the abundance of points, ensures an exceptional and very rare fit to Benford, with SSD = 0.03. First digits come out as:

Exponential 0.407% Growth - {30.16, 17.64, 12.35, 9.70, 7.94, 6.70, 5.82, 5.11, 4.59}
Benford's Law First digits − {30.10, 17.61, 12.49, 9.69, 7.92, 6.69, 5.80, 5.12, 4.58}

Let us now construct a histogram of all the quantities falling on (1, 10). Such a construction would omit the time dimension altogether of course. We construct exactly 27 bins for the entire width of (10 – 1) = 9 units, resulting in 9/27 = 0.333 width for each bin. Hence the first bin is on [1.000, 1.333) and there we encounter 71 data points. The second bin is on [1.333, 1.666) and there we encounter 55 data points. For the very last bin on [9.666, 10.000) we encounter only 9 points.

The series of the count of the data points falling within each bin in the entire range is: {71, 55, 45, 38, 33, 29, 26, 23, 21, 20, 18, 17, 16, 15, 14, 13, 13, 12, 12, 10, 11, 10, 10, 9, 9, 8, 9}. Figure 31 depicts the histogram.

The series of midpoints on the x-axis for each bin are: {1.167, 1.500, 1.833, 2.167, 2.500, 2.833, 3.167, 3.500, 3.833, 4.167, 4.500, 4.833, 5.167, 5.500, 5.833, 6.167, 6.500, 6.833, 7.167, 7.500, 7.833, 8.167, 8.500, 8.833, 9.167, 9.500, 9.833}.

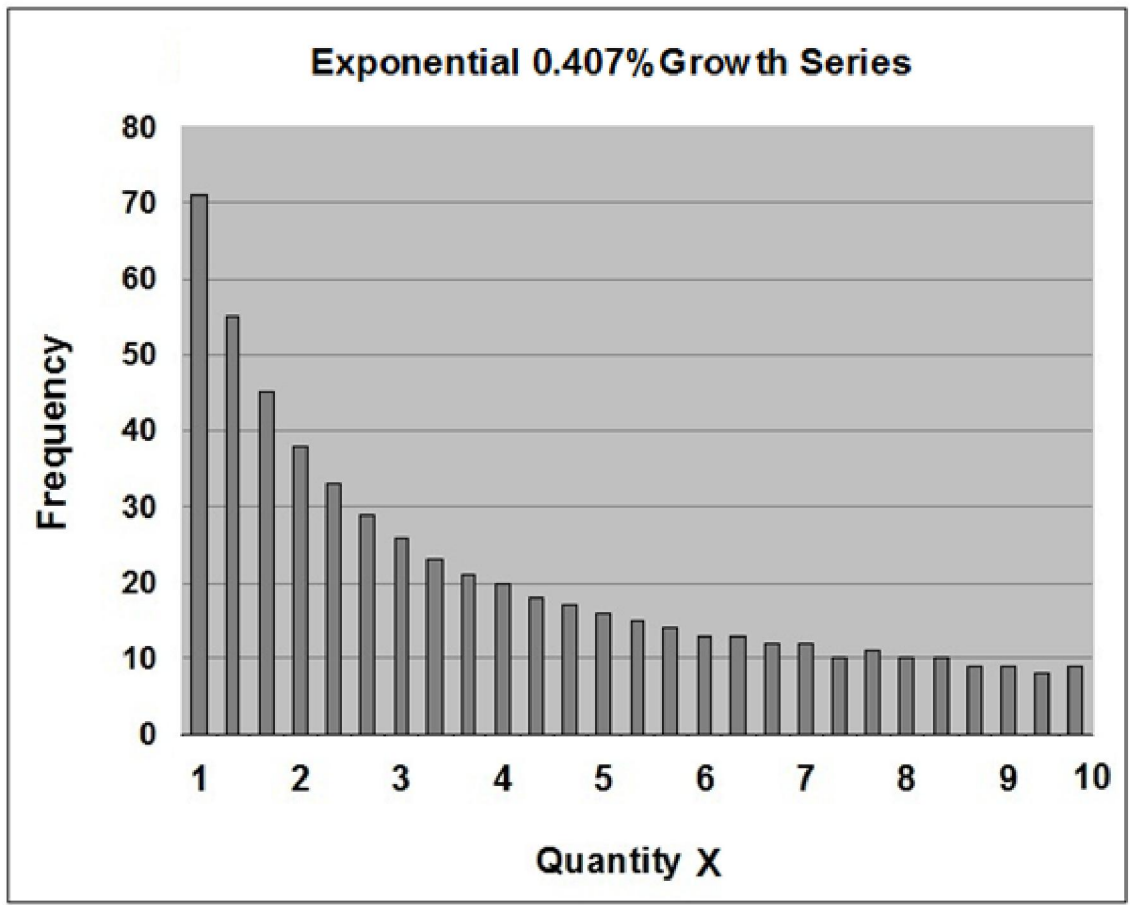

**Figure 31:** Histogram of Discrete Exponential 0.407% Growth from 1 to 10

Let us first roughly check compatibility with k/x distribution. To do this we shall use an essential property of the k/x distribution which states that by doubling quantity x we reduce its frequency (density height) by half. Surely for transition such as: x → 2x the density of k/x always diminishes by half as in: k/x → k/(2x) = (k/x)/2. Very roughly, the two data points of (1.167, 71) and (2.500, 33) on the histogram of the exponential growth series represent approximately the doubling of quantity x, which induces 71 − 33 = 38 reduction in the frequency, namely 38/71 or 53.5% reduction – being comparable approximately to the theoretical 50%. Very roughly, the two data points of (2.167, 38) and (4.500, 18) represent approximately the doubling of quantity x, which induces 38 − 18 = 20 reduction in the frequency, namely 20/38 or 52.6% reduction – being comparable approximately to theoretical 50%. Very roughly, the data points of (4.167, 20) and (8.500, 10) represent the doubling of quantity x, which induces 20 − 10 = 10 reduction in the frequency, namely 10/20 reduction – corresponding nicely to the theoretical 50% expectation.

We shall now look for some k/x distribution that could best fit the histogram of Figure 31. Let us find the most fitting parameter k value, one that would minimize the sum of squared errors (SSE) between such an ideal k/x distribution and the histogram of the 0.407% monthly growth series of Figure 31 [using the midpoints]. We shall set the derivative of SSE with respect to k equal to zero and solve for k. We therefore obtain:

$$\text{SSE} = (k/1.167 − 71)^2 + (k/1.500 − 55)^2 + (k/1.833 − 45)^2 + \ldots + (k/9.833 − 9)^2$$

$$\text{d(SSE)/d(k)} = 2(k/1.167 − 71)/1.167 + 2(k/1.500 − 55)/1.500 + 2(k/1.833 − 45)/1.833 + \ldots + 2(k/9.833 − 9)/9.833$$

$$0 = 2(k/1.167 − 71)/1.167 + 2(k/1.500 − 55)/1.500 + 2(k/1.833 − 45)/1.833 + \ldots + 2(k/9.833 − 9)/9.833$$

Dividing both sides of the equation by 2 we get:

$$0 = (k/1.167 − 71)/1.167 + (k/1.500 − 55)/1.500 + (k/1.833 − 45)/1.833 + \ldots + (k/9.833 − 9)/9.833$$

Rearranging terms involving k and those not involving k, we get:

$$0 = (k/1.167)/1.167 + (k/1.500)/1.500 + (k/1.833)/1.833 + \ldots + (k/9.833)/9.833$$
$$- (71)/1.167 - (55)/1.500 - (45)/1.833 - \ldots - (9)/9.833$$

$$0 = k*[(1/1.167)/1.167 + (1/1.500)/1.500 + (1/1.833)/1.833 + \ldots + (1/9.833)/9.833]$$
$$- (71)/1.167 - (55)/1.500 - (45)/1.833 - \ldots - (9)/9.833$$

$$0 = k*[ 2.67 ] - 220.28$$

$$k = (220.28)/(2.67) = \mathbf{82.4}$$

Figure 32 depicts the excellent fit of 82.4/x distribution with the discrete exponential 0.407% growth series. Points on the line of 82.4/x were constructed as in (x, 82.4/x) using the series of 27 midpoints given earlier for the 27 x values.

**82.4 /** 1.167 = 70.6     and 71 for the exponential growth series
**82.4 /** 1.500 = 54.9     and 55 for the exponential growth series
**82.4 /** 1.833 = 44.9     and 45 for the exponential growth series
**82.4 /** 2.167 = 38.0     and 38 for the exponential growth series

   … etc. ….

**82.4 /** 9.833 = 8.4     and 9 for the exponential growth series

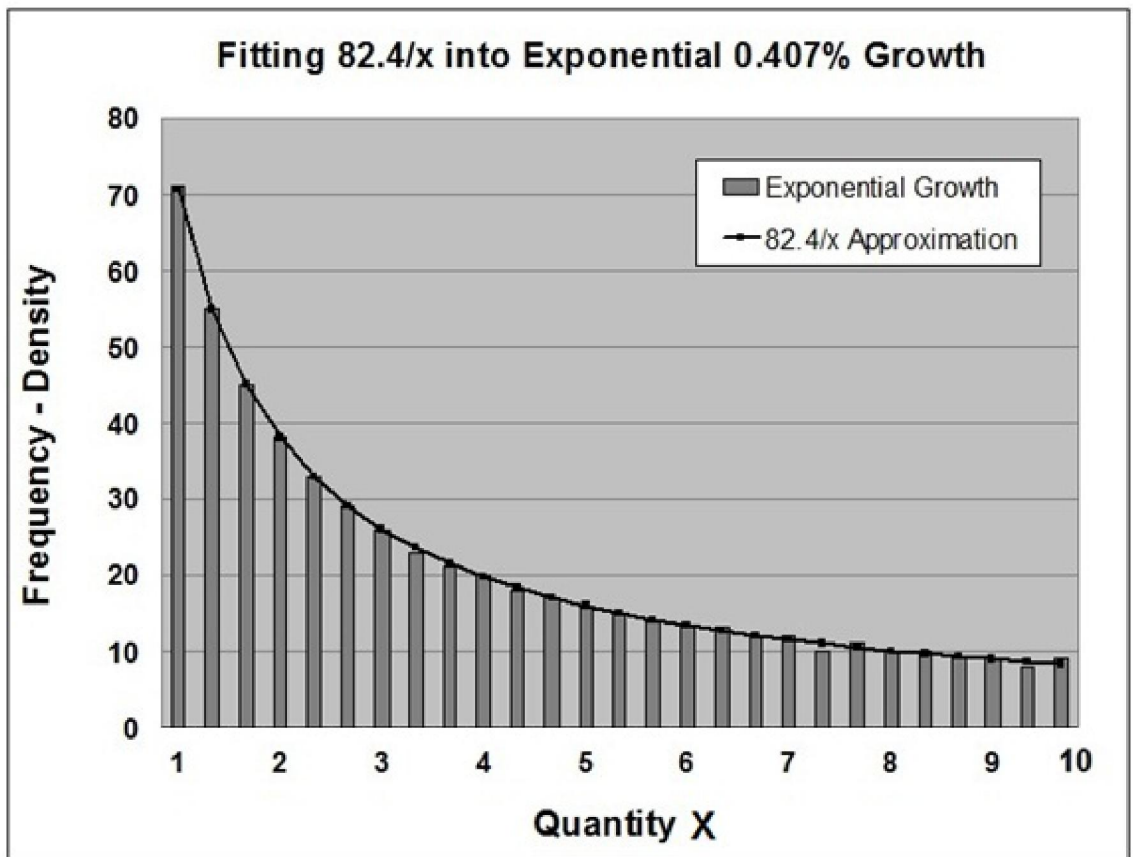

**Figure 32:** The Excellent Fit of 82.4**/**x Distribution into Exponential 0.407% Growth

One straightforward and decisive mathematical confirmation of the intimate relationship between exponential growth series and the k/x distribution in the context of Benford's Law can be given by the comparison of the expression $\mathbf{X(T) = BF^T}$ of the continuous form of exponential growth series and $\mathbf{X(Y) = 10^Y}$, Y being the continuous Uniform Distribution, which according to Proposition II in Kossovsky (2014) is of the k/x form. Proposition II and Corollary II state that if Y is uniformly distributed over [R, S], with the length S - R being an integer, then $X = 10^Y$ over $[10^R, 10^S]$ is distributed as in f(x) = (1/[(S-R)*ln(10)])*(1/x), and it is perfectly logarithmic. Indeed, a moment thought would convinced anybody that the very definition of Benfordness as uniformity of mantissa is also what is seen in the expression $X(Y) = 10^Y$, because whenever X varies between adjacent integral powers of ten, such as from 1 to 10 say, its mantissa Y varies and cycles fully and only once on its entire space of (0, 1).

Let us summarize:

The k/x distribution   $\rightarrow$   $\mathbf{X(Y) = 10^Y}$
Exponential growth   $\rightarrow$   $\mathbf{X(T) = BF^T}$

The striking similarity between these two expressions is found in the fact that both are written in an exponential form, with a fixed constant raised to a variable power of the Uniform type. The expression for k/x can be considered as in exponential growth with B = 1 and F = 10.

As discussed in chapter 15 regarding continuous growth, the continuous time T variable for all practical purposes can be discretely selected every minute or every second for the evaluation of the growing quantity as time progresses. Surely time selection must be made fairly, evenly, and uniformly, without focusing on certain time intervals at the expense of other time intervals. By simply choosing to evaluate the series consistently every minute say, fairness and uniformity are ensured by default. The uniform and even dispersions of the projections onto the horizontal time axis of the points in Figures 28, 29, and 30 visually and clearly demonstrate this principle. For example, it would be wrong to take more measurements of the growing quantity in January say each 1/2 a minute, take fewer measurements in February each minute, and take even fewer measurements in March only each 5 minutes.

If the exponent difference of the exponential series $X(T) = F^T$ for a particular T range is approximately of an integral value, then the series is certainly nearly Benford, and multiplying the series by any factor/base B to obtain $X(T) = BF^T$ does not significantly change its logarithmic property as predicated by the Scale Invariance Principle. This implies that the value of B = 1 in the k/x should not affect logarithmic behavior. In addition, surely the high F = 10 value in the k/x case should not preclude logarithmic behavior whatsoever; rather it only implies that growth is with the very high rate of 900% per period. In general, a change in base for the k/x distribution does not change its nature, but such a change in base implies that the end points of the Y range must adjust in order to obtain perfect Benford behavior. Hence $X(Y) = e^Y$ or $X(Y) = 8^Y$ are also Benford, given corresponding proper spans of Y where the ranges lead to perfect logarithmic behavior. There is nothing special in having base 10 for the k/x distribution as in $X(Y) = 10^Y$.

# [19]  The Last Term of Random Exponential Growth Series is Benford

An innovative explanation of the prevalence of Benford's Law in the physical world has been made relatively recently by **Kenneth Ross** utilizing the collection of the last term of numerous exponential growth series having random growth rate and random initial base.

In Ross (2011) titled "Benford's Law, A Growth Industry", a rigorous mathematical proof  is presented showing that a collection of the last term of numerous exponential growth series $\{BF^N\}$ is logarithmic in the limit as N goes to infinity - given that both B and F are randomly selected from suitable Uniform distributions.

The generic expression $BF^N$ signifies exponential growth series, where base B is the initial quantity at time zero, F is the growth factor per period, and N is the number of periods.

Roughly speaking, Ross chooses B and F as in the continuous Uniforms on (1, 10], where all digits 1 to 9 enjoy equitable proportions on the interval, thus avoiding a priori any possible bias against any particular set of digits.

Ross assigns the same value of N (i.e. length) to all the growth series. In other words, Rose is considering the set of growth series, all with equal length of time.

Ross model can be summarized as: $\lim_{N \to \infty}$ $(Uniform_1)*(Uniform_2)^N$.

One real-life data set that could obvious be relevant to his model is population census data on all the cities and towns in a given country, and which is actually almost always perfectly Benford – except for very small countries with few population centers. Hence, his model can be thought of as a census population snapshot of a large collection of relatively old and well established towns and cities, all being established simultaneously in the same period, and with inter-migration strictly prohibited. In addition, Ross has provided a second model that allows for F to be randomly selected anew for each period, and this model constitutes a more realistic model of cities and towns where the growth rate is generally not constant, but rather varies with the period. Monte Carlo simulations [with a large number of cities] for his first model strongly corroborate his assertion, giving very satisfactory results even after only 20 growth periods. Monte Carlo simulations for his second model show even more rapid convergence, where excellent fit to Benford is obtained even after only 5 or 6 growth periods.

The first model can be succinctly expressed as:

**$Uniform_1[1, 10]*Uniform_2[1, 10]^N$**

The second model can be succinctly expressed as:

**$Uniform_1[1, 10]*[ Uniform_2[1, 10]*Uniform_3[1, 10]*$  … N times …  $*Uniform_{N+1}[1, 10] ]$**

Admittedly, these Uniform distributions are of low order of magnitude of merely 1, and in general multiplication processes with variables of low orders of magnitude are not conducive to Benford behavior. Nonetheless, since the expressions here contain only purely multiplicative terms without any disruptive additions, there exists no tug of war whatsoever between addition and multiplication, and in the limit as N gets large enough and numerous such multiplications are applied, resultant order of magnitude of the entire scheme becomes sufficiently high.

The enormous significance of such logarithmic result for the final elements of random exponential growth series is that it can serve as a logarithmic model for all entities that spring into being gradually via random growth, be it a set of cities and towns with growing populations, rivers forming and enlarging gradually along an incredibly slow geological time scales, biological cells growing, stars and planets formations, and so forth. The potential scope covered by Ross' results is enormous! Hence, instead of having to become a specialist in each and every scientific discipline, laboriously looking for components and physical measurements that may serve as multiplicands within any scientific data set that is Benford, or attempting to fit the physical setup into some data aggregation scheme, Ross lets us skip all this, and obtain extremely general results [for entities that grow] relevant to so many scientific disciplines.

Using the notation $B_J$ for one simulated realization from the random distribution of B, and $F_J$ for one realization from the random distribution of F, his <u>first</u> model can be expressed as the data set $\{B_1F_1^N, B_2F_2^N, B_3F_3^N, B_4F_4^N, \ldots, B_MF_M^N\}$ where N is sufficiently large for the desired convergence and accuracy, and M is the number of simulated cities or entities. Since each $B_J$ and each $F_J$ are realized from the Uniform on the interval $[1, 10)$ and can be expressed as $U_I$, this can be written more succinctly as: $\{U_1U_2^N, U_3U_4^N, U_5U_6^N, U_7U_8^N, \ldots, U_{2M-1}U_{2M}^N\}$.

For N fixed as 50 for example, his <u>second</u> model can be succinctly expressed as:
$\{U_1U_2U_3 \ldots U_{50}U_{51}, \quad U_{52}U_{53}U_{54}\ldots U_{101}U_{102}, \quad U_{103}U_{104}U_{105}\ldots U_{152}U_{153}, \quad \ldots \text{ M times } \ldots \}$

Little reflection is needed to realize that Ross second model be interpreted as repeated multiplications of random Uniform distributions, and which is distributed almost exactly as Lognormal with high shape parameter and high order of magnitude [hence Benford] as predicated by the Multiplicative Central Limit Theorem [assuming the multiplications of numerous such Uniforns]. His first model can also be interpreted in terms of the well-known general result in the field, namely that any random data set X (logarithmic or non-logarithmic) transformed by raising each value within the data to the Nth power, converges to the logarithmic as N gets large. That is, the set $\{X_1, X_2, X_3, \ldots, X_M\}$ of M realizations from any variable X transformed into $\{X_1^N, X_2^N, X_3^N, \ldots, X_M^N\}$ is logarithmic as N gets large.

In Ross models, both B and F are random variables, while the value of N is identical and constant for all cities. Surprisingly, inverted models where both B and F are fixed as constants being identical for all cities, while only N varies randomly by city, also converge to Benford! In one particular Monte Carlo computer simulation, the model is of a country with 242 cities, all being randomly established at different years - ranging uniformly from 1 to 300 years ago; all having an identical growth rate of 11%; all starting from population of 1 (i.e. a single person). Here we allow for fractional values of persons, since the model represents the generic case of quantitative growth, not necessarily only of integral values. In summarizing the model, the country is said to have cities of varying (random) age; some are very old and established cities, some are not as old, and some are more recent modern cities. This model can be succinctly expressed as $(1)*1.11^{\text{UNIFORM }\{1,2,3,\ldots,300\}}$. The variable under consideration is the current snapshot of the populations of all existing 242 cities and towns (only the last elements of the growth series, not including historical population records). Six different simulation runs as well as their average digit distribution came out as follows:

Fixed B & F, Varied N – {27.7, 19.4,  9.9, 10.7, 9.9, 6.6, 7.4, 5.4, 2.9}
Fixed B & F, Varied N – {34.3, 18.6,  9.1,  9.5, 6.6, 5.4, 7.4, 4.5, 4.5}
Fixed B & F, Varied N – {28.5, 24.4, 10.3,  6.6, 7.0, 4.5, 6.2, 5.0, 7.4}
Fixed B & F, Varied N – {28.1, 16.5, 16.1,  7.9, 6.6, 9.1, 6.6, 5.0, 4.1}
Fixed B & F, Varied N – {29.3, 21.5, 12.8,  7.4, 7.4, 7.9, 5.4, 4.1, 4.1}
Fixed B & F, Varied N – {29.3, 17.4, 13.2, 12.4, 4.5, 5.8, 6.6, 7.0, 3.7}

**AVG - fixed BF varied N – {29.5, 19.6, 11.9,  9.1, 7.0, 6.5, 6.6, 5.2, 4.5}**
**Benford's Law 1st digits – {30.1, 17.6, 12.5,  9.7, 7.9, 6.7, 5.8, 5.1, 4.6}**

SSD value for this average digit distribution is 6.6, indicating strong conformity to Benford.

In another Monte Carlo computer simulation scheme, the model is of a country with 242 cities; all being established at the same time 1700 years ago (i.e all are of the same age of 1700); all starting from population of 1 (i.e. a single person); all having a fixed but distinct growth rate between 0% and 11% as F is chosen from the continuous Uniform(1.00, 1.11). The variable under consideration is the current snapshot of the populations of all existing 242 cities and towns after 1700 years (the last elements of the growth series, not including historical population records). In summarizing the model, the country is said to have cities of identical age N as well as identical initial population base B, but with distinct random F growth factors. This model can be succinctly expressed as: $(1)*\text{Uniform}(1.00, 1.11)^{1700}$. Seven different simulation runs as well as their average digit distribution came out as follows:

Fixed B & N, Varied F  – {36.8, 14.9, 10.3, 11.2, 5.0,  7.0, 7.4, 4.1, 3.3}
Fixed B & N, Varied F  – {30.6, 16.9, 16.5,  9.5, 8.7, 2.9, 6.2, 5.0, 3.7}
Fixed B & N, Varied F  – {31.8, 20.7, 12.8,  7.9, 6.6, 6.6, 7.4, 3.7, 2.5}
Fixed B & N, Varied F  – {31.0, 18.6,  9.5,  9.9, 10.3, 7.9, 4.1, 5.4, 3.3}
Fixed B & N, Varied F  – {30.6, 15.3, 13.6,  7.0,  9.9, 6.2, 6.6, 7.4, 3.3}
Fixed B & N, Varied F  – {33.5, 13.6,  9.5, 11.2,  9.1, 9.5, 4.1, 5.4, 4.1}
Fixed B & N, Varied F  – {27.7, 20.7, 13.2,  9.5,  6.2, 6.2, 7.0, 5.4, 4.1}

**AVG - fixed BN varied F  – {31.7, 17.2, 12.2, 9.4, 8.0, 6.6, 6.1, 5.2, 3.5}**
**Benford's Law First digits – {30.1, 17.6, 12.5, 9.7, 7.9, 6.7, 5.8, 5.1, 4.6}**

SSD value for this average digit distribution is 4.2, indicating strong conformity to Benford.

Falling well below 1700 years (periods) yields worsening digital results. For example, one run of the above model with only 300 periods yields {34.7, 16.1, 12.4, 8.7, 7.0, 7.4, 4.5, 4.1, 5.0} and its higher SSD value of 28.6 is an indication of diminished conformity to Benford. Yet even this digital configuration is quite similar to the logarithmic overall. A comparison between this scheme and Ross first model points to the fact that fixing base B as a constant (i.e. reducing randomness) requires the system to have many more periods of growth before a convergence to Benford is achieved. Ross first model achieves near-Benford digit configuration quickly even after 20 or so periods, while this model with only factor F varying randomly requires at least 1000 to 2000 periods for convergence.

In another Monte Carlo simulation scheme with all three variables {B, F, N} chosen randomly, results are nearly Benford even with the shorter time span of 100 years/periods. This is so since there is 'more randomness' in the system, or rather because this model is also covered under Ross first model. The model is of a country with 242 cities; all being randomly established at different years ranging from 1 to 100 and selected from the discrete Uniform {1, 2, 3, …, 100}; all having random growth rate between 0% and 14% as F is chosen from the continuous Uniform(1.00, 1.14); all starting from a random population base chosen from the continuous Uniform[1, 10). In summarizing the model, the country is said to have cities of varying random age, random growth rate, and random initial population count. This model can be succinctly expressed as: **Uniform[1, 10)\*Uniform(1.00, 1.14)$^{\text{UNIFORM }\{1,2,3, …, 100\}}$**. The variable under consideration is the current snapshot of the populations of all 242 existing cities. Six different simulation runs and their average digit distribution came out as follows:

Varied B, N, F  –  {33.9, 18.2,  9.5,  7.9,  8.7, 7.0, 6.2, 4.1, 4.5}
Varied B, N, F  –  {29.8, 15.7,  7.9,  9.1,  9.9, 9.1, 5.8, 7.4, 5.4}
Varied B, N, F  –  {31.4, 16.5, 12.0,  9.9,  9.1, 7.0, 4.5, 5.0, 4.5}
Varied B, N, F  –  {27.7, 16.1, 16.1,  8.7,  7.4, 7.0, 6.6, 4.1, 6.2}
Varied B, N, F  –  {26.0, 17.8, 16.1, 11.2,  5.4, 8.7, 4.1, 6.2, 4.5}
Varied B, N, F  –  {28.1, 14.0, 10.7, 12.8, 10.3, 8.3, 5.8, 5.8, 4.1}

**AVG Varied B, N, F –  {29.5, 16.4, 12.1,  9.9, 8.5, 7.9, 5.5, 5.4, 4.9}**
**Benford's L. 1st digits – {30.1, 17.6, 12.5, 9.7, 7.9,  6.7, 5.8, 5.1, 4.6}**

SSD value for this average digit distribution is 4.1, indicating strong conformity to Benford.

Surprisingly, Ross first model where only F and B vary, while N is fixed and identical for all (i.e. cities begin their existence simultaneously at time zero and all are of the same age) yields better results for short N time span such as 20 periods, 100 periods, and such – as compared with this 'more random' model above. The obvious explanation for this apparent paradox is that for a very short time span all the series should start early on [as in Ross model], so that they would have sufficient lengths to even begin resembling growth series. By varying N randomly [as in this last model] we unfortunately allow into the system some very short and immature 'series' of say 3 or 7 years, 'series' which hamper convergence. Empirically, for very long N time span, both models yield the logarithmic equally, and this fact is consistent with this explanation. In addition, another reasonable explanation is that here we restrict the range F to the Uniform(1.00, 1.14) which has a much lower order of magnitude than the order of magnitude of the Uniform(1, 10] for Ross first model.

In contrast to all the successful models above, the model where F and N are fixed as constants while only base B varies randomly, does <u>not</u> converge to Benford in any way. In one Monte Carlo computer simulation scheme, the model is of a country with 242 cities; all being established simultaneously 2500 years ago (i.e. all are of the same ancient age); all having the same 5% growth rate; all starting from a varied population base B randomly selected from the continuous Uniform(0, 100). This 'less random' model can be succinctly expressed as: **Uniform(0, 100)*1.05$^{2500}$**. The variable under consideration is the current snapshot of the populations of all 242 existing cities after 2500 years. First digit distribution here came out as {10.3, 9.1, 9.5, 14.0, 11.6, 11.2, 14.0, 12.8, 7.4}. Surely if F and N are fixed and are equal for all the cities, then the variation in the final population values is due only to the variable base B, being a random multiplicative factor, hence the Uniform(0, 100)*1.05$^{2500}$ model is equivalent to the model Uniform(0, 100)*Constant, which is of course not Benford at all, but rather of the symmetrical Uniform distribution configuration.

The essential feature leading to Benford convergence in all the above successful models is the multiplicative form in how population [or the generic quantity under consideration] is increasing. If one views [multiplicative] exponential growth as repeated additions, then one sees larger and larger quantities being added as the years pass. A constant 5% population increase implies an addition of 5 at the year when the population is at 100, an addition of 50 at the year when the population is at 1000, and an addition of 500 at the year when the population is at 10000. On a profound level this is almost always the case in the natural world. The larger the city, the more it attracts new inhabitants. A small village attracts very few newcomers. The larger the forming star, the more mass it attracts gravitationally and adds to itself. In sharp contrast, an additive type of growth is one where the quantity being added each year is fixed, not being proportional to the present size of the population. Let us consider the model of a country with 242 cities; all being randomly established at different years ranging uniformly from 1 to 100; all starting from a variable population base B, randomly selected from the continuous Uniform(0, 50); all experiencing fixed and identical additions of quantity D per period. This model can be succinctly expressed as: **Uniform(1, 50) + D*Uniform{1,2,3, ..., 100}**. The variable under consideration is the current snapshot of the populations of all 242 existing cities and towns. Three different simulation runs with distinct D values are shown:

Varied B & N, fixed D = 20 Additive − {47.5, 6.2, 8.7, 7.9, 7.4, 5.4, 2.9, 8.7, 5.4}
Varied B & N, fixed D = 33 Additive − {38.8, 33.5, 12.4, 1.7, 2.5, 2.1, 3.7, 2.1, 3.3}
Varied B & N, fixed D = 8 Additive − {17.4, 10.3, 11.2, 13.2, 14.5, 14.0, 11.6, 6.2, 1.7}

Certainly this model is not Benford at all, being a linear combination of sorts of two Uniform distributions.

**Acknowledgement:**

The late mathematician **Ralph Raimi** who has eloquently written some of the first mathematically rigorous articles on Benford's Law in the 1960's and 1970's has also written about the existence of such anomalous series, using the term 'reentrant series' - in an article titled "The Peculiar Distribution of First Digit", Scientific America, 1969. The author is still quite content to re-invent this old wheel by himself and to add several new features to it. The author wishes to convey his strong sense of affinity and rapport with Raimi for thinking along the same lines in the quest for a thorough understanding of the behavior of exponential growth series in the context of Benford's Law. A lasting impression for the author was a 2-day visit to Raimi's residence in Rochester, New York, in April 2013, discussing Benford's Law in general, and in particularly the discovery of Digital Development Pattern which Raimi enthusiastically endorsed. Funnily, anomalous (reentrant) series were never discussed during the entire visit!? Raimi passed away on January 2, 2017, at the age of 92.

The author wishes to thank the mathematician **Kenneth Ross** for his correspondence regarding the two logarithmic growth models of his article.

PATENT:

The U.S. Patent Office  # 9,058,285.  Inventor: Alex Ely Kossovsky.
Date Granted: June 16, 2015.   http://www.google.com/patents/US20140006468
Titled: "Method and system for Forensic Data Analysis in fraud detection employing a digital pattern more prevalent than Benford s Law".

Alex Ely Kossovsky, Jan 6, 2019
akossovsky@gmail.com